\documentclass[final]{amsart}   
\usepackage[all]{xy}

\def\CC{\mathbb{C}}
\def\RR{\mathbb{R}}

\def\NN{\mathbb{N}}

\def\dbar{{d\!\!\!\;{\vrule height 5.8 pt depth -5.4 pt width 3 pt}}}

\DeclareMathOperator{\Tr}{Trace}

\DeclareMathOperator{\End}{End}
\DeclareMathOperator{\Hom}{Hom}

\DeclareMathOperator{\res}{res}

\DeclareMathAlphabet{\mathbi}{OT1}{cmr}{bx}{it}

\theoremstyle{plain}
\newtheorem{thm}{Theorem}[section]
\newtheorem{lemma}[thm]{Lemma}
\newtheorem{prop}[thm]{Proposition}

\theoremstyle{definition}
\newtheorem{defn}{Definition}[section]

\theoremstyle{remark}
\newtheorem{remark}{Remark}[section]
\newtheorem{example}{Example}

\begin{document}
\title[Pseudodifferential Operators on Orbifolds]
{
Seeley's Theory of Pseudodifferential Operators on Orbifolds
}
\author{Bogdan Bucicovschi}
\address{
Bogdan Bucicovschi\\
Mathematics Department\\
Ohio State University\\
231 W 18th Avenue\\
Columbus, OH 43210
}
\email{bogdanb@math.ohio-state.edu}
\date{\today}

\begin{abstract}
In this paper we extend the results of
Robert Seeley concerning the complex powers of elliptic pseudodifferential 
operators and the residues of their zeta function to operators acting on
vector orbibundles over orbifolds.

We present the theory of pseudodifferential operators acting on a 
vector orbibundle over an orbifold, construct the zeta function of 
an elliptic pseudodifferential operator and show the existence of a 
meromorphic extension to $\CC$ with at most simple poles. We give 
formulas for generalized densities on the orbifold whose
integrals compute the residues of the zeta function.
\end{abstract}

\maketitle
\setcounter{section}{-1}
\tableofcontents

\vfill\eject

\section{Introduction}

In this paper we extend the results of
Robert Seeley concerning the complex powers of elliptic pseudodifferential 
operators and the residues of their zeta function to operators acting on
vector orbibundles over orbifolds.
We present the theory of pseudodifferential operators acting on a vector 
orbibundle over an orbifold, construct the zeta function of an elliptic 
pseudodifferential operator and show the existence of a meromorphic extension 
to $\CC$ with at most simple poles. We give formulas for generalized 
densities on the orbifold whose integrals compute the residues of the 
zeta function.
\par
Recall that a pseudodifferential operator $A$ of order $d$ acting on a 
hermitian vector bundle $E \overset{p}\to M$ over a Riemannian manifold $M$
of dimension $m$
becomes an unbounded operator on the space of $L^2$ sections $L^2(M;E)$. 
If $A$ is an elliptic pseudodifferential
operator with spectrum in $(\epsilon , \infty)$ for a
sufficiently small $\epsilon > 0$, by using the functional calculus one can 
define the complex powers $A^s$, $s\in \mathbb C$, as
\begin{equation}
A^s=\frac{1}{2\pi i}\int_{\gamma}\lambda^s(\lambda - A)^{-1} \,
d\lambda \qquad \text{when }Re(s) < 0
\end{equation}
(where $\gamma$ is a contour in the complex plane obtained by joining
two parallel half-lines to the negative real axis by a circle around
the origin) and
\begin{equation}
A^s=A^{s-k}A^k \qquad \text{for } Re(s)\ge 0
\end{equation}
for a large enough $k\in \mathbb Z$ that makes $s-k < 0$.
\par
The condition on the spectrum of $A$ can be relaxed; in order to apply 
the above construction, it is sufficient for 
the spectrum of $A$ to be inside an angle with the vertex the origin and 
away from the negative real axis. 
\par
In particular one can define the complex powers of a selfadjoint positive 
pseudodifferential operator $A$. Seeley showed in his paper \cite{seeley} 
that $A^s$ are pseudodifferential operators of complex order $sd$
and gave a local description of their symbols.
For $s$ in the half-plane $Re(s)<\frac{-m}{d}$, the operators $A^s$ are
of trace class. Seeley also showed that the zeta function of $A$ defined as
\begin{equation}
\zeta_A(s)=Tr(A^s)\qquad\text{ for } Re(s) < -\frac{m}{d}
\end{equation}
has a meromorphic extension to the complex plane with at most simple poles
at $\frac{-m+k}{d}$, with $k\in \{0,1,\dots\}$ and residues computable as 
integrals
on $M$ of quantities that depend only on the total symbol of the operator $A$.
These results will be referred to as Seeley Theory.
\par
An $m$-dimensional orbifold $M$ is locally the quotient of $\RR^m$ by a finite 
group of diffeomorphisms $\Gamma$ and globally it is the quotient of a 
manifold $\Tilde{M}^{m+q}$ by a compact Lie group $G$ of dimension $q$ with
finite isotropy groups. 
The dimension of the Lie group can be a priori quite large, like for 
example $\frac{m(m+1)}{2}$. 
Similarly, a rank $k$ vector orbibundle $E\overset{p}\to M$ can be represented 
locally as the quotient of a rank $k$ vector bundle by a finite group of 
bundle diffeomorphisms $\Gamma$, and globally as the quotient of 
a $G$ vector bundle of rank $k$,
$\Tilde{E}\overset{\tilde{p}}\to \Tilde{M}$, by a compact Lie group $G$ 
(cf. Proposition
\ref{gsection} and Theorem \ref{represent} ). We will  refer to the local 
description of the orbifolds and orbibundles
using the quotient by a finite group $\Gamma$ as  {\it Perspective 1} and 
to the global description using the quotient by a Lie group $G$ as 
{\it Perspective 2}. 
\par
From these perspectives the space of smooth sections in
a vector orbibundle can be thought of as the space of invariant smooth 
sections in
a genuine vector bundle over a manifold. One defines a pseudodifferential 
operator $A$
in a vector orbibundle as an operator acting on the space of invariant
sections in $\Tilde{E}\overset{\tilde{p}}\to \Tilde{M}$ which is the 
restriction
of a $\Gamma$ (respectively $G$) equivariant pseudodifferential operator 
$\Tilde{A}$ 
acting on $\Tilde{E}\overset{\tilde{p}}\to \Tilde{M}$. In the case of a finite 
group $\Gamma$, the operator $\Tilde{A}$ is unique up to a smoothing
operator, cf. Proposition \ref{unique}. If $dim(G)\ge 1$ then the operator
$\Tilde{A}$ which induces $A$ is far from being unique. 
\par
For the definition and the discussion of the basic properties of the 
pseudodifferential operators we will use
perspective $1$, while for some global properties and spectral theory of the
elliptic pseudodifferential operators we will use perspective $2$.
Seeley theory is about elliptic operators which is a global theory. However,
to formulate and prove the results concerning the pseudodifferentiability 
of the complex powers and the 
residues of the zeta function we use the local theory. 
Of particular importance are the explicit formulae in coordinate charts
for the residues of the zeta function. This makes clear the need of 
perspective $1$.
Nevertheless, the use only of perspective $1$ to extend Seeley theory in the 
context of orbibundles and orbifolds will require the repetition of the entire 
theory of elliptic pseudodifferential -- a rather tedious and very lengthy 
work. 
We will considerably short-cut this work by the use of perspective $2$.  
However, getting explicit formulae using perspective $2$ 
would be an almost impossible task
because of the large number of parameters introduced when constructing the 
operator $\Tilde{A}$ and because of the 
non-unique nature of $\Tilde{A}$. 
In fact even the extension of Seeley theory to incorporate the presence of a 
compact Lie group is rather complicated and very far from providing explicit 
formulas (see \cite{bh}).
\par
We give a brief description of the results in our paper.
\par
In the first section we consider a pseudodifferential operator of classical type 
$A$ on the space of smooth sections of a vector bundle $E\overset{p}\to M$ 
endowed with the a smooth action of a finite group $\Gamma$
such that $A$ commutes with this action and extend Seeley's results 
to this case. 
This extension, which is of independent interest, is also crucial for the study
of the pseudodifferential operators acting on vector orbibundles.
The results are probably implicit in the existing literature, but
not explicit enough for our needs in connection with the operators acting 
in orbibundles.
\par
The $\Gamma$ action on $E\overset{p}\to M$ induces the decomposition:
\begin{equation}
\label{decomp}
C^{\infty}(M;E)=\bigoplus_{i=1}^l C^{\infty}(M;E)_i 
\end{equation}
with 
$C^{\infty}(M;E)_i\overset{\sim}=V_i \otimes \Hom_{\Gamma}(V_i; 
C^{\infty}(M;E))$, where $(V_i, \rho_i)_{i=1}^l$ is a complete set of 
irreducible representations of $\Gamma$, with $\rho_1$ being the trivial
$1$-dimensional representation. The $\Gamma$ equivariance of $A$ induces 
the decomposition 
\begin{equation}
A=\bigoplus_{i=1}^l A_i
\end{equation}
with $A_i:C^{\infty}(M;E)_i \to C^{\infty}(M;E)_i$. This decompositions 
will be referred 
to as the decomposition along the irreducible representations of $\Gamma$. 
We restrict our attention 
to an elliptic positive pseudodifferential $A$ of order $d>0$. We are 
interested in the traces of the complex powers $(A^s)_i=(A_i)^s$ which 
can be recovered  from
$Tr(A^s\cdot \mathcal{T}_{\gamma})$ (as prescribed by formula 
\eqref{equation:i}), where $\mathcal{T}_{\gamma}$ is the operator
by which $\gamma\in \Gamma$ acts on $C^{\infty}(M;E)$. Although 
$\mathcal{T}_{\gamma}$
are not pseudodifferential operators, we showed that the trace functional
$\zeta_{A,\gamma}(s)=Tr(A^s\cdot \mathcal{T}_{\gamma})$, and therefore 
$\zeta_{a,i}(s)=Tr(A^s_i)$, defined for $Re(s)<-\frac{m}{d}$ can be extended
to a meromorphic function on $\CC$ with at most simple poles.
We construct a sequence of smooth densities $\eta^{\gamma}_k$
on the fixed point sets $M^{\gamma}=\{x\in M\,|\,\gamma x=x\}$ such that 
the residue of the meromorphic extension of $\zeta_{A,\gamma}$ at 
$s=\frac{-m+k}{d}$, $k\in\{0,1,\dots\}$, is computed as an integral on 
$M^{\gamma}$ of the density $\eta^{\gamma}_k$ (see Theorem 
\ref{residues:g} and 
\ref{residues:i}). Here $m=dim(M)$ and $d=ord(A)$. These densities can 
be interpreted as Dirac-type generalized densities on $M$(as in 
\cite{gs}, Chap. VI). For $k<(dim(M)-dim(M^{\gamma}))$ we have 
$\eta^{\gamma}_k=0$.
As a direct consequence we will show that the trace functionals 
$\zeta_{A,i}(s)=Tr(A_i^s)$  
can be extended to meromorphic functions with at most simple
poles at $s=\frac{-m+k}{d}$, $k\in\{0,1,\dots\}$ and compute the residues in
Theorem \ref{residues:g}.
Of particular interest is $Tr(A_1^s)$ -- the component 
corresponding to the trivial one dimensional representation $\rho_1$, because
the $\Gamma$ equivariant operator $A$ induces an elliptic pseudodifferential 
operator acting on sections 
of the orbibundle $E/\Gamma\overset{\overline{p}}\to M/\Gamma$ which can
be identified with $A_1$. The trace functional $Tr(A_1^s)$ is the zeta 
function of this operator.
\par
In the second section we present the basic definitions and properties
of the orbifolds and orbibundles.
It is a known fact that any 
orbifold $M^m$ is the quotient of a smooth manifold $\Tilde{M}^{m+q}$ by a 
compact Lie group $G$ of dimension $q$ acting on $\Tilde{M}$ with finite 
isotropy groups.
We complete this result by proving that a vector orbibundle 
$E\overset{p}\to M$ 
is the quotient by $G$ of a $G$ vector bundle $\Tilde{E}\overset{\Tilde{p}}
\to \Tilde{M}$ in Theorem \ref{represent}. This will allow us to give a global 
characterization for the spaces of smooth sections in orbibundles in 
Proposition \ref{gsection}. We also describe 
the space of generalized Dirac-type densities on orbifolds and the
canonical stratification associated with an orbifold structure.
\par
In the third section we show how one can extend the theory of 
pseudodifferential operators acting on smooth sections in a vector bundle 
over a manifold to pseudodifferential operators 
acting on smooth sections in a vector orbibundle over an orbifold. 
The basic definitions and elementary properties were stated, though with 
some gaps, in \cite{GN1} and \cite{GN2}. Next we discuss the complex powers
and the zeta function of such an operator which is selfajoint and elliptic 
and extend Seeley theory to this context.
We show that the zeta function is well defined and has a meromorphic 
extension to the complex plane with at most simple poles located at 
$s=\frac{-m+k}{d}$, $k\in\{0,1,\dots\}$. In Theorem \ref{zeta:Dirac:density},
using the tools presented in the previous sections, we construct a sequence of
generalized Dirac-type densities $\eta_k$, $k\in\{0,1,\dots\}$, on the 
base orbifold $M$
and compute the residues of the zeta function as integrals of these densities. 
In Theorem \ref{zeta:strata:density} we reinterpret the computations we make 
as integrals of genuine smooth densities on the strata of the canonical 
stratification of $M$. 
\par
The results of this paper provide the analytic foundation for the study of 
indices, signatures and torsion for orbifolds equipped with a Riemannian 
metric (Riemannian orbifolds), which will provide the topic of a subsequent 
paper. We recall the reader that many of the moduli spaces  have a canonical 
structure of an orbifold (rather that manifold)
and their topological and geometric invariants are of legitimate interests.
\par
Finally, I want to thank D. Burghelea for the help and advice in preparing 
this work. In many respects he is the coauthor of this paper. 


\section{$\Gamma$ Equivariant Pseudodifferential Operators}
\subsection{$\Gamma$ Vector Bundles}
\par
\begin{defn}
Let $(\Gamma,\cdot)$ be a finite group of order $|\Gamma|$. A real or
complex smooth vector bundle $E\overset{p}\to M$ endowed with a left 
smooth $\Gamma$ action by bundle isomorphisms is 
called a $\Gamma$ vector bundle.
\end{defn}
\par
For each $\gamma \in \Gamma$ we have a diffeomorphism $t_\gamma : M\to M$ 
of the base space and a linear isomorphism $(T_\gamma )_{x}: E_x \to 
E_{t_\gamma (x)}$ between the fiber above $x$ and the fiber above $t_\gamma
(x)$. As usual, we will denote $t_\gamma (x)$ by $\gamma x$. $\Gamma$ 
acts on the space of sections $C^{\infty}(M;E)$ by $(\gamma \cdot f)
(x)=( T_\gamma )_{\gamma ^{-1} x}f(\gamma ^{-1} x)$. 
We will denote the action of $\gamma $ on the space of sections by 
$\mathcal T_\gamma $. 
\par
\begin{example}
\label{gammalinear}
A particular case of a $\Gamma_0$ vector bundle is the trivial vector bundle
$U\times V \overset{pr_1}\to U$ endowed with  a product $\Gamma_0$ action 
$\mu\times \rho$ on $U\times V$ and with the $\Gamma_0$ action $\mu$ on $U$ 
where $\mu:\Gamma_0\times U \to U$ is the restriction of a linear action $\mu:
\Gamma_0\times\mathbb{R}^m\to \mathbb{R}^m$ to an open invariant neighborhood 
of the origin $U\subset  \mathbb{R}^m$ and $\rho:\Gamma_0\times V\to V$ is a 
linear representation on the vector space $V$.
\end{example}
\par
\begin{example}
\label{gammalinear_ext}
If $\Gamma_0\subset \Gamma$ is an inclusions of groups, then the bundle
$\Gamma\times_{\Gamma_0}(U\times V)\overset{pr_1}\to \Gamma\times_{\Gamma_0}U$
is a $\Gamma$ vector bundle, with $\Gamma$ acting on the total space and
base space by left translations.
\end{example}
\par 
The following proposition states that any $\Gamma$ vector bundle is 
locally diffeomorphic to a  $\Gamma$ vector bundle as above.
 
\begin{prop}
\label{local}
Let $E\overset{p}\to M$ be a $\Gamma$ vector bundle and $x\in M$. Let 
$\Gamma_x$ be the isotropy group of $x$. Then there exists $U$ a $\Gamma_x$ 
invariant neighborhood of $x$ in $M$ and $O$ a neighborhood of the origin in 
$\mathbb{R}^m$ such that the restriction bundle $E_{|\Gamma U}\overset{p}\to 
\Gamma U$ is isomorphic to the $\Gamma$ vector bundle $\Gamma\times_{\Gamma_x}
(O\times V)\overset {p_1}\to \Gamma\times_{\Gamma_x}O$  by a $\Gamma$ 
equivariant isomorphism.
\end{prop}

The proof of this proposition is elementary and rather standard. For the sake 
of completeness we included it in the Appendix A.

Let $E\overset{p}\to M$ be a $\Gamma$ vector bundle. Let $(V_1,\rho_1), 
(V_2,\rho_2),...,(V_l,\rho_l)$ be a complete set of irreducible 
non-isomorphic complex representations of $\Gamma$, with $\rho_1$ being the 
trivial one dimensional representation. Then $C^{\infty}(M;E)$ 
decomposes as the direct sum of multiples of $V_i$ 
\begin{equation}
C^{\infty}(M;E)=\bigoplus_{i=1}^l C^{\infty}(M;E)_i
\end{equation}
$C^{\infty}(M;E)_i$ is spanned by submodules of the 
form $h(V_i)$ where $h: V_i\to C^ {\infty}(M;E)$ are $\Gamma$ equivariant 
maps. We have the isomorphism 
\begin{equation}
V_i \otimes \Hom_{\Gamma}(V_i; C^{\infty}(M;E))
\overset{e}\to C^{\infty}(M;E)_i
\end{equation}
where $e$ is the evaluation map, $e(v,\phi)=\phi(v)$.
If $M$ is a closed Riemannian manifold endowed with a $\Gamma$ invariant 
metric, and we have a $\Gamma$ invariant Hermitian structure on the bundle 
$E\overset{p}\to M$, then $\Gamma$ acts on the space of $L^2$ integrable 
sections $L^2(M;E)$, and we get an analogous decomposition 
\begin{gather}
L^2(M;E)=\bigoplus_{i=1}^{l}L^2(M;E)_i\\
\intertext{with}
L^2(M;E)_i=V_i\otimes \Hom_{\Gamma}(V_i; L^2(M;E))
\end{gather}
We will denote by $pr_i$ the projection on the $i^{th}$ factor $L^2(M;E)_i$ 
with respect to the above decomposition.

We will need the following statement that describes the projection $pr_i$ in 
terms of the action and characters of the group $\Gamma$.

\begin{prop}
\label{pri}
If $\chi_i$, $i=1,\dots,l$ is the complete set of irreducible characters of the
group $\Gamma$ corresponding to the representations $(V_i,\rho_i)$ and $k_i=
dim_{\mathbb{C}}V_i$ are the corresponding dimensions, then the projection on 
the $i^{th}$ factor is given by
\begin{equation}
pr_i=\frac {k_i}{|\Gamma|} \sum_{\gamma \in \Gamma} \chi_i(\gamma ^{-1})\cdot 
\mathcal{T}_\gamma 
\end{equation} 
\end{prop}

The proof can be found in \cite{serre}, Part 1, Theorem 8. 

\subsection{$\Gamma$ Equivariant Pseudodifferential Operators}

\begin{defn}

A pseudodifferential operator $A$ acting 
on the space of sections $C^{\infty}(M;E)$ of a $\Gamma$ vector bundle is 
called $\Gamma$ equivariant if $A$ commutes with the action of $\Gamma$ on 
$C^{\infty}(M;E)$ i.e. 
$\mathcal{T}_\gamma \cdot A= A\cdot \mathcal{T}_\gamma $ for any $\gamma \in 
\Gamma$.
\end{defn}

Throughout this chapter we suppose that $A$ is elliptic 
$\Gamma$ equivariant pseudodifferential operator of classical type 
(as described in \cite{shubin}, Section 3.7), of positive order $d$. 
We suppose that $\pi$ is an Agmon angle for $A$, i.e. there exists 
$\varepsilon>0$ such that the spectrum of $A$ is disjoint from the region
in the complex plane $\{z\,|\, arg(z)\in(\pi-\varepsilon, \pi+\varepsilon)\}
\cup\{z\,|\,|z|<\varepsilon\}$.
Then $A$ and all its complex powers $A^s$, $s\in\mathbb C$, will 
preserve the decompositions $C^{\infty}(M;E)=\bigoplus_{i=1,l}
C^{\infty}(M;E)_i$ and $L^2(M;E)=\bigoplus_{i=1}^{l}L^2(M;E)_i$. 
Consequently, we can consider the restrictions $A^s_i:C^{\infty}(M;E)_i\to
C^{\infty}(M;E)_i$ and $A^s_i:L^2(M;E)_i\to L^2(M;E)_i$ for $s\in\mathbb C$.
The goal of this section is the study of the trace of these operators 
$\zeta_{A,i}(s)=Tr (A^s_i)$. First observe that $A^s_i=A^s\circ 
pr_i$. Using Proposition \ref{pri}, we get
\begin{equation}
\label{equation:i}
A^s_i=\frac {k_i}{|\Gamma|}\sum_{\gamma \in \Gamma}\chi_i(\gamma ^{-1}) 
A^s\circ\mathcal{T}_\gamma 
\end{equation}
so, in order to study the trace of $A^s_i$, it is 
convenient to study the trace of $A^s\circ\mathcal{T}_\gamma$. We will 
show that these operators are of trace class for all complex numbers $s$ in a 
half-plane $Re(s)<-K$, and that the associated trace function, 
$Tr(A^s\circ\mathcal{T}_\gamma)$,
as a function of $s$, has a meromorphic extension to $\CC$ with at most 
simple poles. 
\par
Let us fix an element 
$\gamma\in \Gamma$ and denote by $ T = T_\gamma $, $t=t_\gamma 
$ and $\mathcal{T}=\mathcal{T}_\gamma $ 
respectively the action of $\gamma $ on the total space, on the base 
space and on the space of sections of the bundle. Let $M^{\gamma}$ be 
the fixed point set of the
diffeomorphism $t$ and $M^{\gamma}=\bigcup_{i\in I}M^{\gamma}_i$ be the 
decomposition in connected components of dimensions respectively $n_i$. 
Let $m=dim(M)$.
\par
\begin{thm}
\label{residues:g}
The operators $A^s\circ\mathcal{T}$ are of trace class for $s$ in 
the half-plane $Re(s)<\frac {-m}{d}$. The associated zeta function 
$\zeta_{A,\gamma}(s)=Tr (A^s\circ\mathcal{T})$, defined for 
$Re(s)<\frac {-m}{d}$, has a meromorphic continuation to the whole 
complex plane $\mathbb{C}$ with at most simple poles at $\frac{-m+k}{d}$ for
$k\in\{0,1,2,...\}$. One can construct positive numbers $d^{\gamma}_i$ and 
densities $\{\eta^{\gamma}_{i,k}\}$ for  $k\in\{0,1,2,\dots\}$ on the submanifolds 
$M^{\gamma}_i$ such that $\eta^{\gamma}_{i,k}=0$ for $k<m-n^{\gamma}_i$ 
and the residue of the function $\zeta_{A,\gamma}$ at $s=\frac{-m+k}{d}$ is 
equal to
\begin{equation}
\label{res:g}
res_{|s=\frac{-m+k}{d}}\zeta_{A,\gamma}=\sum_{i}d^{\gamma}_i\int_{M^{\gamma}_i}
\eta^{\gamma}_{i,k}
\end{equation}
\end{thm}
\par
The details of the proof of this Theorem are contained in Appendix B. 
We will give the description of the coefficients $d^{\gamma}_i$ and the 
densities $\eta^{\gamma}_{i,k}$.
\par
For $\gamma\in \Gamma$ and $M^{\gamma}_i$ a connected component of the fixed 
point set $M^{\gamma}$, let $x\in M^{\gamma}_i$. Using a $\Gamma$ invariant 
metric in the tangent space $T(M)$ one decomposes $T_x(M)=T_x(M^{\gamma}_i)
\oplus T_x(M^{\gamma}_i)^{\perp}$ and the action of $\gamma$ on $T_x(M)$
as $Id\oplus \overline{t}$. Then we define: 
\begin{equation}
\label{defn:d.gamma}
d^{\gamma}_i=|det(\overline{t}-Id)|^{-1}.
\end{equation}
The quantity above does not depend on the particular choice of 
$x\in M^{\gamma}_i$
and decomposition $T_x(M)=T_x(M^{\gamma}_i)\oplus T_x(M^{\gamma}_i)^{\perp}$.
\par
In order to define the densities $\eta_{i,k}^{\gamma}$ near $x\in M^{\gamma}_i$
we choose  a coordinate chart $\phi:(O,0)\to (\phi(O),x)$ in an open 
neighborhood $\phi(O)$ of $x\in M$ such that the induced action 
of $\Gamma_x$ on $O\subset \RR^m$ is given by a linear orthogonal maps. 
This can be realized with the help of 
a $\Gamma$ invariant metric on $M$ and the associated exponential map at $x$.
Denote the action of $\gamma$ on $O$ by $t$. 
Let $(x_1,x_2)$ be coordinates given by this 
chart, $x_1\in O^t$--the fixed point set of $t$, and $x_2\in O^{\perp}$. 
Let $(\xi_1,\xi_2)$ be the
corresponding coordinates in the cotangent bundle space. Then $t=Id\oplus 
\overline{t}$. Observe that $d^{\gamma}_i=|det(\overline{t}-Id)|^{-1}$.
Let $a_s(x_1,x_2,\xi_1,\xi_2)$ be the total symbol of $A^s$ in $O$ and
\begin{equation}
a_s(x_1,x_2,\xi_2,\xi_2)\sim\sum_{k\ge 0}a_{s,k}(x_1,x_2,\xi_1,\xi_2)
\end{equation}
be its asymptotic expansion (as defined in 
\cite{shubin}, Section 3.7). The component $a_{s,k}$ is homogeneous in
$(\xi_1,\xi_2)$ of degree of homogeneity $sd-k$. 
\par
Let $\overline{a}_{s,k}(x_1,w,\xi_1\xi_2)=a_{s,k}(x_1,(\overline{t}-Id)^{-1}w,
\xi_1,\xi_2)$.
\par
Consider the homogeneous symbol $\Tilde{b}_{s,j}(x_1,\xi_1)$ of degree of 
homogeneity $sd-j$ given by
\begin{equation}
\Tilde{b}_{s,j}(x_1,\xi_1)=\sum_{|\alpha|+k=j}\frac{1}{\alpha !}(D^{\alpha}_w
\partial^{\alpha}_{\xi_2}\overline{a}_{s,k})(x_1,0,\xi_1,0)
\end{equation}
Then the density on $O^t$ whose integral computes the residue of 
$\zeta_{A,\gamma}$ at $s=\frac{-m+k}{d}$ is given by
\begin{align}
\label{defn:eta.def1}
\eta^{\gamma}_{i,k}(x_1)&=-\frac{1}{d}Tr(\int_{S^{n-1}}\Tilde{b}_{s,n_i-m+k}
(x_1,\overline{\xi})\,\dbar\overline{\xi}\circ T)\,dx_1&&\text{ if } 
k\ge m-n_i\\
\label{defn:eta.def2}
\eta^{\gamma}_{i,k}&=0 &&\text{ if } k< m-n_i
\end{align}
Here $n_i=dim(O^t)=dim(M^{\gamma}_i)$ and $T=T_{\gamma,x}$ is the map by 
which $\gamma$ acts in the fiber $E_x$ above $x\in M$. The $(n-1)$ form 
$\dbar\overline{\xi}$ is the canonical volume form on $S^{n-1}$ induced from
$\RR^n$, rescaled by a factor of $(2\pi)^{-n}$. $dx_1$ is the canonical volume 
form on $O^t$.

Using Theorem \ref{residues:g} and the formula \eqref{equation:i}
linking the operators $A^s_i:C^{\infty}(M;E)_i\to C^{\infty}(M;E)_i$ and 
$A^s\circ\mathcal{T}_{\gamma}$, for $\gamma \in \Gamma$, we can formulate 
the following theorem

\begin{thm}
\label{residues:i}
The operator $A^s_i$ is of trace class and the trace functional
\begin{equation}
\zeta_{A,i}(s)=Tr(A^s_i)
\end{equation}
is a holomorphic function in $s$ on the half-plane $Re(s)<-\frac{m}{d}$.
The function $\zeta_{A,i}$ has a meromorphic continuation to the whole
complex plane with at most simple poles situated at $s=\frac{-m+k}{d}$ 
for $k\in\{0,1,2,\dots\}$. For each $k\in\{0,1,2,\dots\}$ there exist smooth
densities $\eta^{\gamma}_k$ on the fixed point set set $M^{\gamma}$ 
such that the residue of $\zeta_{A,i}$ at $s=\frac{-m+k}{d}$ is equal to
\begin{equation}
\label{res:i:g}
res_{|s=\frac{-m+k}{d}}\zeta_{A,i}=\frac {k_i}{|\Gamma|}
\sum_{\gamma \in \Gamma}
\chi_i(\gamma ^{-1})\int_{M^{\gamma}}\eta^{\gamma}_k
\end{equation}
The smooth density $\eta^{\gamma}_k$ depends only on a finite number
of terms in the asymptotic expansion of the total symbol of
the operator $A^s$ and the action of $\Gamma$ on a neighborhood 
of $M^{\gamma}$.
\end{thm}

\begin{proof}
Formula \eqref{equation:i} implies: 
\begin{equation}
\zeta_{A,i}(s)=\frac {k_i}{|\Gamma|}\sum_{\gamma \in \Gamma}\chi_i
(\gamma ^{-1}) \zeta_{A,\gamma}(s)
\end{equation}
The first part of the statement in the theorem follows directly from Theorem 
\ref{residues:g}. If we define $\eta^{\gamma}_k=\sum_i d_i^{\gamma}
\eta^{\gamma}_{i,k}$ then \eqref{res:i:g} is a direct consequence of 
the equality \eqref{res:g}

\end{proof}

\subsection{Dirac-Type Densities on $\Gamma$ Manifolds}
 
Let us consider $C^{\infty}(M)$  the space of smooth functions on a compact 
manifold $M$ and endow it with the topology of uniform convergence together 
with a finite number of derivatives on $M$.
\begin{defn} Let $N$ be a smooth submanifold which is a closed subset of $M$. 
A smooth density $\eta$ on $N$ defines a continuous functional on the 
space $C^{\infty}(M)$ by
\begin{equation}
\label{pairing}
<\eta, f>=\int_N f_{|N}\eta \qquad \text{ for } f \in C^{\infty}(M)
\end{equation}
We call such a functional a Dirac-type distribution on $M$. 
\end{defn}

The singular support of this distribution is equal to $N$
if $dim(N)< dim(M)$ and it is empty if $dim(N)=dim(M)$.

If $N=M$ then a Dirac-type distribution on $M$ is given by smooth density 
on $M$. If $N$  is a proper  subset of $M$ then a Dirac-type distribution 
on $N$ can be described locally by a density that is zero on open sets 
disjoint from $N$ and by a smooth density on $N$. This will not be a 
smooth density on $M$ anymore, the singular set being exactly $N$. We 
call this density a Dirac-type density and denote it with the same letter 
as the smooth density on $N$.
 
If $\phi:M\to M$ is a diffeomorphism and $\eta$ is a Dirac-type density 
with singular support $N$, then the push-forward $\phi_*(\eta)$ is a 
Dirac-type density associated with the smooth density on $\phi(N)$ which 
is the push-forward of the smooth density $\eta$ on $N$ by the 
diffeomorphism $\phi_{|N}:N \to\phi(N)$. 
We have the following identity
\begin{equation} 
<\phi_*(\eta), \phi_*(f)>=<\eta ,f>.
\end{equation}

\begin{defn} 
\label{gamma:density}
If $\Gamma\times M \to M$ is a $\Gamma$ smooth manifold, a $\Gamma$ 
Dirac-type density on $M$ is a collection of Dirac-type densities $\eta=\{
\eta^{\gamma}\}_{\gamma\in\Gamma}$
indexed by the elements of the group $\Gamma$ with singular support 
respectively the fixed point sets $M^{\gamma}$ and such that 
$\gamma'_*(\eta^{\gamma})=
\eta^{\gamma'\gamma \gamma'{}^{-1}}$ for any $\gamma, \gamma'\in \Gamma$
(observe that, in general,  $\gamma'(M^{\gamma})=(M^{\gamma' \gamma
\gamma'{}^{-1}})$).
\end{defn}

\begin{defn} If $\eta=\left\{\eta^{\gamma}\right\}_{\gamma\in\Gamma}$ is a 
$\Gamma$ Dirac-type density on a $\Gamma$ manifold $M$ and $\chi=\chi_i$ 
is the character of an irreducible representation of $\Gamma$ we define 
the associated distribution $\eta^{\chi}$ by
\begin{equation}
<\eta^{\chi}, f>=\frac{k_i}{|\Gamma|}\sum_{\gamma \in \Gamma}\chi_i
(\gamma ^{-1})\cdot \int_{M^{\gamma}}f_{|M^{\gamma}} \eta^{\gamma} 
\end{equation}
\end{defn}

Using the above definitions, Theorem \ref{residues:i} can be reformulated as 
follows:
\begin{thm} 
\label{residues:d}
Let $A$ be a $\Gamma$ equivariant 
elliptic pseudodifferential operator acting on the smooth sections of a 
$\Gamma$ bundles manifold and $A_i$ be the restriction of $A$ to the 
component of the space of sections $C^{\infty}(M;E)$ corresponding to 
the irreducible representation $(V_i,\rho_i)$ with the character $\chi_i$. 
The trace functional 
$\zeta_{A,i}(s)=Tr(A_i^s)$ defines a holomorphic function on the half-plane
$Re(s)<\frac{-m}{d}$ which has a meromorphic extension to the whole 
complex plane with at most simple poles situated at $s=\frac{-m+k}{d}$, for 
$k\in\{0,1,2,\dots\}$.
There exists a family $\{\eta_k\}_{k=0}^{\infty}$ of $\Gamma$ Dirac-type 
densities on $M$ with $\eta_k=\{\eta^{\gamma}_k\}_{\gamma\in\Gamma}$
so that the residue of $\zeta_{A,i}$ at $s=\frac{-m+k}{d}$ is equal to 
$<\eta^{\chi_i}_k,1>$.
\end{thm}

\begin{proof}
We only have to prove the existence of the representation of the residues 
of $\zeta_{A,i}$ as stated in the theorem.

For $k\in\{0,1,2,\dots\}$ and $\gamma\in\Gamma$ let $\eta_k^{\gamma}$ be the 
Dirac-type distribution given by the sum of the smooth densities 
$\sum_{i}d_i^{\gamma} \eta_{i,k}^{\gamma}$ each defined on the connected 
component $M^{\gamma}_i$ of the fixed point set $M^{\gamma}$ as described 
in Theorem \ref{residues:g} and by the formulas \eqref{defn:d.gamma},
\eqref{defn:eta.def1} and \eqref{defn:eta.def2}.
We have $\gamma^{-1}(M^{\gamma\gamma'\gamma^{-1}})=M^{\gamma'}$ for any
$\gamma, \gamma'\in\Gamma$, and after a 
convenient reindexing of the connected components of each fixed point set, we
can suppose that $\gamma'{}^{-1}(M^{\gamma'\gamma\gamma'{}^{-1}}_i)=
M^{\gamma}_i$ as well. Then $d_i^{\gamma'\gamma\gamma}=d_i^{\gamma}$ 
because, as defined by formula \eqref{defn:d.gamma}, they are determinants of 
two linear maps conjugated by a diffeomorphism defined at the tangent 
space level by $\gamma'$. Because the operator $A$ is $\Gamma$ equivariant, 
a straightforward computation 
shows that the smooth densities $\eta_{i,k}^{\gamma}$ and $\eta_{i,k}^{\gamma'
\gamma\gamma'{}^{-1}}$ are conjugated by the the map induced at the cotangent
space level by $\gamma'$. Then the collection $\{\eta_k^{\gamma}\}_{\gamma
\in\Gamma}$ is a $\Gamma$ Dirac-type distribution, as described in the 
definition \ref{gamma:density}. 
Theorem \ref{residues:i} states that the residue of $\zeta_{A,i}$ at
$s=\frac{-m+k}{d}$ is equal to 
\begin{equation}
\frac{k_i}{|\Gamma|}\sum_{\gamma\in\Gamma}\chi_i(\gamma^{-1}) 
(\sum_{i}d_i^{\gamma} \int_{M_i^{\gamma}}\eta_{i,k}^{\gamma})
\end{equation}
which can be rewritten as $<\eta^{\chi_i}_k, 1>$.
\end{proof}

\vfill\eject

\section{Orbifolds and Orbibundles}
\subsection{Orbifolds}
Let $M$ be a Hausdorff topological space. 
\begin{defn}
An orbifold chart on $M$ is
given by $\mathcal{R}=(\Tilde U, \Gamma, \mu, U, \pi)$
where $U$ is an open set in $M$, $\Tilde U$ is an open set in the 
$n$-dimensional
Euclidean space $\RR^n$, $\Gamma$ is a finite group, $\mu:\Gamma\times
\Tilde{U}\to\Tilde{U}$ is a faithful smooth action of $\Gamma$ on $\Tilde{U}$ 
and $\pi :\Tilde U \to U$
is a continuous map that factors through a homeomorphism
$\overline\pi$ from the orbit space $\Tilde U /\Gamma$ to $U$ and
makes the following diagram commutative:
$$
\xymatrix{
\Tilde{U} \ar[rr]^{\pi} \ar @{->>}[dr]& &U\\
&\Tilde{U}/\Gamma\ar @{-->}[ur]^{\overline\pi}_{\sim}&}
$$
The open set $U$ is called a coordinate neighborhood.
\end{defn} 

If the group action $\mu$ is the restriction of a linear representation of 
$\Gamma$ on $\RR^n$ to a neighborhood of the origin $\Tilde{U}$ and
$x=\pi(0)$, we call $\mathcal{R}$ a linear chart at $x$.

\begin{defn} Two orbifold charts $\mathcal{R}_i=(\Tilde U_i, \Gamma_i, \mu_i,
U_i, \pi_i)$, $i=1,2$ are called compatible if for any two points 
$\tilde x_i \in \Tilde U_i$ such that $\pi_1(\tilde x_1)=\pi_2(\tilde x_2)=x
\in U_1\cap
U_2$ there exists a diffeomorphisms $h$ from a neighborhood of $x_1$ in
$\Tilde{U}_1$ to a neighborhood of $x_2$ in $\Tilde{U}_2$ such that 
$\pi_2\circ h=\pi_1$.
\end{defn} 

\begin{remark} 
\label{linchart}
If $\mathcal{R}=(\Tilde U, \Gamma,\mu,U, \pi)$ is an orbifold
chart and $x\in U$, there exist a linear orbifold chart
$\mathcal{R}_x$ at $x$ which is compatible with $\mathcal{R}$. The proof follows 
the same idea as the one in the proof of Proposition \ref{local}, 
contained in Appendix A. Let $\tilde x \in \Tilde{U}$ 
such that  $\pi(\tilde x)=x$ and $\Gamma_{\tilde x} \subset\Gamma$ be 
the isotropy group of $\tilde x$.
Let $\Tilde{V}$ be a $\Gamma_{\tilde x}$ invariant neighborhood of
$\tilde x$ such that $\Tilde{V}\cap \gamma\cdot\Tilde{V} = \emptyset$
for any $\gamma \in \Gamma \backslash \Gamma_{\tilde x}$.  Then 
$(\Tilde{V}, \Gamma_{\tilde x}, \mu, V=\pi(\Tilde{V}), \pi_{|\Tilde{V}})$
is an orbifold chart which is compatible with $\mathcal{R}$.
The linearization of  the action of $\Gamma_{\tilde x}$ on $\Tilde{V}$
gives us an action $\mu'$ of $\Gamma_{\tilde x}$ on $T_{\tilde x}(\Tilde{V})$,
the tangent space of
$\Tilde{V}$ at $\tilde x$ . The exponential map $exp:T_{\tilde x}(\Tilde{V})
\to \Tilde{V}$ associated with a $\Gamma_{\tilde x}$ 
invariant metric on $\Tilde{V}$ is $\Gamma_{\tilde x}$ equivariant.
One can choose a smaller $\Tilde{V}$ so that
$exp$ is an equivariant diffeomorphism between a neighborhood 
$\Tilde{W}$ of $0$ in $T_{\tilde x}(\Tilde{V})$ and $\Tilde{V}$. Then
$(\Tilde{W}, \Gamma_x, \mu', V, \pi_{|\Tilde{V}}\circ exp)$ 
is an linear orbifold chart at $x$ which is compatible with $\mathcal{R}$.
\end{remark}

\begin{defn} An orbifold atlas $\mathcal{A}$ on $M$ is a collection of 
compatible orbifold charts on $M$ such that the corresponding coordinate 
neighborhoods form an open cover of $M$.
Two orbifold atlases $\mathcal{A}_i,\;i=1,2$ on $M$ are compatible if their 
reunion is an orbifold atlas of $M$.
\end{defn}

\begin{defn} An orbifold structure on $M$ is given by a maximal orbifold
atlas $\mathcal{A}$ on $M$.
\end{defn}

Though an orbifold is given by a topological space and a maximal atlas, 
in the future we will drop the atlas from the notation of an orbifold and use
only the letter designated for the underlying topological space. In order to 
keep our notations simple we will also drop the representation $\mu$ 
from the notation of an orbifold chart.
 
\begin{example}
If $W$ is an open subset of $M$, then $W$ inherits an orbifold
structure whose charts are all orbifold charts $(\Tilde{U}, \Gamma, U, 
\pi)$ of $M$ for which $U\subset W$. 
\end{example}

\begin{example}
\label{cano} 
Let $M$ be a differentiable manifold and $\Gamma$ a
finite group of diffeomorphisms of $M$. Let $\pi$ be the projection map
onto the orbit space $M/\Gamma$. Then 
$M/\Gamma$ has a canonical orbifold structure (cf. \cite{haef}). An atlas of 
$M/\Gamma$ can be obtained as follows:
let $x\in M$ and $\overline{x}\in M/\Gamma$ its orbit. Let $\Gamma_x
\subset\Gamma$ be the isotropy group of $x$. Because $\Gamma$ is
finite, there exists an open neighborhood $\Tilde{U}$ of $x$ in $M$ 
which is $\Gamma_x$ invariant and $\Tilde{U}\cap \gamma\cdot\Tilde{U}
= \emptyset$ for $\gamma \in \Gamma\backslash \Gamma_x$. The set
$U=\Tilde{U}/\Gamma_x \subset M/\Gamma$ is an open neighborhood of
$\overline{x}$. We call a neighborhood like $\Tilde{U}$ a slice at $x$. Choose $\Tilde{U}$ 
small enough such that $\Tilde{U}\overset{\phi} 
\to O\subset \RR^n$ is a smooth chart for $M$. The action of the group
$\Gamma$ can be transported to a smooth action on $O$ via $\phi$
and the collection $(O, \Gamma_x, U, \pi\circ\phi^{-1})$ is an
orbifold chart of $M/\Gamma$. All the charts obtained as above are compatible 
(cf.  \cite{haef}) and the maximal atlas containing them defines the
 canonical orbifold structure on $M/\Gamma$.
\end{example}

We have a generalization of the previous example:
 
\begin{prop}
\label{liecano}
Let $M$ be a differentiable manifold and $\mu:G\times M \to M$ a smooth action 
of a compact Lie group $G$ with finite isotropy groups $G_x\subset G$ for any 
$x\in M$. Then the quotient topological space $M/G$ has a canonical structure 
of an orbifold.
\end{prop}

\begin{proof}
Let us fix a $G$ invariant metric on $M$. This can be done by averaging any 
metric over the compact group $G$.

Let $\overline{x}\in M/G$ be a point in the quotient space and $x\in M$ such 
$\overline{x}=Gx$. Let $G_x$ be the isotropy group at $x$. Then $G_x$ acts on the 
tangent space $T_x(M)$ and keeps invariant the tangent space $T_x(Gx)$ to the 
$G$-orbit of $x$. Let $V\subset T_x(M)$ be the orthogonal complement of $T_x(Gx)$
in $T_x(M)$, $V\oplus T_x(Gx)=T_x(M)$. Because the metric on $M$ is $G$ 
invariant, $G_x$ will act on $V$ by restriction. Also, all  its translations 
$gV\subset T_{gx}(M)$ will be $G_{gx}$ invariant and $gV\oplus T_{gx}(Gx)=
T_{gx}(M)$. 

Let $T(M)_{|Gx}\to Gx$ be the restriction of the tangent vector bundle to the orbit 
$Gx$, and $\mathcal{V}\to Gx$ be the subbundle whose fiber above $gx$ is 
$gV$. This subbundle has a natural $G$ vector bundle structure coming from 
the action of $G$ on $M$.

Let us consider the principal bundle $G_x\hookrightarrow G\to G/G_x$, where 
$G_x$ acts by right translations on $G$, and the associated vector bundle
$V\hookrightarrow G\times_{G_x}V \to G/G_x$. This vector bundle has a natural 
$G$ vector bundle structure, $G$ acting by left translations on $G$ and 
$G/G_x$.  

The maps $\Phi: G\times_{G_x} V \to \mathcal{V}$, $\Phi(g,v)=gv$  and  
$\phi:G/G_x \to Gx$, $\phi(gG_x)=gx$ give a $G$-equivariant 
isomorphism $(\Phi,\phi)$ between the two $G$ vector bundles considered above.

If $exp:T(M) \to M$ is the exponential map associated with the $G$ 
invariant metric
on $M$, then $exp$ realizes a $G$ equivariant diffeomorphism between a 
neighborhood $\mathcal{U}$ of the zero section in $\mathcal{V}$ and an 
open tubular neighborhood 
$N$ of the orbit $Gx$ in $M$. Because $G$ is compact one can find an open 
$G_x$ invariant neighborhood $U$ of the origin in $V$ such that 
$exp\circ \Phi: G\times_{G_x} U\to N$ is a $G$ equivariant diffeomorphism. 
We will then take $\mathcal{U}=\Phi(G\times_{G_x} U)$.
Passing to the $G$-orbit spaces, we get a homeomorphism $\overline
{exp\circ \Phi}:(G\times_{G_x}U)/G \to N/G$, where $N/G$ is an open 
neighborhood of $\overline{x}=Gx$ in $M/G$. We will construct an orbifold chart
over $N/G$ around $\overline{x}$.

The map $\iota: U\to G\times_{G_x}\!
U$, $\iota(u)=(e,u)$ gives a homeomorphism when passing to the orbit spaces
$\overline{\iota}: U/G_x \overset{\sim}\to (G\times_{G_x}\!U)/G$.
Denote by $\pi$ the composition $U\overset{\iota}\to G\times_{G_x}U \overset{
exp\circ\Phi}\longrightarrow N \overset{proj}\longrightarrow N/G$. Then
$(U, G_x, N/G, \pi)$ is a linear orbifold chart at $\overline{x}$, where the 
action of $G_x$ on  $U\subset V$ is the restriction of the linear 
representation of $G_x$ on $V$. 
As shown above, the induced map to the orbit spaces $\overline{\pi}:
U/G_x \to N/G$ is a homeomorphism.

We have to show that any two different charts defined as above are compatible. 

Let $\mathcal{R}_i=(U_i, G_{x_i}, N_i/G, \pi_i)$ be two orbifold 
charts around $\overline{x}_i$, $i=1,2$. Let $u_i\in U_i$ such that 
$\pi_1(u_1)=\pi_2(u_2)=\overline{x}
\in M/G$. Then, using the definition of $\pi$ in the construction of the 
charts $\mathcal{R}_i$, one can find $g_i\in G$, $i=1,2$ and $x\in M$,  
$Gx=\overline{x}$, so that we have $exp_{g_1 x_1}(g_1 u_1)=exp_{g_2 x_2}
(g_2 u_2)=x\in M$. By replacing $x$ with $g_1^{-1}x$ we can assume that $g_1=e\in G$. 
Moreover, the map $G\times_{G_{x_i}}\!U_i \ni (g,u)\mapsto
exp_{gx_i}(gu)\in M$ is a $G$ equivariant local diffeomorphism, so one gets 
a local $G$ equivariant diffeomorphism $\Psi$ between a neighborhood 
of $(e,u_1)$ in $G\times_{G_{x_1}}\!U_1$ and a neighborhood of 
$(g_2,u_2)$ in $G\times_{G_{x_2}}\!U_2$. The tangent space to $G\times_{G_
{x_i}}\!U_i$ at $(g_i,u_i)$ is equal to $T_{g_i}(G)\oplus T_{u_i}(U_i)$
and the derivative of $\Psi$ will have a block decomposition corresponding 
to this direct sum as $T_{(g_1,u_1)}(\Psi)=\left[\begin{smallmatrix}
A&B\\C&D \end{smallmatrix}\right]$. Because $G$ acts by left translations 
on the first component of  $G\times_{G_{x_i}}\!U_i$ and $\Psi$ is $G$ 
equivariant, a straightforward computation shows that $A=Id$ and $C=0$, 
so $D$ is a diffeomorphism between $T_{u_1}(U_1)$ and $T_{u_2}(U_2)$.
Then there exist neighborhoods $U'_i\subset U_i$ of $u_i$ such that 
the map $U'_1\ni u \mapsto pr_2\circ\Psi(e,u) \in U'_2$ is a 
diffeomorphism, whose derivative at $u_1$ is $D$. Denote this map by $h$.
The horizontal dotted lines in the following diagram represent locally defined 
maps that are local diffeomorphisms:
\begin{equation}
\xymatrix{
&G\times_{G_1}\!U_1 \ar@{-->}[rr]^{\Psi}_{\sim}\ar[rd]_{exp\circ\Phi_1}
&&G\times_{G_2}\!U_2 \ar[ld]^{exp\circ\Phi_2}\ar@{-->}[rddd]^{pr_2}&\\ 
&&M\ar@{->>}[d]^{proj}&&\\
&&M/G&&\\
U_1\ar[uuur]^{\iota_1} \ar[urr]_{\pi_1} \ar@{-->}[rrrr]^{h}_{\sim}&&&&
U_2\ar[llu]_{\pi_2}
}
\end{equation}

The map $pr_2$ is locally defined in a neighborhood of $(g_2,u_2)\in G\times_{G_2}\!U_2$
as $pr_2(g,u)=u$.
Using the fact that the upper triangle, the left and right quadrilaterals and 
the large trapezoid are commutative, we get that
$\pi_1=\pi_2\circ h$ on a neighborhood near $u_1$. Because $u_1$ was chosen 
arbitrarily, we conclude that the orbifold charts $\mathcal{R}_1$ and $\mathcal
{R}_2$ are compatible. 

\end{proof}

Later we will show that any connected orbifold can be obtained as the orbit space
of a smooth manifold endowed with a smooth action of a compact Lie group. 

\begin{remark}
If $M$ is an orbifold and $x\in M$ consider an orbifold chart 
$(\Tilde{U}, \Gamma, U, \pi)$ in a neighborhood of $x$ 
and $\Tilde{x}\in \Tilde{U}$ such that $\pi(\Tilde{x})=x$. The isomorphism class
of the isotropy group $\Gamma_{\Tilde{x}}$ depends only on $x$ and not
on a particular choice of $\Tilde{x}$ or chart around $x$.
We will denote this isomorphism class by $\mathcal{G}_x$.
\end{remark}

\begin{defn} A point $x\in M$ is called smooth if $\mathcal{G}_x$ is the isomorphism
class of the trivial group $(1)$ and it is called singular otherwise. 
\end{defn}

We will denote by $M_{reg}$ the set of regular points of $M$ and by $M_{sing}$ the set of 
singular points. $M_{reg}$ is an open and dense subset of $M$ whose induced orbifold 
structure is a genuine manifold structure.

\begin{defn} Let $M$ and $N$ be two orbifolds. A map $\phi:M\to N$ is 
an orbifold diffeomorphism if $\phi$ is a homeomorphism and for any
orbifold chart $\mathcal{R}=(\Tilde{U}, \Gamma, U, \pi)$ of $M$, the 
collection  $\phi(\mathcal{R})=(\Tilde{U}, \Gamma,\phi(U),
\phi\circ\pi)$ is an 
orbifold chart of $N$ and for any orbifold chart $\mathcal{R}'=(\Tilde{U}',
\Gamma', U', \pi')$ of $N$, $\phi^{-1}(\mathcal{R}')=(\Tilde{U'}, 
\Gamma', \phi^{-1}(U'), \phi^{-1}\pi')$ is an orbifold chart of $M$.
\end{defn}

To see whether a homeomorphism $\phi$ is a diffeomorphisms one must
check that $\phi$ and $\phi^{-1}$ take orbifold charts of an atlas (not
necessarily the {\it maximal} atlas) into orbifold charts.

\begin{example} If $M$ is an orbifold and $(\Tilde{U}, \Gamma, U, \pi)$ is
an orbifold  chart, then $\Tilde{U}/\Gamma$ has a canonical orbifold
structure as described in Example \ref{cano} and $U$ has an induced orbifold
structure from the orbifold structure of $M$. Then the induced map
$\overline{\pi}:\Tilde{U}/\Gamma \to U$ is an orbifold diffeomorphism. 
\end{example}

\begin{defn} Let $M$ and $N$ be two orbifolds. A continuous map $f:M\to N$ is
smooth if for any $x \in M$ one can find orbifold charts $\mathcal{R}=
(\Tilde{U}, \Gamma, U, \pi)$ around $x$ and $\mathcal{R}'=(\Tilde{U}',
\Gamma', U',\pi')$ around $f(x)$ and a smooth map $\Tilde{f}:\Tilde{U}\to
\Tilde{U}'$ such that the following diagram is commutative:
\begin{equation}
\xymatrix{
\Tilde{U}\ar[r]^{\Tilde{f}}\ar[d]_{\pi}&\Tilde{U}'\ar[d]^{\pi'}\\
U\ar[r]^{f} &U'
}
\end{equation}
\end{defn}

\begin{remark} The spaces $\mathbb{F}^k$ with $\mathbb{F}=\RR$ or
$\CC$ have a canonical structure of an orbifold, given by  the atlas 
with a unique chart $\mathcal{R}=(\mathbb{F}^k, \Gamma=(e), \mathbb
{F}^k, \pi=Id)$. Then a continuous map $f:M\to \mathbb{F}^k$  is called 
smooth  if it is smooth as a map between orbifolds.
\end{remark}

\begin{remark}
\label{glue}
If $M_i$, $i=1,2$ are two orbifolds, $U_i\subset M_i$ are open subsets
and $\phi:U_1\overset{\sim}\to U_2$ is an orbifold diffeomorphism, then
the topological space $M_1\amalg M_2/\sim$ where we identify $x\in U_1$ 
with $\phi(x)\in U_2$ has a canonical structure of an orbifold, and 
the inclusion maps $\iota_i:M_i\to M_1\amalg M_2/\sim$ are smooth.
An atlas for $M_1\amalg M_2/\sim$ is given by the reunion of two atlases
for respectively $M_1$ and $M_2$.

This procedure allows us to create an orbifold by gluing orbifold charts
along orbifold diffeomorphisms.
\end{remark}

We will need the following statement, whose proof can be found in \cite{GN1},
Theorem 2.1.

\begin{prop} 
Let $M$ be a smooth orbifold and $\{U_{\alpha}\}$ an open cover
of $M$. Then there exists a countable partition of unity $\{h_i,\;i\in
\mathbb{N}\}$ subordinated to $\{U_{\alpha}\}$ such that each $h_i$ is a 
smooth function with compact support.
\end{prop}

\subsection{Vector Orbibundles}

Let $p:E\to M$ be a continuous map between two orbifolds. Let $V$ be a fixed 
vector space over $\RR$ or $\CC$.

\begin{defn} A vector orbibundle chart is a collection 
$\mathcal{R}=(\Tilde{U}, V, \Gamma,U, \Pi, \pi)$ 
such that $(\Tilde U, \Gamma, U, \pi)$ is an orbifold chart for $M$,  
$(\Tilde U\times V, \Gamma, p^{-1}(U), \Pi)$ 
is an orbifold chart for $E$ such that the induced action of $\Gamma$ on 
the trivial vector bundle $(\Tilde U\times V\overset{pr_1}\to \Tilde{U})$ 
is by vector bundle isomorphisms,  and the following diagram is commutative:
$$
\xymatrix{
\Tilde{U}\times V \ar[r]^{\Pi} \ar[d]_{pr_1}& p^{-1}(U)
\ar[d]^{p}\\
\Tilde{U} \ar[r]^{\pi}& U}
$$
\end{defn}

Though every vector orbibundle chart comes with the actions of a group $\Gamma$ on
$\Tilde{U}$ and $\Tilde{U}\times V$, for the sake of simplicity we will not 
show them explicitly in the notation of the chart.

\begin{defn}
A vector orbibundle chart $\mathcal{R}=(\Tilde{U}, V, \Gamma, U, \Pi, \pi)$
is called a linear chart at $x\in M$ if $(\Tilde U, \Gamma, U, \pi)$ is a 
linear orbifold chart at $x\in M$ with $\mu:\Gamma\times \Tilde{U}\to \Tilde{U}$
the restriction of a linear representation of $\Gamma$ and there exists a linear 
representation $\rho:\Gamma\times V\to V$ such that the action of $\Gamma$ on the trivial
bundle $\Tilde{U}\times V \overset{pr_1}\to \Tilde{U}$ is equal to the 
diagonal action $\mu\otimes\rho$.
\end{defn}

\begin{defn} Two vector orbibundle charts $\mathcal{R}_i=
(\Tilde U_i, V, \Gamma_i, U_i, \Pi_i, \pi_i)$, 
$i=1,2$,  are compatible if for any two points $\tilde x_i \in \Tilde U_i$ 
such that $\pi(\tilde x_1)=\pi(\tilde x_2)=x \in U_1\cap U_2$ there exists a 
neighborhood $\Tilde W_1$ of $\tilde x_1$ in $\Tilde U_1$, a neighborhood 
$\Tilde W_2$ of $\tilde x_2$ in $\Tilde U_2$ and a vector bundle 
diffeomorphism $(H,h)$ between $(\Tilde W_1\times V \overset{pr_1}\to 
\Tilde W_1)$ and 
$(\Tilde W_2\times V \overset{pr_1}\to \Tilde W_2)$ such that 
$\pi_2\circ h=\pi_1$ and $\Pi_2 \circ H = \Pi_1$.
\end{defn}

\begin{defn} A vector orbibundle atlas on $E\overset{p}\to M$ is a collection
of compatible vector orbibundle charts such that the corresponding coordinate
neighborhoods of $M$ form an open cover of $M$.
A structure of vector orbibundle on $E\overset{p}\to M$ is given by a maximal
vector orbibundle atlas.
\end{defn}

Using the ideas in the proof of Proposition \ref{local} and Remark 
\ref{linchart}, one can prove that for any vector orbibundle chart 
$\mathcal{R}= (\Tilde{U}, V, \Gamma, U, \Pi, \pi)$ and $x\in U$ there exist a 
linear vector orbibundle chart at $x$, $\mathcal{R}_x$, which is compatible 
with $\mathcal{R}$. As a consequence, any vector orbibundle has an atlas 
consisting of linear vector orbibundle charts.

\begin{remark} A vector orbibundle $E\overset{p}\to M$ is usually not a 
vector bundle.
If $x\in M_{\text{sing}}$ is a singular point, then $p^{-1}(x)$  
might not have a vector space structure. If $(\Tilde U, V, \Gamma,
U, \Pi, \pi)$ is a chart around $x$ then $p^{-1}(x)$ is the 
quotient of $V$ by the action of the isotropy group $\Gamma_{\Tilde{x}}$ 
with $\pi(\Tilde{x})=x$. The isomorphism 
class of the representation of $\Gamma_{\Tilde{x}} \in \mathcal{G}_x$ on $V$ depends
only on $x$ and not on a particular choice of vector orbibundle chart and $\Tilde{x}$.
Denote this isomorphism class by $\mathcal{V}_x$. Then the restriction of 
$E\overset{p}\to M$ to $\{x\in M\,|\, \mathcal{V}_x \text{ is trivial }\}$ is a 
genuine vector bundle.
\end{remark}

\begin{example} 
\label{vcano}
Let $E\overset{p}\to M$ be a smooth vector bundle
and suppose the finite group $\Gamma$ acts on the vector bundle by bundle
diffeomorphisms. Denote by $\Pi$, resp. $\pi$, the canonical projections 
onto the orbit spaces $E\overset{\Pi}\to E/\Gamma$ and $M\overset{\pi} \to
M/\Gamma$. For $x\in M$, let $\overline{x}=\Gamma x\in M/\Gamma$ be its orbit and 
$\Gamma_x$ the isotropy group of $x$.  Choose $\Tilde{U}$ a $\Gamma_x$ invariant 
neighborhood of  $x$ in $M$, as we did in the Example \ref{cano}, such 
that $\Tilde{U}\cap\gamma\Tilde{U}=\emptyset$ for $\gamma\in\Gamma\backslash
\Gamma_x$. The group $\Gamma_x$ acts by diffeomorphisms on the restriction 
of the initial bundle $E_{|\Tilde{U}}=p^{-1}(\Tilde{U})\overset{p}\to\Tilde{U}$
and Proposition \ref{local} provides us with a linear vector orbibundle chart at 
$x$, $\mathcal{R}_x=(O, E_x, \Gamma_x, \Tilde{U}/\Gamma_x, \Pi,\pi)$.
\end{example}

We also have a generalization of the above example, analogous to Proposition 
\ref{liecano}

\begin{prop}
\label{lievcano}
Let $E\overset{p}\to M$ be a smooth vector bundle endowed with a
smooth action of a compact Lie group $G$ such that any $x\in M$ has a finite 
isotropy group $G_x\subset G$. Then 
$E/G \overset{\overline{p}}\to M/G$ has a canonical structure of a vector
orbibundle.
\end{prop}

\begin{proof}
 Let us fix a $G$ invariant Riemannian metric on the base space and a
$G$ invariant linear connection $\nabla$ in the vector bundle.

For a fixed $\overline{x}=Gx\in M/G$ with $x\in M$, consider, as described in 
the proof of Proposition \ref{liecano}, the isotropy group $G_x$ and the direct sum
decomposition of the tangent space at $x$ as $G_x$ modules $V\oplus T_x(Gx)=
T_x(M)$. Denote by $\mu_x:G_x\times V\to V$ the action of $G_x$ on the vector 
space $V$. Let $\mathcal{V}\to Gx$ be the restriction to $Gx$ of the subbundle 
of the tangent bundle to $M$
whose fiber above $gx$ is $gV$. The map $\Phi:G\times_{G_x}\!V\to \mathcal{V}$
, $\Phi(g,v)=gv$ is $G$ equivariant and realizes an isomorphism between the 
vector bundles $G\times_{G_x}\!V \to G/G_x$ and $\mathcal{V}\to Gx$. Let 
$U\subset V$ be an open $G_x$ invariant neighborhood of the 
origin such that the map $G\times_{G_x}\! U\ni(g,u) \mapsto exp\circ\Phi(g,u)=
exp_{gx}(gu)\in M$ is a $G$ equivariant map and realizes a diffeomorphism onto 
its image $N$. We 
showed in the proof of Proposition \ref{liecano} that $(U, G_x, N/G, \pi)$ is a
linear orbifold chart at $\overline{x}$
for $M/G$. We will construct a vector orbibundle chart for $E/G \overset
{\overline{p}}\to M/G$ at $\overline{x}$.

Let $\Tilde{\mathcal{V}}\overset{\Tilde{p}}\to \mathcal{V}$ be the 
pull-back of the vector bundle $E_{|Gx}\overset{p}\to Gx$ via the natural 
projection map $T(M)_{|Gx}\overset{proj}\longrightarrow Gx$ restricted to
$\mathcal{V}$.
The vector bundle 
$\Tilde{\mathcal{V}}\overset{\Tilde{p}}\to\mathcal{V}$ has a  
$G$ vector bundle structure coming from the action of $G$ on $E$ and $M$; 
indeed \linebreak
$\Tilde{\mathcal{V}}=\{(v,Y)| v\in \mathcal{V},Y\in E_{|Gx}\text{ with }
proj(v)=p(Y)\in Gx\}$ and $g\in G$ acts by the diagonal action
$g\cdot(v,Y)=(gv,gY)$. 
 
The group $G_x$ fixes the point $x$, so it acts on the fiber above $x$ by
a linear representation $\rho_x:G_x\times E_x\to E_x$. Let $G_x$ act on
$V\times E_x$ by the product action $\mu_x\otimes\rho_x$. 
Consider the
vector bundle $G\times_{G_x}\!(V\times E_x)\overset{pr_1}\to G\times_{G_x}\!
V$ with fiber $E_x$. The group $G$ acts on this bundle by left translations.  
We already showed that $\Phi:G\times_{G_x}\!V\to\mathcal{V}$ given by 
$\Phi(g,v)=gv \in \mathcal{V}_{gx}$ is a $G$ equivariant diffeomorphism.
 
Let $\Tilde{\Phi}:G\times_{G_x}\!(V\times E_x)\to\Tilde{\mathcal{V}}$ 
defined as $\Tilde{\Phi}(g,v,Y)=(gv,gY)$. Then we have the following 
commutative diagram in which the horizontal maps are $G$ equivariant 
diffeomorphisms
\begin{equation}
\xymatrix{
G\times_{G_x}\!(V\times E_x)\ar[rr]^{\Tilde{\Phi}}_{\sim}\ar[d]_{pr_1} &&
\Tilde{\mathcal{V}}\ar[d]^{\Tilde{p}}\\
G\times_{G_x}\!V\ar[rr]^{\Phi}_{\sim} && \mathcal{V}}
\end{equation}

Indeed $\Tilde{p}\circ\Tilde{\Phi}(g,v,Y)=\Tilde{p}(gv,gY)=gv=
\Phi\circ pr_1(g,v,Y)=\Phi(g,v)$. $\Tilde{\Phi}$ is surjective because 
any $(v',Y')\in \Tilde{\mathcal{V}}$ with $proj(v')=p(Y)=gx\in Gx$
is of the form $\Tilde{\Phi}(g, g^{-1}v', g^{-1}Y')$. Also, if
$\Tilde{\Phi}(g,v,Y)=\Tilde{\Phi}(g',v',Y')$ then $proj(gv)=proj(g'v')
\in Gx$, so $g^{-1}g'=g''\in G_x$ and $g'=gg''$. But $gv=g'v'=gg''v'$ and
$gY=g'Y'=gg''Y'$ so $v=g''v'$, $Y=g''Y'$ with $g''\in G_x$, $gg''=g'$. 
We conclude  that $(g,v,Y)=(g',v',Y')\in G\times_{G_x}(V\times E_x)$, so 
$\Tilde{\Phi}$ is injective as well. The pair $(\Tilde{\Phi}, \Phi)$ 
defines a $G$ equivariant isomorphism of $G$ vector bundles.

In the proof of Proposition \ref{liecano} we considered the exponential map
$exp:\mathcal{V}\to M$ with respect to the $G$ invariant metric on $M$, 
which realizes a $G$ equivariant diffeomorphism
between a neighborhood $\mathcal{U}$ of the zero section in $\mathcal{V}\to
Gx$ and $N$--a tubular neighborhood of $Gx$ in $M$. 
For $v\in gV\subset T_{gx}(M)$ and $X\in E_{gx}$, let $s(t)$, $t\in [0,1]$
be the path that realizes the parallel transport in $E\overset{p}\to M$
with respect to the $G$ invariant connection $\nabla$ above the path
$exp_{gx}(tv)\in M$, with $s(0)=X$. Then we define $\widetilde{exp}(v,X)=s(1)$,
$\widetilde{exp}:\Tilde{\mathcal{V}}\to E$. The map $\widetilde{exp}$ is
$G$ equivariant and $\widetilde{exp}(v,\cdot)$ is a linear isomorphism.
Then the restriction of $\widetilde{exp}$ to ${\Tilde{p}}^{-1}(\mathcal{U})
\subset \Tilde{\mathcal{V}}$ together with the restriction of $exp$ to
$\mathcal{U}$ give us a $G$ equivariant isomorphism between the restriction 
of the bundle $\Tilde{\mathcal{V}}\overset{\Tilde{p}}\to\mathcal{V}$ to 
$\mathcal{U}$ and the $G$ vector bundle $E_{|N}\overset{p}\to N$. 
If we choose a small enough $G_x$ equivariant neighborhood $U$ of the 
origin in $V$ as in Proposition \ref{liecano}, the pair of maps $(\widetilde{exp}, 
exp)\circ (\Tilde{\Phi}, \Phi)$ gives us a $G$ equivariant isomorphism
between the $G$ vector bundles $G\times_{G_x}\!(U\times E_x)\overset{pr_1}
\to G\times_{G_x}\!U$ and $E_{|N}\overset{p} \to N$. 
Passing to the $G$ orbits we get a homeomorphism between $G\times_{G_x}\!(
U\times E_x)/G=(U\times E_x)/G_x\overset{\overline{pr}_1}\to
G\times_{G_x}\!U/G=U/G_x$ and $E_{|N}/G=\overline{p}^{-1}(N/G)
\overset{\overline{p}}\to N/G$. 

We will describe a linear vector orbibundle chart around $\overline{x}=
Gx\in M/G$. Let $\iota:U\to G\times_{G_x}U$, $\iota(u)=(e,u)$ and $\pi$ be 
the composition $U\overset{\iota}\to G\times_{G_x}\!U\overset{exp\circ \Phi}
\longrightarrow N \to N/G$. Also let $I:U\times E_x\to G\times_{G_x}
(U\times E_x)$, $I(u,X)=(e,u,X)$ and $\Pi$ be the composition
$U\times E_x\overset{I}\to G\times_{G_x}\!(U\times E_x) \overset{\widetilde{exp}
\circ \Tilde{\Phi}}\longrightarrow E_{|N}\overset{proj}\longrightarrow E_{|N}/G$.
Then $(U, E_x, G_x, N/G, \Pi, \pi)$ is a linear vector orbibundle chart around
$Gx$. 

As in Proposition \ref{liecano}, we have to prove that any two vector 
orbibundle 
charts defined above are compatible.

Let  $\mathcal{R}_i=(U_i, E_{x_i}, G_{x_i}, N_i/G, \Pi_i, \pi_i)$, $i=1,2$ 
be two 
vector orbibundle charts around $\overline{x}_i$. Let $u_i\in U_i$ such that 
$\pi_1(u_1)=\pi_2(u_2)=\overline{x}\in M/G$. As described in Proposition 
\ref{liecano}, we can find $g_2\in G$ and $x\in M$ so that
$exp_{x_1}(u_1)=exp_{g_2 x_2}(g_2 u_2)=x$ and $Gx=\overline{x}$.
The maps $G\times_{G_{x_i}}\!(U_i\times E_{x_i})\ni(g,u,X)\mapsto 
\widetilde{exp}\circ
\Tilde{\Phi}(g,u,X)\in E$ for $i=1,2$ define $G$ equivariant diffeomorphisms. 
Then 
one can choose $G$ invariant neighborhoods $W_i$ of $(g_i,u_i)$ in $G\times_
{G_{x_i}}\!U_i$ so that the composition of the previous diffeomorphisms 
defines a $G$ equivariant diffeomorphism $\Tilde{\Psi}$ between 
${\overline{pr}}^{-1}_1(W_1)\subset G\times_{G_{x_1}}\!(U_1\times E_{x_1})$ 
and 
${\overline{pr}}^{-1}_1(W_2)\subset G\times_{G_{x_2}}\!(U_2\times E_{x_2})$. 
$\Tilde{\Psi}$ induces a $G$ equivariant diffeomorphism $\Psi$ between 
$W_1$ and $W_2$. Using the fact that $\Tilde{\Psi}$ is $G$ equivariant, 
one can show, as in the proof of Proposition \ref{liecano}, that the 
composition
\begin{equation}
H=(U_1\times E_{x_1}\overset{I}\to G\times_{G_{x_1}}\!(U_1\times E_{x_1})
\overset{\Tilde{\Psi}}\to G\times_{G_{x_2}}\!(U_2\times E_{x_2})
\overset{pr_2}\to U_2\times E_{x_2})
\end{equation}
is a local diffeomorphism which together with the induced local diffeomorphism
$h:U_1\to U_2$ between neighborhoods of $u_1$ in $U_1$ and $u_2$ in $U_2$
define a vector bundle diffeomorphism $(H,h)$ which realizes the compatibility 
of the charts $\mathcal{R}_1$ and $\mathcal{R}_2$ in a neighborhood of $u_1$.
\end{proof}

\begin{defn}
Let $E_i\overset{p_i}\to M_i$, $i=1,2$, be two vector orbibundles. 
The pair $(\Phi,\phi)$
defines a vector bundle diffeomorphism if $\Phi:E_1\to E_2$ and 
$\phi:M_1\to M_2$ \linebreak
are smooth diffeomorphisms, $p_2 \Phi= \phi p_1$, for any vector orbibundle
chart \linebreak
$\mathcal{R}=(\Tilde{U},V,\Gamma, \Pi, \pi)$ of 
$E_1\overset{p_1}\to M_1$  
the collection
$(\Phi,\phi) (\mathcal{R})=(\Tilde{U}, V, \Gamma, \phi(U), 
\Phi\Pi, \phi\pi)$ \linebreak 
is a vector orbibundle chart of $E_2\overset{p_2}\to M_2$ 
and for any vector orbibundle chart \linebreak 
$\mathcal{R'}=(\Tilde{U'},V,\Gamma', \Pi', \pi')
$ of $E_2\overset{p_2}\to M_2$ the collection \newline
\centerline{$(\Phi^{-1},\phi^{-1}) (\mathcal{R'})=(\Tilde{U'}, V', \Gamma', 
\phi^{-1}(U), \Phi^{-1}\Pi, \phi^{-1}\pi)$}\newline
is a vector orbibundle chart of $E_1\overset{p_1}\to M_1$.
\end{defn}

\begin{example}
If $E\overset{p}\to M$ is a vector orbibundle and 
$(\Tilde{U},V,\Gamma,U,\Pi,\pi)$ is a chart, then 
$((\Tilde{U}\times V)/\Gamma \overset{pr_1}\to \Tilde{U}/\Gamma)$ is a 
vector orbibundle which is isomorphic to the restriction to $U$ of the initial 
orbibundle $E_{|U}\overset{p}\to U$ by an isomorphism induced  by 
$(\Pi,\pi)$ to the $\Gamma$ orbit spaces.
\end{example}

\begin{defn}
Let $E_i\overset{p_i}\to M_i$, $i=1,2$ be two vector orbibundles. The pair
$(\Phi, \phi)$ defines a vector orbibundle morphism if $\Phi:E_1\to E_2$ 
and $\phi:M_1\to M_2$ are smooth, $p_2 \Phi= \phi p_1$, and for any 
$x\in M_1$ one can find vector orbibundle charts $\mathcal{R}=(\Tilde{U},
V,\Gamma,U,\Pi,\pi)$ around $x$ and $\mathcal{R}'=(\Tilde{U}',V',
\Gamma',U',\Pi',\pi')$ around $\phi(x)$ and a smooth vector bundle map 
$(\Tilde{\Phi},\Tilde{\phi})$ from $\Tilde{U}\times V \overset{pr_1}\to
\Tilde{U}$ to $\Tilde{U}'\times V' \overset{pr_1}\to \Tilde{U}'$ such 
that $(\Pi',\pi')\circ(\Tilde{\Phi},\Tilde{\phi})=(\Phi, \phi)\circ
(\Pi, \pi)$.
\end{defn}

\begin{remark} If $E_i\overset{p_i}\to M_i$, $i=1,2$ are two vector 
orbibundles and $U_i\subset M_i$ are two open sets such that there 
exists a vector orbibundle diffeomorphism $(\Phi, \phi)$ between 
$E_1{}_{|U_1}\overset{p_1}\to U_1$
and $E_2{}_{|U_1}\overset{p_2}\to U_2$ then $E_1\amalg E_2/{}_{\sim}\, 
\overset{p}\to M_1\amalg M_2/{}_{\sim}$ has a canonical structure 
of a vector orbibundle, an atlas being given by the reunion of two 
atlases $\mathcal{A}_i$ for $E_i\overset{p_1}\to M_i$, $i=1,2$.

In view of the above observation, one can apply to vector orbibundles 
the usual 
operations that were done on ordinary vector bundles: duality, direct 
sum, tensor
product, symmetric and exterior powers. If $E_i \overset{p_i}\to M_i$, 
$i\in I$,
is a collection of vector orbibundles and if $\mathcal{R}_{i,\alpha}$ 
are vector 
orbibundle charts, then we can apply the operations on the vector bundles
$\Tilde{U}_{i,\alpha}\times V_i\overset{pr_1}\to\Tilde{U}_{i,\alpha}$ 
and extend the action of the group $\Gamma_{i,\alpha}$ to these 
vector bundles. We will get new vector orbibundle charts and 
the vector orbibundles of orbit spaces associated to these new 
charts can be glued together by diffeomorphisms prescribed by the 
way the initial vector orbibundle charts were glued to yield the 
initial vector orbibundles. Because the initial spaces were Hausdorff, the 
underlying topological spaces for the resulting orbibundle will be 
Hausdorff as well.
\end{remark}

\begin{example}
Let $M$ be an orbifold. The tangent vector orbibundle $T(M)$ can be 
constructed  
starting from an atlas $\mathcal{A}=(\mathcal{R}_i)_i$. If $\mathcal{R}_i=
(\Tilde{U}_i,\Gamma_i, U_i, \pi_i)$ is an orbifold chart, with $\Tilde{U}
\subset\RR^m$, then $(\Tilde{U}_i, \RR^m, \Gamma_i, U_i, \pi'_i,\pi_i)$
will be a vector orbibundle chart for $T(M)$, where the action of $\Gamma_i$ 
on $\Tilde{U}\times\RR^m\overset{\sim}=T(\Tilde{U}_i)$ is the induced 
action to the tangent space level. The vector orbibundles $\Tilde{U}_i
\times\RR^m/\Gamma_i \overset{pr_1}\to\Tilde{U}_i/\Gamma_i$ are glued together 
along diffeomorphisms prescribed by the identifications of $\Tilde{U}_i/\Gamma
\overset{\sim}=U_i$ inside $M$ to form the tangent vector orbibundle $T(M)$. 
\end{example}

\begin{defn}
A metric on the orbifold $M$ is a smooth map $\mu:T(M)\otimes T(M)\to
\RR$ such that for each orbifold chart $(\Tilde{U},\Gamma,U,\pi)$
the map $\mu$ is induced by a $\Gamma$ invariant metric $\Tilde{\mu}:
T(\Tilde{U})\otimes T(\Tilde{U})\to \RR$.

An orbifold endowed with a metric is called a Riemannian orbifold.
\end{defn}

On any paracompact orbifold one can construct a metric by first defining 
metrics on orbifold charts and then gluing them together using a 
smooth partition of unity.

\begin{defn} 
A hermitian structure on a complex vector orbibundle $E\overset{p}\to M$
is given by a smooth map $<,>:E\otimes E \to \CC$ such that 
for each vector orbibundle chart $(\Tilde{U}, V,\Gamma, \Pi, \pi)$ there 
exists a $\Gamma$ invariant hermitian structure on $(\Tilde{U}\times 
V \overset{pr_1}\to \Tilde{U})$ which induces $<,>$ above 
$\Tilde{U}/\Gamma\overset{\sim}=U$.
\end{defn}

\begin{remark} We can think of a vector orbibundle as of a smooth family 
of linear representations of groups $\Gamma_x\in \mathcal{G}_x$ on 
$\Tilde{E}_x$ parameterized by $x\in M$. Then $E_x=\Tilde{E}_x/\Gamma_x$. 
A hermitian structure is a smooth family of hermitian products in 
$\Tilde{E}_x$ which are $\Gamma_x$ invariant and induce $<,>_x$ on 
$E_x=\Tilde{E}_x/\Gamma_x$.
\end{remark}

\begin{defn}
Let $E\overset{p}\to M$ be a vector orbibundle. A map $f:M\to E$ is a smooth 
section if $p\circ f= Id_M$ and if $f$ is smooth as a map between the orbifolds
$M$ and $E$. We will denote the space of smooth sections by $C^{\infty}(M;E)$.
\end{defn}

\begin{remark}
If $\mathcal{A}=({\mathcal{R}_i})$ is a vector orbibundle atlas with $\mathcal{R}_i
=(\Tilde{U}_i, V, \Gamma_i, U_i, \Pi_i, \pi_i)$ then the vector orbibundles 
$\Tilde{U}_i\times V/\Gamma_i\overset{pr_1}\to \Tilde{U}/\Gamma_i$ and 
$E_{|U_i}\overset{p}\to U_i$ are isomorphic via the map induced by $(\Pi_i,\pi_i)$ to 
the orbit spaces. A map $f:M\to E$ is a smooth section if the restriction to each
$U_i$ can be identified via the above vector bundle diffeomorphisms to sections 
induced to the $\Gamma_i$ orbit spaces by $\Gamma_i$ invariant smooth sections $\Tilde{f}_i$
of the $\Gamma_i$ vector bundles $\Tilde{U}_i\times V\overset{pr_1}\to \Tilde{U}_i$.
\end{remark}

In Proposition \ref{lievcano} we showed that starting with a $G$ vector bundle 
$E\overset{p}\to M$ where the Lie group $G$ acts on $M$ with finite isotropy groups 
and passing to the $G$ orbit spaces we obtain a vector orbibundle $E/G \overset
{\overline{p}}\to M/G$. In this case we have a global characterization of the 
smooth sections of the vector orbibundle. 

\begin{prop}
\label{gsection}
The canonical map 
\begin{equation}
Sect(M;E)^G\to Sect(M/G;E/G)
\end{equation}
identifies the $G$ invariant 
smooth sections of the $G$ vector bundle $E\overset{p}\to M$ with the smooth sections 
of the vector orbibundle of $G$ orbits $E/G \overset{\overline{p}}\to M/G$.
\end{prop}
\begin{proof}

It is obvious that any $G$ invariant section $f\in Sect(M;E)^G$ defines a section 
$\overline{f}\in Sect(M/G;E/G)$ by $\overline{f}(Gx)=Gf(x)$. 

Let $f_1$ and $f_2$ be two $G$ invariant sections in $Sect(M;E)^G$ which induce 
equal sections $\overline{f}_1=\overline{f}_2\in Sect(M/G;E/G)$. If $f_1\ne f_2$
then there exists $x\in M$ such that $f_1(x)\ne f_2(x)$. Nevertheless, $\overline
{f}_1(\overline{x})=\overline{f}_2(\overline{x})\in E/G$. Then there exists $g\in G$ 
such that $f_1(x)=gf_2(x)\in E_x$. But $gf_2(x)=f_2(gx)\in E_x$ so $gx=x$. So
$f_1(x)=f_1(g^{-1}x)=g^{-1}f_1(x)=g^{-1}gf_2(x)=f_2(x)$. We reached a contradiction so
the map $Sect(M:E)^G\to Sect(M/G;E/G)$ is injective.

Let $h\in C^{\infty}(M/G;E/G)$ be a smooth section. We need to construct a $G$ invariant 
smooth section $f\in C^{\infty}(M;E)^G$ such that $\overline{f}=h$.
We will construct $f$ from local data. As we showed in Proposition \ref{lievcano}, 
the local model of passing from a $G$ vector bundle 
to a vector orbibundle is given by passing from the $G$ vector bundle $G\times_{G'}
\!(U\times V) \overset{p}\to G\times_{G'}\! U$ to the vector orbibundle
$(U\times V)/G' \overset{\overline{p}}\to U/G'$. $G'\subset G$ is a finite group 
which acts on the left on the open set $U\subset \RR^n$ and on a vector 
space $V$ and by right translations on $G$. The action of $G'$ on $G\times U$ and
$G\times(U\times V)$ is given by $g' \cdot (g,u)=(g g'{}^{-1}, g'u)$ and 
$g' \cdot (g,u,v)= (gg'{}^{-1}, gu, gv)$. In fact the vector bundle $G\times(U\times V)
\overset{\Tilde{p}}\to G\times U$ is a $G'$ vector bundle whose associated vector
orbibundle of $G'$ orbits is the local model $G\times_{G'}\!(U\times V)\overset{p}
\to G\times_{G'}\!U$, which happens to be a genuine vector bundle because the action
of $G'$ is free. A section $h\in C^{\infty}(U/G';
(U\times V)/G')$ is induced by a smooth section $\Tilde{h}\in C^{\infty}(U; U\times V)$.
Let $\Tilde{f}\in C^{\infty}(G\times U; G\times (U\times V))$ be defined by 
$\Tilde{f}(g,u)=(g,h(u))$. This is a $G$ invariant smooth section which is also $G'$ 
invariant and it induces a smooth section $f\in C^{\infty}(G\times_{G'}\! U;
G\times_{G'}\!(U\times V))$. It is obvious that the map $f$ induces $h$ when passing 
to the $G$ orbit spaces. Because of the first part of the proof such an $f$ is
unique, so we can glue different sections $f$ to get a global $G$ invariant smooth
section in $C^{\infty}(M;E)$ which will induce $h\in C^{\infty}(M/G; E/G)$.

The problem of proving that if $f\in Sect(M;E)^G$ is smooth then its image
$\overline{f}\in Sect(M/G;E/G)$ is smooth is a local one. A smooth 
section $f\in C^{\infty}(G\times_{G'}\! U;G\times_{G'}\!(U\times V))$ is induced by a $G'$ 
invariant smooth section $\Tilde{f}\in C^{\infty}(G\times U; G\times(U\times V))^{G'}$. 
If $f$ is $G$ invariant then $\Tilde{f}$ is $G$ invariant so 
there exists a $G'$ smooth section $h$ of the $G'$ vector bundle $U\times 
V\overset{pr_1}\to U$ such that $\Tilde{f}(g,u)=(g,h(u))$.
The induced section $\overline{f}$ of the vector orbibundle $(U\times V)/G' \overset
{\overline{p}}\to U/G'$ is the same as the one induced by $h$. 
So $\overline{f}$ is smooth. 
\end{proof}

\begin{remark}
\label{globalize}
Let $\Tilde{E}\overset{\Tilde{p}}\to \Tilde{M}$ be a $\Gamma$ vector 
bundle and $E\overset{p}\to M$ be the associated vector orbibundle of
$\Gamma$ orbits. If $\Tilde{U}\subset\Tilde{M}$ is an open subset and
$\Gamma'\subset\Gamma$ a subgroup such that $\Tilde{U}$ is $\Gamma'$ invariant
and $\gamma\Tilde{U}\cap\Tilde{U}=\emptyset$ for $\gamma\in\Gamma
\setminus\Gamma'$, then $\Tilde{E}_{|\Tilde{U}}\overset{p}\to \Tilde{U}$ is
a $\Gamma'$ vector bundle and the associated vector orbibundle of 
$\Gamma'$ orbits is a subbundle of $E\overset{p}\to M$. Let $U=\Tilde{U}/
\Gamma'$. Denote by
$\iota$ the inclusion map between the two vector orbibundles.
$\iota$ induces a map between the spaces of compactly supported sections 
$\iota_*: C^{\infty}_0(U;E_{|U})\to C^{\infty}_0(M;E)$ given by extension 
by zero on $M\setminus U$ and a restriction map $\iota^* :C^{\infty}
(M;E)\to C^{\infty}(U;E_{|U})$. We have $C^{\infty}(M;E)\overset{\sim}=
C^{\infty}(\Tilde{M};\Tilde{E})^{\Gamma}$ and $C^{\infty}(U;E_{|U})\overset
{\sim}=C^{\infty}(\Tilde{U};\Tilde{E}_{|\Tilde{U}})^{\Gamma'}$. 
Then the map $\iota^*$ is induced by the restriction map $C^{\infty}
(\Tilde{M};\Tilde{E})\to C^{\infty}(\Tilde{U};\Tilde{E}_{|\Tilde{U}})$.
The map $\iota_*$ between $C^{\infty}(\Tilde{U}; \Tilde{E}_{|\Tilde{U}})$ and
$C^{\infty}(\Tilde{M};\Tilde{E})$ is given by 
\begin{align}
\iota_*(\Tilde{f})&=\sum_{i=1}^l \gamma_i \Tilde{f}\\
\intertext{where $\{\gamma_1=e,\gamma_2,\dots,\gamma_l\}$ with $l=[\Gamma:\Gamma']$
is a complete system of left coset representatives for $\Gamma/\Gamma'$, or}
\iota_*(\Tilde{f})&=\frac{1}{|\Gamma'|}\sum_{\gamma\in\Gamma}\gamma \Tilde{f}
\end{align}
for $\Tilde{f}\in C^{\infty}(\Tilde{U};\Tilde{E}_{|\Tilde{U}})^{\Gamma'}$.
\end{remark}

\subsection{The Representation Theorem}

In this section we will state and prove the reciprocal of 
Proposition \ref{lievcano}.

\begin{thm}
\label{represent}
Any vector orbibundle $E\overset{p}\to M$ is isomorphic to the vector orbibundle
of $G$ orbits associated with a $G$ vector bundle $\Tilde{E}\overset{\Tilde{p}}\to 
\Tilde{M}$, with $G$ a compact Lie group acting on $M$ with finite isotropy groups. 
In particular any orbifold $M$ is diffeomorphic to the orbifold of $G$ orbits 
associated with a $G$ manifold $\Tilde{M}$ where the action of the compact Lie  group
$G$ has finite isotropy groups. 

\end{thm}

\begin{proof}
Let us fix a metric $\mu$ on the base space $M$. Let $G=O(m)$ be the 
orthogonal group, $m=dim(M)$. We will describe the orthonormal frame 
bundle $\mathcal{F}(M)\overset{proj}\longrightarrow M$.
We will show that $\mathcal{F}(M)$ is a manifold and that the group 
$G$ acts on it with finite isotropy groups. The $G$ vector bundle we need 
to construct will be essentially the pull-back of $E\overset{p}\to M$ 
to $\mathcal{F}(M)$ via the canonical projection.

The description of the space of orthonormal frames uses an atlas 
$\mathcal{A}=(\mathcal{R}_i)$ of the vector orbibundle $E\overset{p}\to M$.
Let $\Tilde{\mu}_i$ be the $\Gamma_i$ invariant metric on 
$\Tilde{U}_i\subset \RR^m$ which induces the metric $\mu$ via the map 
$\Tilde{U}_i\overset {\pi_i}\to U_i$. Let $\mathcal{F}(\Tilde{U}_i)$ be the 
space of orthonormal frames of $\Tilde{U}_i$. Then the orthogonal group
$O(m)$ acts by right translations on $\mathcal{F}(\Tilde{U}_i)$. Because 
the metric $\Tilde{\mu}_i$
is $\Gamma_i$ invariant, the action of $\Gamma_i$ on $\Tilde
{U}_i$ extends to $\mathcal{F}(\Tilde{U}_i)$. Observe that this action 
is free. Indeed, if $\gamma\in \Gamma_i$ and $(x,\mathcal{X})\in \mathcal{F}
(\Tilde{U}_i)$ such that $\gamma(x,\mathcal{X})=(x,\mathcal{X})$ then 
$\gamma\in \Gamma_x$ and $\gamma$ acts on the tangent space $T_x(\Tilde{U}_i)$
as the identity map. But the order of $\gamma$ is finite so $\gamma=
e\in \Gamma_i$. Then the space of orbits $\mathcal{F}(\Tilde{U}_i)/
\Gamma_i$ is a smooth manifold. The local diffeomorphisms between open
sets $\pi_i^{-1}(U_i\cap U_j)\subset \Tilde{U}_i$ and $\pi_j^{-1}
(U_i\cap U_j)\subset \Tilde{U}_j$ preserve the metrics so they will 
induce local diffeomorphisms between $\mathcal{F}\pi_i^{-1}(U_i\cap
U_j)$ and $\mathcal{F}\pi_i^{-1}(U_i\cap U_j)$ and a diffeomorphism 
between their images in $(\mathcal{F}(\Tilde{U}_i))/\Gamma_i$ respectively
$\mathcal{F}(\Tilde{U}_j)/\Gamma_j$. The space of orthonormal
frames $\mathcal{F}(M)$ is obtained by gluing the manifolds $(\mathcal{F}
(\Tilde{U}_i))/\Gamma_i$ along these diffeomorphisms and it has a natural
structure of smooth manifold. The right action of $O(m)$ on $\mathcal{F}
(\Tilde{U}_i)$ commutes with the action of $\Gamma_i$ and with the local
diffeomorphisms induced between different $\mathcal{F}(\Tilde{U}_i)$ 
and $\mathcal{F}(\Tilde{U}_j)$ so we will get an induced right action of 
$O(m)$ on $\mathcal{F}(M)$. 

Let us consider the $\Gamma_i$ vector bundle $\Tilde{U}_i\times V 
\overset{pr_1}\to\Tilde{U}_i$ corresponding to the chart $\mathcal{R}_i$ 
and the pull-back to the space of orthonormal frames 
$\mathcal{F}(\Tilde{U}_i)\times V \overset{pr_1}\to\mathcal{F}
(\Tilde{U}_i)$ together with the induced left action of the group 
$\Gamma_i$ on it and the right action of group $O(m)$. The action of 
$\Gamma_i$ is free and commutes with the action of $O(m)$ so if we pass 
to the spaces of $\Gamma_i$ orbits we
get a genuine vector bundle $\mathcal{F}(\Tilde{U}_i)\times V/\Gamma_i
\overset{\overline{pr}_1} \to\mathcal{F}(\Tilde{U}_i))/\Gamma_i$ endowed
with the action of $O(m)$. By gluing these vector bundles together using 
diffeomorphisms induced by the identifications made on the initial vector orbibundle
charts we get an $O(m)$ vector bundle $\Tilde{E}\overset{\Tilde{P}}\to \mathcal{F}(M)$. 

We will show that the action of $O(m)$ on $\mathcal{F}(M)$ has finite isotropy 
groups and if we pass to the space of $O(m)$ orbits we get a vector orbibundle 
isomorphic with the initial vector orbibundle $E\overset{p}\to M$. 

Any point in $\mathcal{F}(M)$ is the $\Gamma_i$ orbit of a point $(x,\mathcal{X})
\in \mathcal{F}(\Tilde{U}_i)$ for some index $i$. The group $O(m)$ acts freely 
on $\mathcal{F}(\Tilde{U}_i)$.
If $g\in O(m)$ fixes the $\Gamma_i$ orbit $\overline{(x,\mathcal{X})}$ then there exists 
$\gamma\in \Gamma_i$ such that $(x,\mathcal{X})g=\gamma(x,\mathcal{X})$. This is 
equivalent to $x=\gamma x$ and $\mathcal{X}g=\gamma\mathcal{X}$. The
group $O(m)$ acts freely and transitively on the frames at $x\in\Tilde{U}_i$.
The group $(\Gamma_i)_x$ acts freely on the frames at $x\in\Tilde{U}_i$ as well.
Then the above equalities imply that the cardinality of the isotropy group
of $\overline{(x,\mathcal{X})}\in\mathcal{F}(M)$ with respect to the action of 
$O(m)$ is equal to the cardinality of the group $(\Gamma_i)_x$, which is finite.

Passing to the vector orbibundle of $O(m)$ orbit spaces $\Tilde{E}/O(m)\overset
{\overline{\Tilde{p}}}\to \mathcal{F}(M)/O(m)$ can be realized by first passing to the 
$\Gamma_i$ vector bundles $\mathcal{F}(\Tilde{U})\times V/O(m)\overset{\overline{pr}_1}
\to \mathcal{F}(\Tilde{U})/O(m)$ and then gluing together the resulting vector 
orbibundles of $\Gamma_i$ orbits. But the $\Gamma_i$ vector bundles 
$\Tilde{U}\times V\overset{pr_1}\to\Tilde{U}$ and $\mathcal{F}(\Tilde{U})\times V/O(m)
\overset{\overline{pr}_1}\to \mathcal{F}(\Tilde{U})/O(m)$ are canonically isomorphic
and the resulting vector orbibundles obtained by gluing the vector orbibundles of
$\Gamma_i$ orbits will be isomorphic. So $E \overset{p}\to M$ and $\Tilde{E}/O(m)
\overset{\Tilde{p}}\to \mathcal{F}(M)/O(m)$ are isomorphic.

\end{proof}

\subsection{Sobolev Spaces}
Consider the Euclidean space $\RR^m$ with the canonical scalar product 
$\left<x,y\right>=\sum_{i=1}^m x_iy_i$ and let the finite group $\Gamma$ act on $\RR^m$ by 
linear isometries. We have an induced action of $\Gamma$ on the space of function 
on $\RR^m$ defined by $\gamma f(x)=f(\gamma^{-1}x)$.

For compactly supported smooth functions $f\in C^{\infty}_0(\RR^m)$ the Fourier 
transform $\Hat{f}:\RR^m\to \RR$ is defined by
$\Hat{f}(\xi)=\int e^{-i<x,\xi>}f(x) \, dx$.
Then $\Hat{(\gamma f)}(\xi)=\int e^{-i<x,\xi>}f(\gamma^{-1}x)\,dx$ and if we replace 
$x$ by $\gamma y$ we get
\begin{equation}
\begin{split}
\Hat{(\gamma f)}(\xi)&=\int e^{-i<\gamma y,\xi>}f(y)det(\gamma)\,dy=\\
&=\int e^{-i<y, \gamma^{-1} \xi>}f(y)\, dy=
\Hat{f}(\gamma^{-1}\xi)=\gamma \Hat{f}(\xi)
\end{split}
\end{equation}
(we used the fact that $\gamma$ was an isometry and $det(\gamma)=1$).

This shows that the Fourier transform commutes with the action of $\Gamma$ on the space 
of functions.

For $s\in\RR$ and $f\in C^{\infty}(\RR^m)$ the Sobolev norm $||\cdot||_s$
is defined by $||f||^2_s=\int (1+|\xi|^2)^s |\Hat{f}(\xi)|^2\,d\xi$. The Sobolev space $H_s
(\RR^m)$ is the completion of $C^{\infty}_0(\RR^m)$ with respect to the norm
$||\cdot||_s$. Because  the group $\Gamma$ acts by isometries on $\RR^m$ and 
the Fourier transform is $\Gamma$ equivariant  we conclude that the action of $\Gamma$ on 
the space of smooth functions preserves the Sobolev norms and extends to an action by 
isometries on the Sobolev spaces $H_s(\RR^m)$. The closure of the space of smooth
$\Gamma$ invariant functions $C^{\infty}_0(\RR^m)^{\Gamma}$ coincides with the
subspace of $\Gamma$ invariant elements $H_s(\RR^m)^{\Gamma}=\{f\in H_s(\RR^m)\,
|\, \gamma f=f \text{ for all } \gamma \in \Gamma \}$.

In order to define the Sobolev norm of the smooth functions on the orbifold $\RR/
\Gamma$ it is convenient to replace the Sobolev norm of the $\Gamma$ invariant 
smooth functions with an equivalent norm. 

\begin{defn}
For $f \in C^{\infty}_0(\RR^m)^{\Gamma}\overset{\sim}=
C^{\infty}_0(\RR^m/\Gamma)$ let $||f||^{\Gamma}_s=\frac{1}{|\Gamma|}||f||_s$.

The Sobolev space $H_s(\RR^m/\Gamma)$ is the completion of $C^{\infty}
(\RR^m/\Gamma)$ with respect to the norm $||\cdot||^{\Gamma}_s$.

If $\Tilde{U}\subset \RR^m$ is a $\Gamma$ invariant neighborhood of the origin 
with compact closure 
let $H_s(\Tilde{U}/\Gamma)$ be the closure of $C^{\infty}_0(\Tilde{U})^{\Gamma}$ 
in $H_s(\RR^m/\Gamma)$.
\end{defn}

We observe that $H_s(\RR^m/\Gamma)\overset{\sim}=H_s(\RR^m)^{\Gamma}$.

Our goal is to define the Sobolev spaces $H_s(M)$ for $M$ a compact orbifold. The above 
definitions will describe the elements of the Sobolev space $H_s(M)$ whose support is
included in the image of a linear chart $(\Tilde{U}, \Gamma, \Tilde{U}/\Gamma\overset
{\sim}=\pi(\Tilde{U}), \pi)$. In order for this description to be independent
of the chosen orbifold chart we need to show that the definition of the spaces
$H_s(\Tilde{U}/\Gamma)$ is compatible with the following two operations occurring 
in the compatibility of the linear orbifold charts:
\begin{enumerate}

\item If $\phi:\Tilde{U}_1\to \Tilde{U}_2$ is a diffeomorphism intertwining the 
two linear actions of the group $\Gamma$ on two open neighborhoods of the origin 
$\Tilde{U}_i$, $i=1,2$ in $\RR^m$ then the two orbifold charts 
$(\Tilde{U}_i, \Gamma, \Tilde{U}_i
/\Gamma, \pi_i)$, $i=1,2$ are compatible with the map $\phi$ inducing a diffeomorphism
between the orbifolds $\Tilde{U}_1/\Gamma$ and $\Tilde{U}_2/\Gamma$.

\item If $\Gamma_x\subset \Gamma$ is the isotropy group of $x\in \RR^m$
and $\Tilde{V}\subset \Tilde{U}$ is a $\Gamma_x$ invariant neighborhood of $x$ such that
$\gamma \Tilde{V}\cap \Tilde{V}=\emptyset$ for $\gamma \in \Gamma \backslash \Gamma_x$ then 
the restricted chart $(\Tilde{V}, \Gamma_x, \pi(\Tilde{V}), \pi_{|V})$ and the initial chart
$(\Tilde{U}, \Gamma, \Tilde{U}/\Gamma, \pi)$ are compatible via the canonical 
inclusion map $\iota:\Tilde{V}\hookrightarrow \Tilde{U}$ inducing a smooth map between
the associated orbifolds $\Tilde{V}/\Gamma_x$ and $\Tilde{U}/\Gamma$.

\end{enumerate}

In the first case the diffeomorphism $\phi$ induces a continuous bijection with a
continuous inverse between the spaces of smooth functions $C^{\infty}_0(\Tilde{U}_i)$, 
for $i=1,2$, endowed with the Sobolev norms, (cf. \cite{gilkey} Lemma [1.3.3]) and permutes 
with the action
of the group $\Gamma$. As a consequence $\phi$ will induce a continuous bijection
with a continuous inverse between the Sobolev spaces $H_s(\Tilde{U}_i/\Gamma)\overset
{\sim}=H_s(\Tilde{U}_i)^{\Gamma}$, $i=1,2$.

In the second case the inclusion $\iota:\Tilde{V}\to \Tilde{U}$ induces the map
$\iota_*:C^{\infty}_0(\Tilde{V})^{\Gamma_x}\to C^{\infty}_0(\Tilde{U})^{\Gamma}$ with
$\iota_*(f)=\sum_{i=1}^k\gamma_i f$ where $\{\gamma_i\}_{i=1}^k$ is a complete system 
of left coset representatives for $\Gamma/\Gamma_x$. We extend the compactly 
supported function $f$ by zero outside $\Tilde{V}$. Then $||\iota_*(f)||_s^{\Gamma}
=\frac{1}{|\Gamma|}||\iota_*(f)||_s=\frac{1}{[\Gamma:\Gamma_x]|\Gamma_x|}
\sum_{i=1}^k ||\gamma_i f||_s=\frac{1}{|\Gamma_x|}||f||_s=||f||_s^{\Gamma_x}$.
So $\iota_*$ induces an isometric embedding of $H_s(\Tilde{V}/\Gamma_x)$ into
$H_s(\Tilde{U}/\Gamma)$.

We will need the following statement as well:

\begin{lemma}
If $f\in C^{\infty}_0(\Tilde{U})^{\Gamma}\overset{\sim}=C^{\infty}_0(\Tilde{U}/\Gamma)$ 
is a compactly supported smooth function then the multiplication operator by $f$ is
continuous operator on $C^{\infty}_0(\Tilde{U}/\Gamma)$ endowed with the Sobolev norm
$||\cdot||_s$.
\end{lemma}

\begin{proof}
The multiplication operator by $f$ is a continuous map on $C^{\infty}_0(\Tilde{U})$
endowed with the Sobolev norm $||\cdot||_s$, being a differential operator of order 
zero. The the restriction to the $\Gamma$
invariant smooth functions remains a continuous operator.
\end{proof}

We will define the Sobolev spaces $H_s(M)$ for a compact orbifold $M$. 

Let $\mathcal{A}=(\Tilde{U}_i, \Gamma_i, U_i, \pi_i)_i$ be a finite atlas of the orbifold $M$
consisting of linear charts and $(\phi_i)_i$ be a smooth partition of unity associated with 
the open cover $(U_i)_i$. For $f\in C^{\infty}(M)$ define the Sobolev norm 
\begin{equation}
||f||^{\mathcal{A},(\phi_i)}_s=\sum_i||\pi_i^*(\phi_i f)||^{\Gamma_i}_s
\end{equation}
where $\pi_i^*(\phi_i f)\in C^{\infty}_0(\Tilde{U}_i)^{\Gamma_i}$ is 
the unique $\Gamma_i$ 
invariant map inducing $\phi_i f\in C^{\infty}(U_i)$. As in the case of 
manifolds, this
norm depends on the choice of the atlas and the partition of unity, but 
two different choices 
lead to two equivalent norms. This is a direct consequence of the 
equivalence of Sobolev norms for equivalent orbifold charts as described above
and the continuity of the multiplication operators by smooth functions. 
The Sobolev space $H_s(M)$ is the completion of the space of smooth functions
with respect to any of the norms $||\cdot||^{\mathcal{A},(\phi_i)}_s$.

In the case of vector orbibundles, we define the Sobolev norms 
$||\cdot||^{\Gamma}_s$
for the local model $(\Tilde{U}\times \RR^k)/\Gamma \overset{pr_1}
\to \Tilde{U}/\Gamma$, with$\Tilde{U}\subset \RR^m$,
by $||f||^{\Gamma}_s=\frac{1}{|\Gamma|}||f||_s$ where $f\in C^{\infty}_0
(\Tilde{U}/\Gamma; (\Tilde{U}\times \RR^k)/\Gamma))\overset{\sim}=
C^{\infty}_0(\Tilde{U};\Tilde{U}\times
\RR^k)^{\Gamma}$ and $||\cdot||_s$ is the usual Sobolev norm restricted 
to the $\Gamma$ invariant sections. If $E\overset{p}\to M$ is a 
vector orbibundle, we use a vector orbifold atlas and a partition of 
unity to define the Sobolev norm of smooth sections as in the case of 
smooth functions. All Sobolev norms defined using different atlases 
and partitions of unity 
will be equivalent and the Sobolev space $H_s(M;E)$ is the 
completion of $C^{\infty}(M;E)$  
with respect to any of these equivalent norms.

\begin{remark}
We showed in Proposition \ref{represent} that any vector orbibundle 
$E\overset{p}\to M$ is isomorphic to the associated vector orbibundle 
of a $G$ vector bundle $\Tilde{E} \overset{\Tilde{p}}\to \Tilde{M}$ 
with the compact Lie group $G$ acting on $\Tilde{M}$
with finite isotropy groups. Then $C^{\infty}(M;E)\overset{\sim}=
C^{\infty}(\Tilde{M};\Tilde{E})^G$. Also 
$H_s(M;E)=H_s(\Tilde{M};\Tilde{E})^G$. To prove this, we use 
a $G$ invariant Sobolev norm $||\cdot||_s$ on $C^{\infty}(\Tilde{M};\Tilde{E})$. 
The fact that the restriction of the Sobolev norm $||\cdot||_s$ to 
the $G$ invariant smooth sections is equivalent to the Sobolev norm on 
$C^{\infty}(M;E)$ can be seen directly when 
looking at the local model $G\times_{G'}\!(U\times V) \overset{p}
\to G\times_{G'}\! U$ 
and passing to the vector orbibundle $(U\times V)/G' \overset{\overline{p}}
\to U/G'$. Then the completion of $C^{\infty}(\Tilde{M};\Tilde{E})^G
\overset{\sim}=C^{\infty}(M;E)$ with respect to the induced Sobolev norm 
$||\cdot||_s$ is equal to $H_s(M;E)$.
\end{remark}

\subsection{Dirac-Type Densities on Orbifolds}

In this chapter we will construct a certain family of Dirac-type densities and 
the associated Dirac-type distributions. 

\begin{defn} 
Let $\Gamma\times \Tilde{M}\to \Tilde{M}$ be a faithful action of a finite 
group on a smooth manifold and $M=\Tilde{M}/\Gamma$ the associated orbifold. 
A Dirac-type density $\eta$ on the orbifold $M$ is given by a $\Gamma$ 
Dirac-type density $\{\eta^{\gamma}\}_{\gamma\in\Gamma}$ on the $\Gamma$ 
manifold $\Tilde{M}$. 

The Dirac-type density $\eta$ defines a continuous functional on $C^{\infty}_0
(M)$ by
\begin{equation}
<\eta,f>=\frac{1}{|\Gamma|}\sum_{\gamma\in\Gamma}
\int_{\Tilde{M}^{\gamma}}\Tilde{f}_{|\Tilde{M}^{\gamma}} \eta^{\gamma}
\end{equation}
where $\Tilde{f}\in C^{\infty}_0(\Tilde{M})^{\Gamma}$ induces $f\in
C^{\infty}_0(M)$. We call the above functional the Dirac-type distribution 
associated with $\eta$.
\end{defn}

We want to extend the above notions to a general orbifold $M$, which is not 
the quotient of a smooth manifold by a finite group.

\begin{remark}
Let $\Gamma_i$, $i=1,2$ be two isomorphic groups acting respectively on two 
smooth manifolds $\Tilde{M}_1$ and $\Tilde{M}_2$ as above. Suppose that 
the actions are conjugated by a diffeomorphism $\Tilde{\phi}:\Tilde{M}_1
\to\Tilde{M}_2$. Then the orbifolds $M_1=\Tilde{M}_1/\Gamma_1$ and
$M_2=\Tilde{M}_2/\Gamma_2$ are diffeomorphic by the diffeomorphism 
$\phi$ and for any Dirac-type density $\eta_1=\{\eta^{\gamma}\}
_{\gamma\in\Gamma_1}$ on $M_1$ the collection $\phi_*(\eta_1)=
\{\Tilde{\phi}_*(\eta^{\gamma})\}_{\gamma\in\Gamma_1}$ is a Dirac-type
density on $M_2$. The corresponding Dirac-type distributions on 
$C^{\infty}_0(M_i)$, $i=1,2$ are conjugated by the isomorphism 
induced by $\phi$ between the spaces of smooth functions.
\end{remark}

\begin{remark}
\label{restrict:density}
Let $\eta$ be a Dirac-type density on the orbifold $M=\Tilde{M}/\Gamma$. 
Let $\Tilde{U}\subset \Tilde{M}$ be an open subset and $\Gamma'\subset\Gamma$ 
a subgroup such that $\Tilde{U}$ is $\Gamma'$ invariant and $\gamma\Tilde{U}
\cap \Tilde{U}=\emptyset$ for $\gamma\in\Gamma\setminus\Gamma'$. Then
we have a $\Gamma$ equivariant diffeomorphism $\Gamma\times_{\Gamma'}\Tilde{U}
\overset{\sim}=\Gamma\Tilde{U}$. The orbifold $U=\Tilde{U}/\Gamma'$ is 
an open suborbifold of $M$. Denote by $\iota:U\to M$ the inclusion map.
The restriction of the Dirac-type densities 
$\eta^{\gamma'}$ on $U$ for $\gamma'\in\Gamma'$ define a $\Gamma'$ 
Dirac-type density on $\Tilde{U}$ and so a Dirac-type density on the orbifold 
$U$ which we denote by $\iota^*(\eta)$ or simply $\eta_{|U}$. 
Let $\iota_* :C^{\infty}_0(U) \to C^{\infty}_0(M)$ be the map given by the 
extension by zero on $M\setminus U$. Then 
\begin{equation}
\label{restrict:dens:f}
<\eta_{|U}, f>=<\eta,\iota_*(f)> \qquad\text{ for any } f\in C^{\infty}_0(U)
\end{equation}
Indeed, if $\Tilde{f}\in C^{\infty}_0(\Tilde{U})$ induces $f$, then 
$\iota_*(\Tilde{f})=\sum_{i=1}^{l}\gamma_i \Tilde{f}$, with $\{\gamma_1,
\gamma_2,\dots,\gamma_l\}$ a complete system of left coset representatives 
for $\Gamma/\Gamma'$, as proved in Remark \ref{globalize}.

The left-hand side of \eqref{restrict:dens:f} is equal to 
\begin{equation}
<\eta_{|U}, f>=\frac{1}{|\Gamma'|}\sum_{\gamma'\in\Gamma'}
\int_{\Tilde{U}^{\gamma'}}\Tilde{f}_{|\Tilde{U}^{\gamma'}} \eta^{\gamma'}
\end{equation}
and the right-hand side is equal to
\begin{align}
<\eta,\iota_*(\Tilde{f})>&=\frac{1}{|\Gamma|}\sum_{\gamma\in\Gamma}
\int_{\Tilde{M}^{\gamma}}
\sum_{i=1}^{l}\gamma_i\Tilde{f}_{|\Tilde{M}^{\gamma}} \eta^{\gamma}=\\
&=\frac{1}{|\Gamma|}\sum_{i=1}^{l}(\sum_{\gamma\in\Gamma}\int_{\gamma_i^{-1}
\Tilde{M}^{\gamma}}f(\gamma_i^{-1}\eta^{\gamma}))=\\
&=\frac{1}{|\Gamma|}\sum_{i=1}^{l}(\sum_{\gamma\in\Gamma}\int_{(\Tilde{M}
^{\gamma_i^{-1}\gamma\gamma_i}\cap \Tilde{U})}f\eta^{\gamma_i^{-1}\gamma
\gamma_i})=\\
&=\frac{1}{|\Gamma|}\sum_{i=1}^{l}(\sum_{\gamma_i^{-1}\gamma\gamma_i\in\Gamma'}
\int_{\Tilde{U}^{\gamma_i^{-1}\gamma\gamma_i}}f
\eta^{\gamma_i^{-1}\gamma\gamma_i})=\\
&=\frac{1}{|\Gamma|}\sum_{i=1}^{l}(|\Gamma'|<\eta_{|U},f>)=\\
&=<\eta_{|U},f>
\end{align}
\end{remark}

Let $\phi:\Tilde{M}_1/\Gamma_1\to\Tilde{M}_2/\Gamma_2$ be a diffeomorphism 
between two 
orbifolds. We will define the action of $\phi$ on Dirac-type densities. 
Because $\phi$ is a diffeomorphism, there exist open covers 
$\{U_{i,\alpha}\}_{\alpha}$ of the two orbifolds, finite groups 
$\Gamma_{\alpha}$ acting on open subsets $\Tilde{U}_{i,\alpha}$ of the 
Euclidean space $\RR^m$ such that $\Tilde{U}_{i,\alpha}/\Gamma_{\alpha}
\overset{\sim}=U_{i,\alpha}$. The diffeomorphism $\phi$ between 
$U_{1,\alpha}$ and $U_{2,\alpha}$ is induced by a diffeomorphism 
$\Tilde{\phi}_{\alpha}:\Tilde{U}_{1,\alpha}\to\Tilde{U}_{2,\alpha}$
which is $\Gamma_{\alpha}$ equivariant.
If $\eta$ is a Dirac-type density on $M_1=\Tilde{M}_1/\Gamma_1$
then the Dirac-type density $\phi_*(\eta)$ on $M_2=\tilde{M}_2/\Gamma_2$
is given locally on each $U_{2,\alpha}$ by the Dirac-type density 
$\phi_{\alpha}{}_*(\eta_{|U_{1,\alpha}})$. Because of the two remarks we made
above, the Dirac-type distributions associated with $\eta$ and $\phi_*(\eta)$
are conjugated by the isomorphism between the spaces of functions 
induced by $\phi$.

\begin{defn}
Let $M$ a smooth orbifold. A Dirac-type density on $M$ is given by a
collection of Dirac-type densities $\eta_{\alpha}$ on the orbifolds 
$U_\alpha=\Tilde{U}_\alpha/\Gamma_{\alpha}$ such  that
$\phi_{\alpha\beta}{}_*(\eta_{\alpha}{}_{|U_{\alpha}\cap U_{\beta}})=
\eta_{\beta}{}_{|U_{\alpha}\cap U_{\beta}}$, where 
$\phi_{\alpha\beta}$ are the transition maps associated with 
an orbifold atlas $\mathcal{A}=(\Tilde{U}_{\alpha},\Gamma_{\alpha},
U_{\alpha}, \pi_{\alpha})_{\alpha}$.
\end{defn}

The distribution associated with the Dirac-type density $\eta$ on 
$M$ is defined using a partition of unity associated with the orbifold atlas;
if $f\in C^{\infty}_0(M)$ then 
\begin{equation}
<\eta, f>=\sum_{\alpha}<\eta_{\alpha},
\psi_{\alpha}f>
\end{equation}
where $\{\psi_{\alpha}\}_{\alpha}$ is a partition of 
unity associated with the above atlas. We will denote $<\eta, f>$ by 
$\int_M f\eta$.

\begin{defn}
The integral of the Dirac-type density $\eta$ on $M$ is equal to \linebreak
$<\eta, 1>$.
\end{defn}

\subsection{The Canonical Stratification of an Orbifold}

The goal of this section to describe the canonical stratification $\mathcal{S}(M)$ 
of an orbifold $M$.
Let $\Gamma$ be a finite group. 
Let $\mathcal{O}(\Gamma)$ be the set of equivalence classes of subgroups 
of $\Gamma$ with respect to the conjugation relation: $\Gamma_1\sim\Gamma_2$ if
there exists $\gamma\in\Gamma$ such that $\gamma\Gamma_1\gamma^{-1}=\Gamma_2$.
Denote by $[\Gamma']$ the class of $\Gamma'\subseteq \Gamma$.
On $\mathcal{O}(\Gamma)$ we consider the order relation $[\Gamma_1]\preceq[\Gamma_2]$
if there exists $\gamma\in\Gamma$ such that $\gamma\Gamma_1\gamma^{-1}\subseteq
\Gamma_2$. Then $(\mathcal{O}(\Gamma), \preceq)$ is a partially order set.

Let $\Gamma\times \Tilde{M}\to \Tilde{M}$ be a faithful action on a manifold 
$\Tilde{M}$. Let $M=\Tilde{M}/\Gamma$ be the associated orbifold of $\Gamma$ orbits. 
Let $\mathcal{O}(M)\subseteq \mathcal{O}(\Gamma)$ given by 
\begin{equation}
\mathcal{O}(M)=\{[\Gamma']\,|\,\text{ there exists }\Tilde{x}\in\Tilde{M}
\text{ such that }\Gamma_{\Tilde{x}}\sim \Gamma'\}
\end{equation}
\begin{defn}
$\mathcal{O}(M)$ is called the set of orbit types of $M$. 
\end{defn}
We have the map $\mathcal{O}:M\to\mathcal{O}(M)$ given by $\mathcal{O}(x)=
[\Gamma_{\Tilde{x}}]$ 
where $\Tilde{x}\in\Tilde{M}$ is a representative for the class $x\in\Tilde{M}/\Gamma$.
The image of $x$ in $\mathcal{O}(M)$ is called the orbit type of $x$. 

We will define the set of orbit types of a general orbifold $M$. Let $\mathcal{A}
=(\mathcal{R}_{\alpha})$ an atlas with $\mathcal{R}_{\alpha}=(\Tilde{U}_{\alpha},
\Gamma_{\alpha},U_{\alpha}, \pi_{\alpha})$. On the disjoint reunion $\coprod_{\alpha}
\mathcal{O}({U}_{\alpha})$ we consider the equivalence relation generated by 
$[\Gamma'_{\alpha}]\sim[\Gamma'_{\beta}]$ if $\Gamma'_{\alpha}=(\Gamma_{\alpha})_
{\Tilde{x}_{\alpha}}\subseteq\Gamma_{\alpha}$ and $\Gamma'_{\beta}=(\Gamma_{\beta})
_{\Tilde{x}_{\beta}}\subseteq\Gamma_{\beta}$ where $\Tilde{x}_{\alpha}\in
\Tilde{U}_{\alpha}$ and $\Tilde{x}_{\beta}\in\Tilde{U}_{\beta}$ are such that
$\pi_{\alpha}(\Tilde{x}_{\alpha})=\pi_{\beta}(\Tilde{x}_{\beta})=x\in U_{\alpha}\cap
U_{\beta}$. Let $\mathcal{O}(M)=\coprod_{\alpha}\mathcal{O}({U}_{\alpha})/\!{\sim}$
be the set of orbit types of $M$. This definition does not depend on the choice of 
the  atlas because any addition of an extra chart to $\mathcal{A}$ will not
change $\mathcal{O}(M)$. There exists a map $\mathcal{O}:M\to\mathcal{O}(M)$ 
which associates to each $x\in M$ its orbit type, given by $\mathcal{O}(x)= 
[(\Gamma_{\alpha})_{\Tilde{x}}]$
where $\Tilde{x}\in U_{\alpha}$ with $\pi_{\alpha}(\Tilde{x})=x$. We have a partial 
order $\preceq$ on $\mathcal{O}(M)$ induced by the partial order relations on 
$\mathcal{O}(U_{\alpha})\subseteq\mathcal{O}(\Gamma_{\alpha})$. 

For each $\upsilon\in\mathcal{O}(M)$ let 
\begin{equation}
M_{\upsilon}=\{x\in M\,|\, \mathcal{O}(x)=\upsilon\}
\end{equation}
Observe that $M_{[(e)]}$ is equal to the set of regular points $M_{reg}$.

\begin{lemma}
For any $\upsilon\in \mathcal{O}(M)$ the set $M_{\upsilon}$ is a 
smooth submanifold of $M$.
\end{lemma}

\begin{proof}
It is sufficient to prove that any point $x\in M_{\upsilon}$ has a neighborhood
$U$ such that $M_{\upsilon}\cap U$ is a submanifold in $U$. Let $(\Tilde{U},\Gamma,
U,\pi)$ be an linear orbifold chart at $x$, $\pi(0)=x$. Observe that $[\Gamma]=\upsilon$.
Then for any $\Tilde{y}\in\Tilde{U}$ we have $\Gamma_{\Tilde{y}}
\subset \Gamma=\Gamma_0$. So $M_{[\Gamma]}\cap U=\pi(\Tilde{U}^{\Gamma})$.
The fixed point set $\Tilde{U}^{\Gamma}\subset\Tilde{U}$ is a submanifold given by 
linear equations on which $\Gamma$ acts trivially. Then $M_{[\Gamma]}\cap 
U\overset{\sim}=\Tilde{U}^{\Gamma}$ is a smooth submanifold of $U$.
\end{proof}

\begin{defn}
The decomposition 
\begin{equation}
M=\bigcup_{\upsilon\in\mathcal{O}(M)} M_{\upsilon}
\end{equation}
is called the canonical stratification on $M$ and will be denoted by $\mathcal{S}(M)$.  
\end{defn}

A strata $M_{\upsilon}$ is not necessarily connected. Let $M_{\upsilon}=\cup_{\alpha}
M_{\upsilon}^{\alpha}$ be the decomposition in connected components of $M_{\upsilon}$. 

\begin{prop}
\label{eta:strata}
Let $\eta$ be a Dirac-type density on $M$. There exist densities 
$\eta_{\upsilon}$ on the strata $M_{\upsilon}$ such that for any $f\in C^{\infty}_0(M)$ 
we have:
\begin{equation}
<\eta, f>=\sum_{\upsilon\in\mathcal{O}(M)}\int_{M_{\upsilon}} f_{|M_{\upsilon}} 
\eta_{\upsilon}
\end{equation}
\end{prop}

\begin{proof}
We will give the construction of the density $\eta_{[\Gamma']}$ on a linear orbifold 
chart $\mathcal{R}=(\Tilde{U},\Gamma, U, \pi)$. Let $\{\eta^{\gamma}\}_{\gamma
\in\Gamma}$ be the $\Gamma$ Dirac-type density on $\Tilde{U}$ that represents
$\eta_{|U}$.

For $\Gamma'\subseteq\Gamma$ let $\Tilde{U}_{\Gamma'}=\{\Tilde{x}\,|\, \Tilde{x}
\in\Tilde{U}\text{ and }\Gamma_x=\Gamma'\}$. Then $\Tilde{U}_{\Gamma'}$ is a 
submanifold of $\Tilde{U}$. We have $\pi(\Tilde{U}_{\Gamma'})=U_{[\Gamma']}$
and $\pi^{-1}(U_{[\Gamma']})=\Gamma\Tilde{U}_{\Gamma'}=\cup_{\gamma\in\Gamma}
\Tilde{U}_{\gamma\Gamma'\gamma^{-1}}$.  

We will need the following lemma:

\begin{lemma}
The restriction of the projection map
\begin{equation}
\pi_{|\Gamma\Tilde{U}_{\Gamma'}}:\Gamma\Tilde{U}_{\Gamma'}\to U_{[\Gamma']}
\end{equation}
is a covering map with $[\Gamma:\Gamma']$ sheets. 
\end{lemma}
\begin{proof}

$\Gamma_{\Tilde{x}}$ is locally constant on $\Gamma\Tilde{U}_{\Gamma'}$, 
being either $\Gamma'$ or a conjugate of it by an element in $\Gamma$. 
For each $\Tilde{x}\in \Tilde{U}_{\Gamma'}$ one can
choose a small neighborhood $\Tilde{V}_{\Tilde{x}}$ of $\Tilde{x}$ in
$\Tilde{U}_{\Gamma'}$ such that $\Gamma'$ acts trivially on $\Tilde{V}_
{\Tilde{x}}$ and $\gamma\Tilde{V}_{\Tilde{x}}\cap \Tilde{V}_{\Tilde{x}}=
\emptyset$ for $\gamma\in\Gamma\setminus\Gamma'$. Let $x=\pi(\Tilde{x})$ 
and $V_x=\pi(\Tilde{V}_{\Tilde{x}})$.
The map $\pi_{|\Tilde{V}_{\Tilde{x}}}:\Tilde{V}_{\Tilde{x}}\to V_x$ is
a homeomorphism. For any other $\Tilde{y}=\gamma\Tilde{x}\in
\pi^{-1}(x)$ the map $\pi_{|\gamma\Tilde{V}_{\Tilde{x}}}:\gamma
\Tilde{V}_{\Tilde{x}}\to V_x$ is a homeomorphism between a neighborhood of 
$\Tilde{y}$ and $V_x$. Observe also that $\pi^{-1}=\Gamma\Tilde{x}$ has exactly
$[\Gamma:\Gamma']$ points.
\end{proof}

Let us fix $\gamma\in\Gamma$. Then $\Tilde{U}^{\gamma}=\{\Tilde{x}\,|\,
\gamma\Tilde{x}=\Tilde{x}\}$ is a smooth submanifold of $\Tilde{U}$ and $
\Tilde{U}^{\gamma}=\bigcup\limits_{\gamma\in\Gamma'}\Tilde{U}_{\Gamma'}$
is a stratification of $\Tilde{U}^{\gamma}$. Then 
\begin{equation}
\bigcup\limits_{\substack{\gamma\in\Gamma' \\ dim(\Tilde{U}^{\gamma})=
dim(\Tilde{U}_{\Gamma'}^{\alpha})}}\Tilde{U}_{\Gamma'}^{\alpha} \subset 
\Tilde{U}^{\gamma}
\end{equation}
is an open and dense set in $\Tilde{U}^{\gamma}$ ($\{\Tilde{U}_{\Gamma'}
^{\alpha}\}_{\alpha}$ is the decomposition in connected components 
of $\Tilde{U}_{\Gamma'}$).
Let $\eta^{\gamma}_{\Gamma', \alpha}$ be the restriction of the density $\eta^
{\gamma}$ to $\Tilde{U}_{\Gamma'}^{\alpha}$ and 
\begin{equation}\label{eta_Gamma'}
\eta_{\Gamma',\alpha}=\frac{1}{|\Gamma'|}\sum_{\substack{\gamma\in\Gamma' \\
dim(\Tilde{U}^{\gamma})=dim(\Tilde{U}_{\Gamma'}^{\alpha})}}\eta^{\gamma}
_{\Gamma',\alpha}
\end{equation}
Let $U_{[\Gamma']}^{\alpha}=\pi(\Tilde{U}_{\Gamma'}^{\alpha})$. We showed in the 
previous lemma that $\pi$ realizes a local diffeomorphism between 
$\Tilde{U}_{\Gamma'}^{\alpha}$ and $U_{[\Gamma']}^{\alpha}$. 
Let $\eta_{[\Gamma']}^{\alpha}(x)=\pi_*(\eta_{\Gamma',\alpha})(\Tilde{x})$
for some $\Tilde{x}\in \Tilde{U}_{\Gamma'}^{\alpha}$ such that
$\pi(\Tilde{x})=x$. 
This definition is independent of the choice of $\Tilde{x}$ and
of $\Gamma'\subset\Gamma$. Indeed, let $\gamma\in\Gamma$ be such that 
$\gamma\Gamma'\gamma^{-1}$ is another 
representative for $[\Gamma']$ and $\pi(\gamma\Tilde{x})=x$.
After a reordering of the indices $\alpha$
we can suppose that $\Tilde{U}_{\gamma\Gamma'\gamma^{-1}}^{\alpha}=
\gamma\Tilde{U}_{\Gamma'}^{\alpha}$. Because $\{\eta^{\gamma'}\}_
{\gamma'}$ is a $\Gamma$ Dirac-type density on $\Tilde{U}$ we have
$\eta_{\gamma\Gamma\gamma^{-1},\alpha}(\gamma\Tilde{x})=\gamma_*(\eta_
{\Gamma',\alpha} (\Tilde{x}))$ and because $\pi$ commutes with 
the action of $\Gamma$ on $\Tilde{U}$ we get the independence 
of $\eta_{[\Gamma']}(x)$ of the choices we made. 

If $\eta_{[\Gamma']}^{\alpha}$ was not defined on $U_{[\Gamma']}^{\alpha}$ 
because there were no $\gamma\in\Gamma$ such that $\Tilde{U}_{\Gamma'}^{\alpha}$ 
is an open submanifold in $\Tilde{U}^{\gamma}$, we take 
$\eta_{[\Gamma']}^{\alpha}=0$.

We define $\eta_{[\Gamma']}$ on $M_{[\Gamma']}$ to be the density whose restriction 
to $U_{[\Gamma']}^{\alpha}$ is equal to $\eta_{[\Gamma']}^{\alpha}$.
We need to show that 
\begin{equation}
\label{eta:strata:eq}
<\eta, f>=\sum_{[\Gamma']\in\mathcal{O}(M)}
\int_{M_{[\Gamma']}} f_{|M_{[\Gamma']}} \eta_{[\Gamma']}\quad
\text{ for any } f\in  C^{\infty}_0(M).
\end{equation}
 
It is sufficient to prove this for $f\in C^{\infty}_0(U)$.
If $\Tilde{f}\in C^{\infty}_0(\Tilde{U})^{\Gamma}$ is such that
$\Tilde{f}=\pi f$ then
\begin{align}
<\eta, f>&=\frac{1}{|\Gamma|}\sum_{\gamma\in\Gamma}\int_{\Tilde{U}^
{\gamma}} \Tilde{f} \eta^{\gamma}=\\
&=\frac{1}{|\Gamma|}\sum_{\gamma\in\Gamma}(\sum_{\substack{\Gamma'
\ni \gamma \\ dim(\Tilde{U}^{\gamma})= dim(\Tilde{U}^{\alpha}_{\Gamma'})}}
\int_{\Tilde{U}^{\alpha}_{\Gamma'}}\Tilde{f}\eta^{\gamma}_
{\Gamma',\alpha})=\\
&=\frac{1}{|\Gamma|}\sum_{\alpha,\Gamma'<\Gamma}(\sum_{\substack{
\gamma\in\Gamma' \\ dim(\Tilde{U}^{\gamma})= dim(\Tilde{U}^{\alpha}_{\Gamma'})}}
\int_{\Tilde{U}^{\alpha}_{\Gamma'}}\Tilde{f}\eta^{\gamma}_{\Gamma',\alpha})=\\
&=\frac{1}{|\Gamma|}\sum_{\upsilon\in\mathcal{O}(\Gamma)}
(\sum_{\alpha, \Gamma'\in\upsilon} \int_{\Tilde{U}^{\alpha}_{\Gamma'}}
\Tilde{f}(|\Gamma'| \eta_{\Gamma',\alpha}))=
\qquad\text{(by \eqref{eta_Gamma'})}\\
&=\frac{1}{[\Gamma:\Gamma']}\sum_{[\Gamma']\in\mathcal{O}(U)}
\int_{\pi^{-1}(U^{\alpha}_{[\Gamma']})} \Tilde{f} \eta_{\Gamma',\alpha}=
\sum_{[\Gamma']\in\mathcal{O}(U)}\int_{U_{[\Gamma']}} f \eta_{[\Gamma']}
\end{align}
which is equal to the right side of \eqref{eta:strata:eq}. In the last step
we used the fact that $\pi:\pi^{-1}(U^{\alpha}_{[\Gamma']})\to 
U_{[\Gamma']}^{\alpha}$ is a covering map with $[\Gamma:\Gamma']$ leaves.
\end{proof}

\vfill\eject


\section{Pseudodifferential Operators in Orbibundles}
\subsection{Pseudodifferential Operators -- The Local Model}

We remind the reader that locally any vector orbibundle over an orbifold is 
obtained as $\Gamma$ orbits of a $\Gamma$ vector bundle, with
$\Gamma$ finite. We begin by analyzing the case where $\Gamma$ acts freely on
a vector bundle $E\overset{p}\to M$, in which case the associated vector
orbibundle of $\Gamma$ orbits $E/\Gamma \overset{\overline{p}}\to B/\Gamma$
is a genuine vector bundle over a manifold.
Let $R:C^{\infty}(M;E)\to C^{\infty}(M/\Gamma;E/\Gamma)$ be defined as
\begin{equation}
R(f)(\overline{x})=\frac1{|\Gamma|}\sum_{x\in \overline{x}}\overline{f(x)}
\qquad\text{ for } f\in C^{\infty}(M;E)\text{ and }\overline{x}\in M/\Gamma 
\end{equation}
and $I:C^{\infty}(M/\Gamma;E/\Gamma)\to C^{\infty}(M;E)$ by
\begin{equation}
I(f)(x)=v \in E_x \quad \text{ for }f\in C^{\infty}(M/\Gamma;E/\Gamma)
\text{ and }x\in M
\end{equation}
where $v\in E_x$ is such that $f(\Gamma x)=\Gamma v$. 
$v$ is unique in $p^{-1}(x)=E_x$ because the action of $\Gamma$ is free.
$I$ identifies $C^{\infty}(M/\Gamma;E/\Gamma)$ with the space of 
$\Gamma$ invariant sections $C^{\infty}(M;E)^{\Gamma}$ and then $R$ is 
the averaging operator over the group $\Gamma$. Obviously, $R\cdot I=Id$ and
$I\cdot R$ is the averaging operator over the group $\Gamma$. If we endow 
the quotient vector bundle with the trivial $\Gamma$ action then 
both $R$ and $I$ are $\Gamma$ equivariant.

Let $\Psi(M;E)$ and $\Psi(M/\Gamma;E/\Gamma)$ be the spaces of 
pseudodifferential operators on $C^{\infty}(M;E)$ respectively 
$C^{\infty}(M/\Gamma;E/\Gamma)$. Let
$\mathbi{R}:\Psi(M;E)\to \Psi(M/\Gamma;E/\Gamma)$ given by 
\begin{equation}
\mathbi{R}(A)=R A I\quad \text{ for }A\in \Psi(M;E)
\end{equation}
$\mathbi{R}(A)$ is a pseudodifferential operator acting on $C^{\infty}
(M/\Gamma;E/\Gamma)$. To see this, observe that if the distributional kernel
of the operator $A$ is equal to $K_{A}(x,y)$ then the distributional kernel of
$\mathbi{R}(A)$ is given by the formula 
\begin{equation}
K_{\mathbi{R}(A)}(\overline{x},\overline{y})=\frac1{|\Gamma|}
\sum_{\substack{x\in 
\overline{x}\\ y\in \overline{y}}}\overline{K_{A}(x,y)}
\end{equation}
In a coordinate chart $(U,\phi)$ the total symbol of the operator 
$\mathbi{R}(A)$ is given by 
\begin{equation}
\overline{a}(\overline{x},\overline{\xi})=\dfrac{1}{|\Gamma|}\sum_
{\substack{x\in\overline{x} \\ \xi\in\overline{\xi}}}a(x,\xi)
\end{equation}
where $a(x,\xi)$ is the total symbol of $A$ in the chart $(p^{-1}U, \phi p)$.
If $B$ is the average 
of $A$ with respect to the action of $\Gamma$ by conjugation on $\Psi(M;E)$, 
$B=\dfrac{1}{|\Gamma|}\sum_{\gamma\in\Gamma}\gamma A \gamma^{-1}$, then
$\mathbi{R}(A)=\mathbi{R}(B)$. Moreover, $\mathbi{R}(B)$ can be seen as the 
restriction of the $\Gamma$ invariant operator $B$ to $C^{\infty}(M;E)^{\Gamma}
\overset{\sim}=C^{\infty}(M/\Gamma;E/\Gamma)$.

\begin{lemma}\label{bold:I}
There exists a $\Gamma$ equivariant operator $\mathbi{I}:
\Psi(M/\Gamma;E/\Gamma)\to 
\Psi(M;E)$ so that $\mathbi{R}\mathbi{I}=Id$ on $\Psi(M/\Gamma;E/\Gamma)$.  
\end{lemma}

\begin{proof}
If we define $\mathbi{I}$ by $\mathbi{I}(A)=IAR$ then
$\mathbi{R}\mathbi{I}=Id$, but $IAR$ will not be, in general, 
a pseudodifferential operator, unless $A$ is smoothing and then 
the smooth kernel $K_{\mathbi{I}(A)}$ of $\mathbi{I}(A)$ is 
defined uniquely by $\overline{K_{\mathbi{I}(A)}(x,y)}=
K_A(\overline{x},\overline{y})$.

We will define $\mathbi{I}$ using a partition
of unity of $M/\Gamma$ with smooth functions $\{\phi_{\lambda}\}_
{\lambda\in \Lambda}$ subordinated to a cover of $M/\Gamma$ with 
coordinate chart neighborhoods $(U_i)_{i\in I}$ which
are slices with respect to the $\Gamma$ action on $M$. We will 
choose $\phi_{\lambda}$  so that for $\lambda, \lambda' \in \Lambda$ 
at least one of the 
following conditions holds:
\begin{enumerate}
\item[(i)]$supp(\phi_{\lambda})\cap supp(\phi_{\lambda'})=\emptyset$
\item[(ii)]$supp(\phi_{\lambda})\cup supp(\phi_{\lambda'})$ is included in
the same coordinate neighborhood.
\end{enumerate}
If $A\in \Phi(M/\Gamma;E/\Gamma)$
then $A=\sum_{\lambda,\lambda'}\phi_{\lambda}A\phi_{\lambda'}$. Denote
$A_{\lambda\lambda'}=\phi_{\lambda}A\phi_{\lambda'}$. We will define 
$\mathbi{I}_{\lambda\lambda'}(A)=\mathbi{I}(A_{\lambda\lambda'})$ and 
then $\mathbi{I}=\sum_{\lambda,\lambda'}\mathbi{I}_{\lambda\lambda'}$.

If $supp(\phi_{\lambda})\cap supp(\phi_{\lambda'})=\emptyset$, then 
$A_{\lambda\lambda'}$ is smoothing and we define
$\mathbi{I}_{\lambda\lambda'}(A)=IA_{\lambda\lambda'}R$.

Otherwise, choose $i\in I$ such that $supp(\phi_{\lambda})\cup 
supp(\phi_{\lambda'}) \subset U_i$. Because $U_i$ is a slice, 
its preimage $\Tilde{U}_i\subset M$ can be 
identified with $\Gamma\times U_i$ endowed with the $\Gamma$ action by 
left translations. We have the induced $\Gamma$ equivariant isomorphism
\begin{equation}\label{block}
C^{\infty}(\Tilde{U}_i ;E_{|\Tilde{U}_i})\overset{\sim}=\Gamma\times 
C^{\infty}(U_i;E/\Gamma_{|U_i})
\end{equation}
The operator $A_{\lambda\lambda'}$ is localized above $U_i$ and we define 
$\mathbi{I}(A_{\lambda\lambda'})\in \Psi(\Tilde{U}_i ;E_{|\Tilde{U}_i})$ 
to be the block-diagonal operator with diagonal cells equal to
$A_{\lambda\lambda'}$, where the block representation is with 
respect to the decomposition given by the isomorphism \eqref{block}.
Then $\mathbi{I}(A_{\lambda\lambda'})$ is a $\Gamma$ invariant 
pseudodifferential operator on $C^{\infty}(M;E)$ localized above 
$\Tilde{U}_i$ and 
$\mathbi{R}\mathbi{I}(A_{\lambda\lambda'})=A_{\lambda\lambda'}$.

Define $\mathbi{I}$ to be the sum $\sum_{\lambda,\lambda'}\mathbi{I}_
{\lambda\lambda'}$. Then $\mathbi{R}\mathbi{I}=Id$  and $Im(\mathbi{I})$ is a 
subset of the $\Gamma$ equivariant pseudodifferential operators in 
$C^{\infty}(M;E)$.
So $\mathbi{I}$ is $\Gamma$ equivariant.
\end{proof}

Thus every pseudodifferential operator $A$ in $C^{\infty}(M/\Gamma;E/\Gamma)$ 
can be represented by the restriction to the $\Gamma$ invariant sections 
of a $\Gamma$ equivariant pseudodifferential operator in $C^{\infty}(M;E)$ 
given by $\mathbi{I}(A)$.

\begin{remark}
\label{gamma:prime}
In the presence of the action of another finite group $\Gamma'$ on the vector
bundle $E\overset{p}\to M$ which is free and commutes with the action of
$\Gamma$, we have an induced natural action of $\Gamma'$ on the 
quotient vector bundle
$E/\Gamma \overset{p}\to M/\Gamma$; $\gamma'\in\Gamma'$ acts on $\Gamma x$ 
by $\gamma'\Gamma x=\Gamma \gamma' x$. This induced action might not be free.
The maps $R$ and $I$ between the spaces of smooth
sections in the two vector bundles are $\Gamma'$ equivariant.
Indeed, the $\Gamma'$ action on the $\Gamma$ invariant sections is
the restriction of the $\Gamma'$ action on all sections, and because
$R$ and $I$ are essentially the projection onto the $\Gamma$ 
invariant sections 
and the inclusion of the invariant sections into the space of all sections, 
the two maps are $\Gamma'$ equivariant.
As a consequence, $\mathbi{R}$ is $\Gamma'$ invariant. The components
$\mathbi{I}_{\lambda,\lambda'}$ of $\mathbi{I}$ for $supp(\phi_{\lambda})\cap 
supp(\phi_{\lambda'})=\emptyset$ are $\Gamma'$ equivariant as well.
To show that $\mathbi{I}_{\lambda,\lambda'}$ are $\Gamma'$ equivariant
for $supp(\phi_{\lambda})\cup supp(\phi_{\lambda'})
\subset U_i$ for a slice $U_i\subset M/\Gamma$, observe that the 
identification of the preimage $\Tilde{U}_i\subset M$ of $U_i$  
with $\Gamma\times U_i$ is
$\Gamma'$ equivariant as well. The map $\mathbi{I}_{\lambda,\lambda'}$
from $\Psi(U_i,E_{|U_i})$ to $\Psi(\Tilde{U}_i ;E_{|\Tilde{U}_i})$ given by 
the block-diagonal construction as above is $\Gamma'$ equivariant. 
So $\mathbi{I}$ is $\Gamma'$ equivariant. 
\end{remark}

We will use the above facts to define the pseudodifferential operators 
acting on the space of smooth sections of a vector orbibundle. We will start 
with the local construction.

Let $\Tilde{U}\times V\overset{pr_1}\to \Tilde{U}$ be a $\Gamma$ vector 
bundle, with $\Gamma$ a finite group. Suppose that $\Tilde{U}\subset 
\RR^m$ is an open subset with compact closure. The group $\Gamma$ 
acts on the space of sections $C^{\infty} (\Tilde{U}, \Tilde{U}\times V)$. 
Let $\Gamma$ act on the space of operators on $C^{\infty}(\Tilde{U}, 
\Tilde{U}\times V)$ by conjugation: if $A$ is an operator then 
$(\gamma A)(f)=\gamma A(\gamma^{-1} f)$ for any $f\in C^{\infty}(\Tilde{U},
\Tilde{U}\times V)$ and $\gamma \in \Gamma$.

Let $(\Tilde{U}\times V)/\Gamma\overset{\overline{pr}_1}\to \Tilde{U}/\Gamma$ 
be the associated vector orbibundle.

\begin{defn} An operator $A$ acting on $C^{\infty}_0(\Tilde{U}/\Gamma; 
(\Tilde{U}\times V)/\Gamma)\overset{\sim}=C^{\infty}_0(\Tilde{U}; 
\Tilde{U}\times V)^{\Gamma}$ is called a pseudodifferential operator 
if $A$ is the restriction to the $\Gamma$ invariant sections 
of a $\Gamma$ invariant pseudodifferential operator $\Tilde{A}$ acting 
on $C^{\infty}_0
(\Tilde{U}; \Tilde{U}\times V)$.
\end{defn}

We will show later in Proposition \ref{unique}, in a greater generality, that 
the operator $\Tilde{A}$ whose restriction to the invariant
sections is equal to $A$ is unique up to smoothing operators. 

We define the order of $A$ to be equal to the order of $\Tilde{A}$.
We denote the space of pseudodifferential operators of order $d$ by 
$\Psi^d(\Tilde{U}/\Gamma; (\Tilde{U}\times V)/\Gamma)$ or simply by $\Psi^d$
and the space of all pseudodifferential operators by $\Psi$. The space $\Psi$
is a filtered algebra with respect to the composition of operators, where the 
filtration is given by the degree of the operators. 
An operator $A$ is smoothing if the operator $\Tilde{A}$ is smoothing. 
The space of smoothing operators will be denoted by $\Psi^{-\infty}$.

\begin{defn}
$A$ is called a classical pseudodifferential operator if $\Tilde{A}$
is classical. (as described in \cite{shubin}, Section 3.7).
\end{defn}

Let $\Tilde{a}(x,\xi)$ be the total symbol of $\Tilde{A}$. $\Tilde{a}$ 
is a section of the 
pull-back of the endomorphism vector bundle to the tangent space, $\Tilde{a}\in
C^{\infty}(T^*(\Tilde{U}), T^*(\Tilde{U})\times End(V))$. The principal symbol 
$\Tilde{a}_{pr}$ is a section of the same vector bundle. 
The group $\Gamma$ acts naturally on the vector bundle $T^*(\Tilde{U})
\times End(V)\overset{pr_1}\to T^*(\Tilde{U})$. Then $\Gamma$ acts 
on the sections of this vector bundle and $\gamma \Tilde{a}_{pr}$ 
is the principal symbol of the operator $\gamma \Tilde{A}$ (cf. Lemma 1.1.3 
\cite{gilkey}). This is the direct consequence of the invariance of 
the principal symbol with respect to changes of coordinates. Because
$\Tilde{A}$ is $\Gamma$ invariant we conclude that $\Tilde{a}_{pr}$ is 
a $\Gamma$ invariant section so it defines a smooth section in the 
vector orbibundle
$(\Tilde{U}\times End{V})/\Gamma\overset{\overline{pr}_1}\to\Tilde{U}/\Gamma$. 
Because $\Tilde{A}$ is unique up to smoothing operators,
$\Tilde{a}_{pr}$ depends only on $A$ and will be called the principal symbol
of the operator $A$. 

In the case when the $\Gamma$ action on the vector bundle $\Tilde{U}\times V
\overset{pr_1}\to \Tilde{U}$ is given by the restriction to $\Tilde{U}\subset 
\RR^m$ of a linear action of $\Gamma$ on $\RR^m$ and 
by a $\Gamma$ representation $\rho:\Gamma\times V\to V$,  we will be 
able to define the total symbol of the pseudodifferential operator $A$. 
This is the consequence of the invariance of the total symbol of 
$\Tilde{A}$ with respect to linear changes of coordinates. The total 
symbol of $A$ is the given by $\Tilde{a}\in C^{\infty}_0((T^*(\Tilde{U})
\times V)/\Gamma; T^*(\Tilde{U}))$ and is unique up to a smoothing symbol. 

If $\Tilde{A}$ is classical and
\begin{equation}
\label{asymptotic:exp}
\Tilde{a}(x,\xi)\sim\sum_{i\ge 0} \Tilde{a}_{d-i}(x,\xi)
\end{equation}
is an asymptotic expansion of the total symbol with $\Tilde{a}(x, \xi)_{d-i}$
being homogeneous symbols of degree of homogeneity $d-i$, then each 
homogeneous component is $\Gamma$ invariant and will define the 
homogeneous component $a_{d-i}$ of the
asymptotic expansion of the total symbol of the operator $A$. The formal 
sum in the right-hand side of the equality \eqref{asymptotic:exp} is 
called the asymptotic symbol of the operator $A$.

If a pseudodifferential operator $A$ of order $d$ acts on the sections of the 
vector orbibundle $(\Tilde{U}\times V)/\Gamma \overset{\overline{pr}_1}
\to \Tilde{U}/\Gamma$ then, as in the case of genuine vector bundles, 
it extends to a continuous operator between Sobolev spaces 
$A:H_s(\Tilde{U}/\Gamma; (\Tilde{U}\times V)/\Gamma)
\to H_{s-d}(\Tilde{U}/\Gamma; (\Tilde{U}\times V)/\Gamma)$ (cf. Lemma 1.2.1 
\cite{gilkey})

\begin{remark} 
\label{smoothing}
If $A\in \Psi^{-\infty}$ then there exists a $\Gamma$ invariant
smoothing operator $\Tilde{A}$ whose restriction to the $\Gamma$ invariant 
sections is 
equal to $A$. Because $\gamma \Tilde{A}=\Tilde{A} \gamma$ for all 
$\gamma\in \Gamma$,
the kernel $\Tilde{K}(\Tilde{x}, \Tilde{y})$ of $\Tilde{A}$ satisfies 
the equality
$\gamma\Tilde{K}(\gamma^{-1}\Tilde{x}, \gamma^{-1}\Tilde{y})=\Tilde{K}
(\Tilde{x},\Tilde{y})$ so
it defines a $\Gamma$ invariant section in the vector bundle $(\Tilde{U}
\times\Tilde{U})\times\End(V)\overset{pr_1}\to (\Tilde{U}\times\Tilde{U})$. 
Here 
$\Gamma$ acts by the diagonal action on $\Tilde{U}\times\Tilde{U}$. 
If we define
\begin{equation}
\Tilde{K}'(\Tilde{x}, \Tilde{y})=\frac {1}{|\Gamma|}\sum_{\gamma\in \Gamma}
\gamma \Tilde{K}(\Tilde{x}, \gamma^{-1}\Tilde{y})
\end{equation}
then $\Tilde{K}'$ is the kernel of the $\Gamma$ invariant smoothing operator
$\Tilde{A}'=\frac{1}{\Gamma}\sum_{\gamma\in \Gamma}\Tilde{A}\gamma$ whose 
restriction to the $\Gamma$ invariant sections induces $A$. $\Tilde{K}'$ 
defines
a $\Gamma\times\Gamma$ invariant smooth section of the vector bundle 
$(\Tilde{U}
\times\Tilde{U})\times\End(V)\overset{pr_1}\to (\Tilde{U}\times\Tilde{U})$ 
where 
$\Gamma\times\Gamma$ acts by the product action on $\Tilde{U}\times\Tilde{U}$.
Then $\Tilde{U}\times\Tilde{U}/\Gamma\times\Gamma\overset{\sim}=U\times U$ and 
$\Tilde{K}'$ defines a smooth section of the endomorphism vector orbibundle
$(\Tilde{U}\times\Tilde{U})\times \End(V)/\Gamma\times\Gamma\overset{pr_1}\to
U\times U$ which will be the smooth kernel of the operator $A$. 
\end{remark}

\begin{defn}
\label{elliptic:local}
The pseudodifferential operator $A$ of order $d$ is called elliptic above 
an open set
$U_1\subset \Tilde{U}/\Gamma$ if the operator $\Tilde{A}$ is elliptic above
the preimage of $U_1$ in $\Tilde{U}$.
\end{defn}

Let $\Tilde{U}_1\subset \Tilde{U}$ be the preimage of $U_1$.
The above definition implies the existence of a parametrix $\Tilde{B}'$ 
for $\Tilde{A}$
above $\Tilde{U}_1$. Then $\Tilde{B}=\frac{1}{|\Gamma|} \sum_{\gamma\in 
\Gamma} \gamma
\Tilde{B}'\gamma^{-1}$ is a $\Gamma$ invariant parametrix for $\Tilde{A}$ 
and defines 
a parametrix $B\in \Psi^{-d}$ for $A$, $BA-Id\in \Psi^{-\infty}$, $AB-Id\in 
\Psi^{-\infty}$.

\subsection{Pseudodifferential Operators in Orbibundles}

Let $E\overset{p}\to M$ be a vector orbibundle and $A$ be a linear operator
acting on the space of smooth sections $C^{\infty}(M;E)$. 
If $\mathcal{R}=(\Tilde{U},\,V,\,\Gamma,U,\Pi,\pi)$ is 
an orbibundle chart for $E\overset{p}\to M$ above $U\subset M$ then 
$\Tilde{U}\times V\overset{pr_1}\to \Tilde{U}$ is a $\Gamma$ vector bundle and 
passing to the orbit space yields the orbibundle 
$\Tilde{U}\times V/\Gamma \overset{\overline{pr}_1}\to \Tilde{U}/\Gamma$. 
This orbibundle is isomorphic to the
restriction of $E\overset{p}\to M$ to $U$ via the orbibundle map 
$(\overline{\Pi},\overline{\pi})$:
$$
\xymatrix{
\Tilde{U}\times V/\Gamma\ar@{-->}[r]^{\overline{\Pi}} \ar[d]_{\overline{pr}_1}
& E_{|U}\ar[d]^{p}
\\
\Tilde{U}/\Gamma\ar@{-->}[r]_{\overline{\pi}}&U}
$$
The orbibundle $E\overset{p}\to M$ can be covered by such orbibundle maps.

Let $\phi$ and $\psi$ two smooth functions on $M$ such that their support 
is included in $U$. We will denote by the same letter the multiplication 
operator by the respective functions. Then the operator $\phi\circ A\circ\psi$
has support in $U$ and takes sections with support in $U$ to sections
with support in $U$. Using the orbibundle isomorphism $\overline{\pi}$, 
one can define the operator $A_{\mathcal{R}}$ acting on the sections 
of $\Tilde{U}\times V/\Gamma \to \Tilde{U}/\Gamma$ as
\begin{equation}
\label{localize:op}
A_{\mathcal{R}}\,f=\overline{\pi} (\phi\circ A\circ\psi){\overline
{\pi}}^{-1}\,f.
\end{equation}

\begin{defn} $A$ is a pseudodifferential operator if $A_{\mathcal
{R}}$ is a pseudodifferential operator for any orbibundle chart $\mathcal{R}$.
\end{defn}

In order for the previous definition to be consistent, we need to show that 
this property is independent of the choice of a chart. More precisely:

\begin{prop} 
\label{pdounique}
Let $\Tilde{E}_i\overset{p_i}\to\Tilde{M}_i$ be two
$\Gamma_i$ vector bundles, $i=1,2$, such that the corresponding orbibundles
are isomorphic via an isomorphism $( T, t)$. Let $A_1$ be an operator 
acting on sections of $\Tilde{E}_1/\Gamma_1\overset{p_1}\to
\Tilde{M}_1/\Gamma_1$ and $A_2= (T,t)_*\circ A_1\circ (T,t)^*$ the 
corresponding 
operator acting on $\Tilde{E}_2/\Gamma_2\overset{p_2}\to
\Tilde{M}_2/\Gamma_2$. 

Then there exists a $\Gamma_1$ equivariant pseudodifferential operator 
$\Tilde{A}_1$ acting on sections of $\Tilde{E}_1\overset{p_1}\to
\Tilde{M}_1$ whose restriction to $\Gamma_1$ 
invariant sections is equal to $A_1$ if and only if  there exists a 
$\Gamma_2$ equivariant 
pseudodifferential operator $\Tilde{A}_2$ acting on sections of 
$\Tilde{E}_2\overset{p_2}\to\Tilde{M}_2$
whose restriction to $\Gamma_2$ invariant sections is equal to $A_2$.
\end{prop}

\begin{proof} For each $x_i \in \Tilde{M}_i/\Gamma_i$ such that $ t(x_1)=
x_2$ choose $\tilde{x}_i \in \Tilde{M}_i$ with $\pi_i(\tilde{x}_i)=x_i$.
Because the orbibundles are isomorphic, there exists a group isomorphism 
$\alpha$ between $(\Gamma_i)_{\tilde{x}_i}$--the corresponding isotropy 
groups of  $\tilde{x}_i$ in $\Gamma_i$ and one can find neighborhoods 
$\Tilde{W}_i$ of $\tilde{x}_i$ in $\Tilde{M}_i$ which are $(\Gamma_i)_
{\tilde{x}_i}$ invariant and an $\alpha$-equivariant bundle diffeomorphism 
$( L, l)$  that makes the following diagram commutative:
\begin{equation}\label{local:diffs}
\xymatrix{
\Tilde{E}_1 \ar[dd]_{p_1}& & & \Tilde{E}_2\ar[dd]^
{p_2}\\
&\Tilde{W}_1\times V \ar@{_{(}->}[ul] \ar[r]^{ L}_{\sim} \ar[dd]_{p_1} & 
\Tilde{W}_2 \times V \ar@{^{(}->}[ur] \ar[dd]^{p_2}& \\
\Tilde{M}_1 \ar[dd]_{\pi_1}& & & \Tilde{M}_2\ar[dd]^{\pi_2}\\
&\Tilde{W}_1 \ar@{_{(}->}[ul] \ar[r]^{ l}_{\sim} \ar@{-}'[d]_{\pi_1}
\ar[dl]_{\pi_1} &
\Tilde{W}_2 \ar@{^{(}->}[ur] \ar[dr]^{\pi_2} \ar@{-}'[d]^{\pi_2}& \\
\Tilde{M}_1/\Gamma_1 \ar[rrr]^{ t}_{\sim} & \ar[d]&\ar[d] & \Tilde{M}_2/
\Gamma_2\\
&\Tilde{W}_1/(\Gamma_1)_{\tilde{x}_1}\ar@{_{(}->}[ul]
\ar[r]^{\overline{l}}_{\sim} & \Tilde{W}_2/(\Gamma_2)_{\tilde{x}_2}
\ar@{^{(}->}[ur]&
}
\end{equation}

Let $\{\phi_{\lambda}\}_{\lambda\in\Lambda}$ be a partition of unity of 
$\Tilde{M}_i/\Gamma_i$ subordinated 
to the open cover $\Tilde{W}_i/(\Gamma_i)_{\tilde{x}_i}$
(because $\Tilde{M}_i/\Gamma_i$ are diffeomorphic and the diffeomorphism $t$
permutes the two open covers, we will denote the two partition of 
unity with the
same letters, thought we will refer to functions on two different, 
but diffeomorphic 
orbifolds). We can choose this partition so that for any $\lambda,\,
\lambda' \in \Lambda$ at least one of the following conditions holds:
\begin{enumerate}
\item[(i)]$supp(\phi_{\lambda})\subset \Tilde{W}_i/(\Gamma_i)_{\tilde{x}_i}$, 
$\quad supp(\phi_{\lambda'})\subset \Tilde{W'}_i/(\Gamma_i)_{\tilde{x}_i'}$ and
\newline
$\Tilde{W}_i/(\Gamma_i)_{\tilde{x}_i}\cap \Tilde{W'}_i/(\Gamma_i)_
{\tilde{x}_i'}
=\emptyset\qquad$
\item[(ii)]$supp(\phi_{\lambda})\cup supp(\phi_{\lambda'})\subset \Tilde{W}_i/
(\Gamma_i)_{\tilde{x}_i}$.
\end{enumerate}
We will show only one implication of the proposition, the other implication 
can be proved similarly.
Let $\Tilde{A}_1$ be a  $\Gamma_1$ equivariant pseudodifferential 
operator acting on the smooth sections of $\Tilde{E}_1
\overset{p_1}\to \Tilde{M}_1$ whose restriction to the 
invariant sections is equal to $A_1$. We will construct a lift $\Tilde{A}_2$ 
as required in the proposition.

Because 
\begin{equation}\label{lambda:decomp}
A_1=\sum_{\lambda,\lambda'\in\Lambda}\phi_{\lambda}A_1 \phi_{\lambda'}
\end{equation}
and because a lift for the operator $\phi_{\lambda}A_1
 \phi_{\lambda'}$ is provided by $\tilde{\phi_{\lambda}}
\Tilde{A}_1\tilde{\phi_{\lambda'}}$, we reduced our problem to the 
operators of the form $\phi_{\lambda}A_1 \phi_{\lambda'}$ which we
will denote by $(A_1)_{\lambda,\lambda'}$ for simplicity (we denoted by 
$\tilde{\phi}$ 
the lift of the smooth function $\phi$).
Denote the lifts $\tilde{\phi_{\lambda}}\Tilde{A}_1\tilde{\phi_{\lambda'}}$ 
by $(\Tilde{A}_1)_{\lambda,\lambda'}$.  
In the first case  the operator $(\Tilde{A}_1)_{\lambda,\lambda'}$ 
is smoothing.
In the second case the operator $(\Tilde{A}_1)_{\lambda,\lambda'}$ 
is localized to the
open set $\Gamma_1\cdot\Tilde{W}_1\overset{\sim}=\Gamma_1\times_
{(\Gamma_1)_{\tilde{x}_1}}\Tilde{W}_1$. We will treat the two cases separately.

In the first case define: 
\begin{equation}\label{bi:average}
\overline{A}_1=\sum_{\gamma\in\Gamma_1}(\Tilde{A}_1)_{\lambda\lambda '}\,\gamma
\end{equation}
Because $(\Tilde{A}_1)_{\lambda\lambda '}$ was smoothing and $\Gamma_1$ 
equivariant, $\overline{A}_1$ is smoothing and $\Gamma_1$ bi-invariant: 
\begin{equation}\label{bi:invariant}
\mu \overline{A}_1=\overline{A}_1\nu=\overline{A}_1\qquad \text{ for any }
\mu,\nu\in\Gamma_1 
\end{equation}
Moreover $\overline{A}_1$ induces the same operator $(A_1)_{\lambda\lambda'}$ 
at the orbibundle level because its action on the $\Gamma_1$ invariant 
sections is not changed. 
If we chooose a $\Gamma_1$ invariant metric on $M_1$ then the 
kernel $\overline{K}_1(\tilde{x},\tilde{y})$ of the operator $\overline{A}_1$
is a smooth section in $C^{\infty}(\Tilde{M}_1\times \Tilde{M}_1; 
\End(\Tilde{E}_1))$
and because of the equalities \eqref{bi:invariant} it satisfies the condition
$\mu K_1(\mu^{-1}\tilde{x},\tilde{y})=\nu K(\tilde{x},\nu^{-1}\tilde{y})$ 
for all $\mu,\nu\in\Gamma_1$. 
This implies that $K_1$ is a $\Gamma_1\times\Gamma_1$ invariant smooth 
section in the 
$\Gamma_1\times\Gamma_1$ vector bundle $\End(\Tilde{E}_1)\to\Tilde{M}_1
\times\Tilde{M}_1$ and so it defines a smooth section in the vector orbibundle
$\End(\Tilde{E}_1/\Gamma_1)\to \Tilde{M}_1/\Gamma_1\times \Tilde{M}_1/
\Gamma_1$. 
Because the vector orbibundles $\Tilde{E}_i/\Gamma_i\overset{p_i}\to 
\Tilde{M}_i/\Gamma_i$, 
$i=1,2$,  are isomorphic, there exists a unique smooth $\Gamma_2\times
\Gamma_2$ invariant
section $\overline{K}_2$ in $\End(\Tilde{E}_2)\to\Tilde{M}_2\times
\Tilde{M}_2$ which induces 
the same section via the isomorphism $(T,t)$ at the orbibundle level as  
$\overline{K}_1$. 
The smoothing operator $\overline{A}_2$ with kernel $\overline{K}_2$ is 
$\Gamma_2$ bi-invariant, 
and so it is $\Gamma_2$ equivariant. The operator $\overline{A}_2$ induces 
at the orbibundle level 
the operator $(T,t)_* \circ (A_1)_{\lambda,\lambda '} \circ (T,t)^*$. 
The statement of the proposition is proved in the case (i). 

In the case (ii) $(\Tilde{A}_1)_{\lambda \lambda'}$ is a $\Gamma_1$ 
equivariant operator in
$\Psi( \Gamma_1\Tilde{W}_1;\Tilde{E}_1{}_{|\Gamma_1\Tilde{W}_1})$. 
Let us denote 
$\Gamma_1'= (\Gamma_1)_{\tilde{x}_1}$. The $\Gamma_1$ 
vector bundle $\Tilde{E}_1{}_{|\Gamma_1\Tilde{W}_1}\to\Gamma_1\Tilde{W}_1$ 
is isomorphic with
$\Gamma_1\times_{\Gamma_1'}(\Tilde{W}_1\times V)\to\Gamma_1\times_{\Gamma_1'}
\Tilde{W}_1$.
Observe that $\Gamma_1\times_{\Gamma_1'}\Tilde{W}_1=(\Gamma_1\times
\Tilde{W}_1)/\Gamma_1'$.
We will use the observations regarding the operators $\mathbi{R}$ 
and $\mathbi{I}$ made at the beginning of 
section 3.1.
Cf. Lemma \ref{bold:I}, there exists a $\Gamma_1'$ equivariant 
pseudodifferential operator 
$\mathbi{I}((\Tilde{A}_1)_{\lambda \lambda'}) \in \Psi(\Gamma_1\times
\Tilde{W}_1;
\Gamma_1\times(\Tilde{W}_1\times V))$ which is $\Gamma_1$ equivariant as well. 
An element $\gamma'\in\Gamma_1'$ acts on $\Gamma_1\times\Tilde{W}_1$ 
by right multiplication on $\Gamma_1$ by $\gamma'{}^{-1}$ and left 
multiplication on $\Tilde{W}_1$ by $\gamma'$.
Let $\mathbi{R}\mathbi{I}((\Tilde{A}_1)_{\lambda \lambda'})$ be the 
induced operator to the vector 
bundle of $\Gamma$  orbits $\Gamma\times(\Tilde{W}_1\times V)/\Gamma
\to \Gamma\times
\Tilde{W}_1/\Gamma$ which is isomorphic to $\Tilde{W}_1\times V\to 
\Tilde{W}_1$. Because
$\mathbi{I}((\Tilde{A}_1)_{\lambda \lambda'})$ was $\Gamma_1'$ equivariant, 
cf. remark \ref{gamma:prime},
the operator $\mathbi{R}\mathbi{I}((\Tilde{A}_1)_{\lambda \lambda'})$ 
is $\Gamma'$ equivariant.
We have the natural isomorphisms between the spaces of invariant sections
\begin{multline}
C^{\infty}(\Tilde{W}_1;\Tilde{W}_1\times V)^{\Gamma_1'}\overset{\sim}=
C^{\infty}(\Gamma_1\times\Tilde{W}_1;\Gamma_1\times(\Tilde{W}_1\times V))^
{\Gamma_1\times\Gamma_1'}\overset{\sim}=\\
\overset{\sim}=C^{\infty}(\Gamma_1\times_{\Gamma_1'}\Tilde{W}_1;\Gamma_1
\times_{\Gamma_1'}(\Tilde{W}_1\times V))^{\Gamma_1}
\end{multline}
which can be identified with the space of smooth sections in the vector 
orbibundle 
$\Tilde{E}_1/\Gamma_1 \overset{p_1}\to \Tilde{M}_1/\Gamma_1$ above the 
open set $\Tilde{W}_1/\Gamma_1'$.
The operators $\mathbi{R}\mathbi{I}((\Tilde{A}_1)_{\lambda \lambda ' })$, 
$\mathbi{I}((\Tilde{A}_1)_{\lambda \lambda ' })$
and $(\Tilde{A}_1)_{\lambda \lambda ' }$ restricted to the above spaces 
of invariant sections will be equal via the natural isomorphisms. 
In particular they induce $(A_1)_{\lambda \lambda '}$ at the level of sections 
in the associated vector orbibundle 
$\Tilde{E}_1/\Gamma_1 \overset{p_1}\to \Tilde{M}_1/\Gamma_1$.

We define the $\Gamma_2'$ equivariant operator
$\overline{A}_2=(L,l)_*\circ\mathbi{R}\mathbi{I}((\Tilde{A}_1)_{\lambda 
\lambda'})\circ
(L,l)^*$ in $\Psi(\Tilde{W}_2;\Tilde{W}_2\times V)$ using the isomorphisms 
$(L,l)$ provided in the diagram \eqref{local:diffs}.
(here, as before, $\Gamma_2'=(\Gamma_2)_{\tilde{x}_2}$). Let 
$\mathbi{I}(\overline{A}_2)\in\Psi(\Gamma_2\times\Tilde{W}_2;\Gamma_2\times
(\Tilde{W}_2\times V))$
be the $\Gamma_2$ equivariant operator which induces $\overline{A}_2$ at the
$\Gamma_2$ orbit level, and let $\mathbi{R}\mathbi{I}(\overline{A}_2)
\in\Psi(\Gamma_2\times_{\Gamma_2'}\Tilde{W}_2;\Gamma_2
\times_{\Gamma_2'}(\Tilde{W}_2\times V))$ be the 
operator induced at the $\Gamma_2'$ orbit level. Cf. remark \ref{gamma:prime}, 
this operator is $\Gamma_2$ equivariant. Using the sequence of isomorphisms:
\begin{multline}
C^{\infty}(\Tilde{W}_2;\Tilde{W}_2\times V)^{\Gamma_2'}\overset{\sim}=
C^{\infty}(\Gamma_2\times\Tilde{W}_2;\Gamma_1\times(\Tilde{W}_2\times V))^
{\Gamma_2\times\Gamma_2'}\overset{\sim}=\\
\overset{\sim}=C^{\infty}(\Gamma_2\times_{\Gamma_2'}\Tilde{W}_2;\Gamma_2
\times_{\Gamma_2'}(\Tilde{W}_2\times V))^{\Gamma_2}
\end{multline}
observe that the operators $\mathbi{R}\mathbi{I}(\overline{A}_2)$ and 
$\overline{A}_2$ act the same way on the 
corresponding spaces of invariant sections and induce $(A_2)_{\lambda 
\lambda '}=(T,t)_*\circ ({A}_1)_
{\lambda \lambda '}\circ (T,t)^*$ at the vector orbibundle level. If we extend 
$\mathbi{R}\mathbi{I}(\overline{A}_2)$
by $0$ outside the neighborhood $\Gamma_2\Tilde{W}_2$ we obtain a lift of 
$(A_2)_{\lambda \lambda '}$.
We proved the proposition in the case (ii) as well.
\end{proof}

We also have the uniqueness up to smoothing operators of the lift of a 
pseudodifferential operator.

\begin{prop} 
\label{unique}
Let $\Tilde{E}\overset{p}\to\Tilde{M}$ a $\Gamma$ vector bundle,
with $\Gamma$ a finite group and $E\overset{\overline{p}}\to M$ the 
associated orbibundle. If $\Tilde{A}_i$, $i=1,2$ are two classical 
pseudodifferential operators on $C^{\infty}(\Tilde{M};\Tilde{E})$ which are
$\Gamma$ equivariant and induce the same pseudodifferential operator on
$C^{\infty}(M;E)$, then $\Tilde{A}_1-\Tilde{A}_2$ is a smoothing operator.
\end{prop}
\begin{proof}
It is sufficient to prove that if $\Tilde{A}$ is a classical 
pseudodifferential operator that induces the zero operator on $C^{\infty}
(M;E)$, then $\Tilde{A}$ is smoothing. 

Let $\tilde{x}\in\Tilde{M}$ be a point such that its projection $\pi
(\tilde{x})$ onto
$M=\Tilde{M}/\Gamma$ is a smooth point. We will show that the total symbol
of the operator $\Tilde{A}$ on a neighborhood of $\tilde{x}$ is smoothing
(though the total symbol of the operator is not well-defined, the property 
of being a smoothing symbol is independent of the chart around $\tilde{x}$).
Let $\Tilde{W}$ be a slice at $\tilde{x}$ with respect to the action of the
group $\Gamma$. Because the point $x=\pi(\tilde{x})$ is smooth, the isotropy
group $(\Gamma)_{\tilde{x}}$ is trivial. Let $\phi$ be a smooth function
on $\Tilde{M}$ with support inside $\Tilde{W}$ which is equal to 1 
on a neighborhood of a smaller slice $\Tilde{V}$ at $\tilde{x}$.
Let $\psi$ be a smooth function on $\Tilde{M}$ which is equal to 1 on
a neighborhood of $supp(\phi)$ and zero outside $\Tilde{W}$.

$\phi\Tilde{A}\phi$ is a pseudodifferential operator which is 
localized on $\Tilde{W}$. We will prove that this operator is smoothing. 
Then $\Tilde{A}$ will be smoothing above $\Tilde{V}$ where $\phi=1$. 

Let $\tilde{\phi}=\sum_{\gamma\in\Gamma}\gamma\phi$ be the $\Gamma$ invariant 
extension of $\phi$ to $\Tilde{M}$. Then $\tilde{\phi}\Tilde{A}\tilde{\phi}$
is a $\Gamma$ equivariant extension to $C^{\infty}(\Tilde{M};\Tilde{E})$ 
of the zero operator on $C^{\infty}(M;E)$. Let $f\in C^{\infty}(\Tilde{M};
\Tilde{E})$ be an arbitrary section with support in a 
neighborhood of $supp(\psi)$ that vanishes outside the set $\{\tilde{x}|
\psi(\tilde{x})=1\}$.
Then $\tilde{f}=\sum_{\gamma\in\Gamma}\gamma f$ is a $\Gamma$ invariant 
section that extends $f$.
Because $\tilde{\phi}\Tilde{A}\tilde{\phi}$ restricted to the $\Gamma$
invariant section is equal to zero, we have:
\begin{align}
0&=\tilde{\phi}\Tilde{A}\tilde{\phi}(\tilde{f})=
\psi\tilde{\phi}\Tilde{A}\tilde{\phi}(\sum_{\gamma\in\Gamma}\gamma f)=\\
&=\psi\tilde{\phi}\Tilde{A}\tilde{\phi}\psi (\sum_{\gamma\in\Gamma}\gamma f)+
\psi\tilde{\phi}\Tilde{A}\tilde{\phi}(1-\psi)(\sum_{\gamma\in\Gamma}
\gamma f)=\\
&=\phi\Tilde{A}\phi(f)+(\sum_{\gamma\in\Gamma}\phi\Tilde{A}\tilde{\phi}(1-\psi)
\gamma)(f)=0
\end{align}
But $\tilde{\phi}(1-\psi)$ has support inside $\displaystyle\cup_{\gamma\ne e}
\gamma\Tilde{W}$ which is disjoint from $\Tilde{W}$ so $\phi\Tilde{A}
\tilde{\phi}
(1-\psi)$ and so $\sum_{\gamma\in\Gamma}\phi\Tilde{A}\tilde{\phi}(1-\psi)
\gamma$ are smoothing. Because $f$ was chosen arbitrary, we conclude that
$\phi\Tilde{A}\phi$ is smoothing so $\Tilde{A}$ is smoothing above $\Tilde{V}$.

Then $\Tilde{A}$ is smoothing on the regular set $\Tilde{M}_{\text{reg}}=
\{\tilde{x}|(\Gamma)_{\tilde{x}}=(e)\}$.
If $\tilde{x}_1\in\Tilde{M}_{\text{sing}}$ is a singular point then, on an
Euclidean chart around $\tilde{x}_1$ the total symbol of $\Tilde{A}$ will 
have the asymptotic expansion
\begin{equation}
a(\tilde{x},\xi)\sim\sum_{k\ge 0}a_{d-k}(\tilde{x},\xi)
\end{equation}
But $\Tilde{M}_{\text{reg}}$ is dense in $\Tilde{M}$ so $a_{d-k}(\tilde{x},
\xi)=0$ on a dense set because $a(\tilde{x},\xi)$ is smoothing on $\Tilde{M}_
{\text{reg}}$. We conclude that $a_{d-k}(\tilde{x},\xi)=0$ on a 
neighborhood of $\tilde{x}_1$ so $\Tilde{A}$ is smoothing on that neighborhood
and thus on the whole manifold $\Tilde{M}$.
\end{proof}

We denote by $\Psi(M;E)$ or simply by $\Psi$ the space of pseudodifferential 
operators acting on the smooth sections of the vector orbibundle 
$E\overset{p}\to M$. Most of the definitions and constructions for the 
pseudodifferential operators acting on genuine vector bundles can be extended
to the operators in $\Psi(M;E)$.

The symbol of a pseudodifferential operator $A$ in an linear vector
orbibundle chart $\mathcal{R}$ is a section in the vector orbibundle  
$\End(E)\overset{p}\to T^*(U)$ and is unique up to a smoothing symbol. 
 
The order of a pseudodifferential operator 
is the maximum of the orders of the restrictions to vector orbibundle charts
and denote by $\Psi^d$ the subspace of operators of order less or equal to $d$.
$\Psi$ becomes a filtered algebra with respect to the compositions of 
operators.

\begin{defn} A pseudodifferential operator $A$ of order $d$ is called elliptic 
if it is elliptic when restricted to all vector orbifold charts.
\end{defn}

An elliptic pseudodifferential operator $A\in \Psi^d$ has a parametrix 
$B\in\Psi^{-d}$  such that $AB-Id$ and $BA-Id$ are both smoothing. The 
construction of $B$ can be done locally, in each vector orbibundle chart, 
as in the case of pseudodifferential operators on manifolds.

In Proposition \ref{represent} we showed that any vector orbibundle 
$E\overset{p}\to M$ can be realized as a map between the orbit spaces 
of a $G$ vector bundle 
$\Tilde{E}\overset{\Tilde{p}}\to \Tilde{M}$, with $G$ a compact Lie group.
If $G$ was finite, Proposition \ref{pdounique} shows that any 
pseudodifferential operator acting on $C^{\infty}(M;E)$ can be realized 
as a $G$ equivariant pseudodifferential operator acting on 
$C^{\infty}(\Tilde{M};\Tilde{E})$. The following proposition extend 
this result to the case of Lie groups. 

\begin{prop}
\label{pdolift}
Let  $\Tilde{E}\overset{\Tilde{p}}\to \Tilde{M}$ be a $G$ vector bundle 
over a compact smooth manifold $M$ with $G$ a compact
Lie group such that any $x\in M$ has a finite isotropy group $G_x\subset G$ 
and denote by $E\overset{p}\to M$ the associated vector orbibundle of 
$G$ orbits. 
For any pseudodifferential operator $A$ acting on $C^{\infty}(M;E)$ there 
exists a
$G$ equivariant pseudodifferential operator $\Tilde{A}$ acting on 
$C^{\infty}(\Tilde{M};
\Tilde{E})$ of the same order as $A$ such that the restriction to the 
$G$ invariant sections
$C^{\infty}(\Tilde{M};\Tilde{E})^G \overset{\sim}=C^{\infty}(M;E)$ is 
equal to $A$.
If $A$ is classical then $\Tilde{A}$ can be chosen classical.
\end{prop}

\begin{proof}
As described in Proposition \ref{lievcano} one can find an open cover 
$(U_{\alpha})
_{\alpha}$ of the orbifold $M$, a family of finite groups 
$(G_{\alpha})_{\alpha}$ 
, $G_{\alpha}\subset G$, and a vector space $V$ such that $E_{|U_{\alpha}}
\overset{p}\to U_{\alpha}$ is isomorphic to the associated vector orbibundle 
of $G_{\alpha}$ orbits of the $G_{\alpha}$ vector bundle $\Tilde{U}_{\alpha}
\times V \overset{pr_1}\to\Tilde{U}_{\alpha}$ and $\Tilde{M}$ can be 
covered with a family of open sets $(N_{\alpha})_{\alpha}$ 
which are $G$ equivariant diffeomorphic to $G\times_{G_{\alpha}}
\!U_{\alpha}$ and such that the restriction of the $G$ vector bundle $\Tilde{E}
\overset{\Tilde{p}}\to \Tilde{M}$ to $N_{\alpha}$ is diffeomorphic with the
$G$ vector bundle $G\times_{G_{\alpha}}\!(U\times V)\overset{pr_1}\to G
\times_{G_{\alpha}}\!U$ by a $G$ equivariant vector bundle diffeomorphism.
The cover $(U_{\alpha})_{\alpha}$ can be chosen to be finite because 
$\Tilde{M}$ is compact. 
\par
Using a partition of unity, any pseudodifferential operator $A$ acting on 
sections 
of the vector orbibundle $E\overset{p}\to M$ can be written  as
$A=\sum_{\alpha}A_{\alpha}+K$ where $K$ is a smoothing operator with 
smooth kernel
which is zero on a neighborhood of the diagonal and $A_{\alpha}$ are
pseudodifferential operators which vanish on sections that are 
supported outside the 
open set $U_{\alpha}$ and take sections with support inside $U_{\alpha}$ 
into sections with support inside $U_{\alpha}$. 
For each $\alpha$ we will construct a $G$ equivariant pseudodifferential 
operator $\Tilde{A}_{\alpha}$ acting on sections of $\Tilde{E}
\overset{\Tilde{p}}
\to \Tilde{M}$ whose restriction to the $G$ invariant sections induces 
$A_{\alpha}$. 
\par
First let us fix a $G$ bi-invariant pseudodifferential operator $Q$ acting 
on $C^{\infty}(G)$ of order $d=ord(A)$ such that $Qf=0$ on constant 
functions $f$. 
In all the situations $Q$ can be taken to be the $d/2$ power of the 
Laplacian $\Delta$ on $G$ with respect to a bi-invariant metric on $G$.
Using Proposition \ref{pdounique} we can find a $G_{\alpha}$ equivariant 
pseudodifferential operator $B_{\alpha}$ acting on the space of sections 
$C^{\infty}(\Tilde{U}_{\alpha};V)$ of
the $G_{\alpha}$ bundle $\Tilde{U}_{\alpha}\times V\overset{pr_1}\to 
\Tilde{U}_{\alpha}$ which  induces $A_{\alpha}$ when restricted to 
$G_{\alpha}$ invariant sections. 
The space of smooth sections of the $G$ vector bundle $G\times
(\Tilde{U}_{\alpha}\times V)\overset{pr_1}\to G\times \Tilde{U}_{\alpha}$ 
is equal to $C^{\infty}(G\times \Tilde{U}_{\alpha}; V)\overset{\sim}=
C^{\infty}(G)\hat{\otimes}C^{\infty}(\Tilde{U}_{\alpha};V)$ so the 
pseudodifferential operator $Q\otimes 1+1\otimes 
B_{\alpha}$ acts on this space and it has order $d$. It is $G$ equivariant. 
The finite group $G_{\alpha}$ acts by a diagonal action on $G\times
(\Tilde{U}\times V)$ and $G\times\Tilde{U}$ so we can replace the above 
pseudodifferential operator with the $G_{\alpha}$ equivariant 
pseudodifferential operator
$\frac{1}{|G_{\alpha}|}\sum_{g\in G_{\alpha}} gQ\otimes 1 + 1
\otimes B_{\alpha}$.
This will induce a pseudodifferential operator on the sections
of the vector bundle of $G_{\alpha}$ orbits $G\times_{G_{\alpha}}\!(\Tilde{U}
_{\alpha}\times V)\overset{pr_1}\to G\times_{G_{\alpha}}\!\Tilde{U}_{\alpha}$.
We will transport it using the isomorphism between the previous $G$ 
vector bundle and 
$\Tilde{E}_{|N_{\alpha}}\overset{\Tilde{p}}\to N_{\alpha}$ to a 
pseudodifferential
operator $\Tilde{A}_{\alpha}$ acting on $\Tilde{E}\overset{\Tilde{p}}\to M$, 
extending it by $0$ outside $N_{\alpha}$. 
\par
We have to compare the actions of $A_{\alpha}$ and $\Tilde{A}_{\alpha}$.
Let $f$ be a $G$ invariant smooth section in $C^{\infty}(\Tilde{M};\Tilde{E})$
which induces $\overline{f}\in C^{\infty}(M;E)$. It is enough to take $f$
with support in $N_{\alpha}$. In this case $f$ is induced by a $G_{\alpha}$ 
invariant 
section $\Tilde{f}\in C^{\infty}(G\times \Tilde{U}_{\alpha}; G\times
(\Tilde{U}_{\alpha}\times V))\overset{\sim}=C^{\infty}(G)\hat{\otimes}
C^{\infty}
(\Tilde{U};V)$ which is $G$ invariant as well. Then $\Tilde{f}(g,u)=gh(u)$
with $h\in C^{\infty}(\Tilde{U};V)$, $h$ induces $\overline{f}\in C^{\infty}
(M;E)$ and $\Tilde{A}_{\alpha}f$ corresponds via the $G$ equivariant 
vector bundle diffeomorphism to $(Q\otimes 1 + 1\times B_{\alpha})
\Tilde{f}=B_{\alpha}(h)$ and this section induces $A_{\alpha}\overline{f}
\in C^{\infty}(M;E)$.
\par
We will construct a smoothing operator $\Tilde{K}$ in $C^{\infty}
(\Tilde{M};\Tilde{E})$ whose action on $G$ invariant sections induces the
smoothing operator $K$. The kernel of $K$ is given by a smooth section 
$S(x,y)$ of the endomorphism bundle $\End(E)\to M\times M$ which 
corresponds to a smooth $G$ invariant section $\Tilde{S}(\Tilde{x},\Tilde{y})$ 
of the endomorphism bundle $\End(\Tilde{E})\to\Tilde{M}\times\Tilde{M}$.
$\Tilde{S}$ is the kernel of a smoothing operator $\Tilde{K}$, which
induces $K$ when restricted to the $G$ invariant sections of $\Tilde{E}
\overset{p}\to \Tilde{M}$. 
\par
The operator $\Tilde{A}$ will be the sum of the operators 
$\Tilde{A}_{\alpha}$ and $\Tilde{K}$.
\par
If the operator $A$ is classical, we can choose $Q$ to be classical, then 
$B_{\alpha}$, $Q\otimes 1+1\otimes
B_{\alpha}$ and  $\Tilde{A}_{\alpha}$ are classical so $\Tilde{A}$ 
is classical.

\end{proof}

\subsection{Zeta Function of an Elliptic Pseudodifferential Operator}

Let $A$ be a pseudodifferential operator acting on the space of smooth sections
in a vector orbibundle $E\overset{p}\to M$. Suppose that $M$ is compact and it 
is endowed with a Riemannian metric $g$ and that we have a hermitian 
structure $<,>$ in the vector
orbibundle. As in the case of hermitian vector bundles over a closed Riemannian
manifold, we define the scalar product on $C^{\infty}(M;E)$ by
the formula
\begin{equation}
<f,g>=\int_M<f_1(x),f_2(x)>_x \,d\,vol \qquad\text{ for }f_1,f_2\in 
C^{\infty}(M;E).
\end{equation}
The integration with respect to the volume form $d\,vol$ can be defined
locally on vector orbibundle charts. If $(\Tilde{U}, \Gamma , U, \pi)$
is such a chart then the metric $g$ and the hermitian structure $<,>$ 
above $U$ are
induced by  $\Gamma$ invariant metric $\Tilde{g}$ and respectively a $\Gamma$ 
invariant hermitian structure $\tilde{<,>}$ in the $\Gamma$ vector 
bundle $\Tilde{U}\times V\overset{pr_1}\to\Tilde{U}$. Let $d\,\Tilde{vol}$ be
the volume form on $\Tilde{U}$ associated with $\Tilde{g}$. Then define
\begin{equation}
\int_U \alpha(x)\,d\,vol=\frac{1}{|\Gamma|}\int_{\Tilde{U}} \Tilde{\alpha}
(\Tilde{x})\,d\,\Tilde{vol}
\end{equation}
for $\Tilde{\alpha}\in C^{\infty}(\Tilde{U})^{\Gamma}$ that induces 
$\alpha\in C^{\infty}(U)$ (we chose $\alpha=<f_1,f_2>$). It is easy to see 
that the above formula is 
consistent with respect to the equivalence of charts.
The $L^2$ completion of $C^{\infty}(M;E)$ with respect to the above scalar 
product will be denoted by $L^2(M;E)$. A pseudodifferential operator $A$
defines an unbounded operator on $L^2(M;E)$.
\par
\begin{defn}
The pseudodifferential operator $A$ is called symmetric if
\begin{equation}
<Af_1,f_2>=<f_1,Af_2> \qquad\text{ for all } f_1,f_2\in C^{\infty}(M;E).
\end{equation}
A is positive if
\begin{equation}
<Af, f>\, \ge 0 \qquad\text{ for any } f\in C^{\infty}(M;E).
\end{equation}
\end{defn}

\begin{remark} If $A$ is symmetric, then the restriction $A_{\mathcal{R}}$ 
to any vector orbibundle chart $\mathcal{R}$ as defined in 
\ref{localize:op} is symmetric. This implies that $A_{\mathcal{R}}$ is 
induced by a $\Gamma$ equivariant pseudodifferential operator $\Tilde{A}_
{\mathcal{R}}$ which can be chosen to be symmetric with respect to the 
$\Gamma$ invariant metric $\Tilde{g}$ and $\Gamma$ invariant hermitian
structure inducing $g$ and $<,>$ above $U$. Indeed, the formal adjoint
$\Tilde{A}^*_{\mathcal{R}}$ of the operator $\Tilde{A}_{\mathcal{R}}$ 
induces $A_{\mathcal{R}}$ as well, so the symmetric operator $\frac{1}{2}
(\Tilde{A}_{\mathcal{R}}+\Tilde{A}^*_{\mathcal{R}})$ induces 
$A_{\mathcal{R}}$. 

Then the principal symbol of $A$, which is induced by the $\Gamma$ invariant
principal symbol of $\Tilde{A}$, is selfadjoint.
\end{remark}

\begin{thm} 
\label{selfadjoint}
Let $A$ be a positive symmetric elliptic pseudodifferential operator acting 
on the space of sections of a vector orbibundle $E\overset{p}\to M$ 
with compact base space. Then the operator $A$ acting on $L^2(M;E)$ 
is essentially selfadjoint and its spectrum is discrete.
\end{thm}

\begin{proof}
Let $m=dim(M)$ and denote by $G$ the orthogonal group $O(m)$.
Let $\Tilde{E}\overset{\Tilde{p}}\to \Tilde{M}$ be the $G$ vector bundle, 
with $\Tilde{M}=\mathcal{F}(M)$, 
such that the vector orbibundle of orbits is canonically isomorphic 
to $E\overset{p}\to M$, as in Proposition \ref{represent}.  Let us 
fix a $G$ bi-invariant metric $g_m$ on $G$ so that $\int_G\,dg_m=1$ and 
a $G$ bi-invariant positive selfadjoint pseudodifferential operator $Q$ 
of order $d$ acting on $C^{\infty}(G)$. As a metric $g_m$ we can choose 
the left translations of the opposite of the Killing form on the
Lie algebra $o(m)$ and $Q=\Delta^{d/2}$ where $\Delta$ is the Laplace
operator associated with the metric $g_m$.

We will construct a $G$ invariant metric on $\Tilde{M}$ which induces 
the given metric on $M$, a $G$ invariant hermitian 
structure on $\Tilde{E}\overset{\Tilde{p}}\to \Tilde{M}$ which 
induce the given hermitian metric on $E\overset{p}\to M$
and a $G$ equivariant elliptic positive selfadjoint pseudodifferential 
operator $\Tilde{A}$ acting on $C^{\infty}(\Tilde{M};\Tilde{E})$ whose 
restriction to the invariant sections is equal to $A$. 
\par
Let us fix a finite atlas $\mathcal{A}=(\mathcal{R}_{\alpha})_{\alpha}$ 
of the vector orbibundle $E\overset{p}\to M$, with
$\mathcal{R}_{\alpha}=(\Tilde{U}_{\alpha}, V, 
\Gamma_{\alpha}, U, \Pi_{\alpha}, \pi_{\alpha})$.
Then $\Tilde{\mathcal{R}}_{\alpha}=(G\times\Tilde{U}_{\alpha}, V, \Gamma
_{\alpha}, G\times_{\Gamma_{\alpha}}\Tilde{U}_{\alpha}, Id\times\Pi_{\alpha}, 
Id\times\pi_{\alpha})$ forms an atlas for $\Tilde{E}\overset{\Tilde{p}}\to 
\Tilde{M}$.
Let $\Tilde{g}_{\alpha}$ be the $\Gamma_{\alpha}$ 
invariant metric on $\Tilde{U}_{\alpha}$ which induces $g$ on $U$. Then the 
collection of $\Gamma_{\alpha}$ invariant metrics $g_m\otimes \Tilde{g}
_{\alpha}$ on $G\times\Tilde{U}_{\alpha}=\mathcal{F}(\Tilde{U})$ induces 
the $G$ invariant 
metric $\Tilde{g}$ on $\Tilde{M}$. Because $\Tilde{E}\overset{\Tilde{p}}
\to \Tilde{M}$ is the pull-back of $E\overset{p}\to M$ with respect 
to the canonical map $\Tilde{M}=\mathcal{F}(M)\overset{pr}\to M$,
we  define the hermitian structure on $\Tilde{E}\overset{\Tilde{p}}\to
\Tilde{M}$ to be the pull-back hermitian structure. It will be 
$G$ invariant by construction and induces the initial hermitian structure
when passing to the $G$ orbit spaces. Consider the $G$ invariant 
scalar product on 
$C^{\infty} (\Tilde{M};\Tilde{E})$ associated with hermitian structure and the 
metric $\Tilde{g}$ constructed above. When restricted to the $G$ invariant 
sections, this scalar product is equal with the scalar product on
$C^{\infty}(M;E)$ constructed with the help of the metric $g$ and the
hermitian structure on $E\overset{p}\to M$. Also $L^2(\Tilde{M};\Tilde{E})^G
\overset{\sim}=L^2(M;E)$.

The pseudodifferential operator $\Tilde{A}$ acting on $C^{\infty}(\Tilde{M};
\Tilde{E})$ constructed as in Proposition \ref{pdolift} will be 
symmetric with respect to the scalar product on $C^{\infty}(\Tilde{M};
\Tilde{E})$. By construction, the principal symbol of $\Tilde{A}$ is equal to 
$q_{pr}\otimes Id+a_{pr}$ where $q_{pr}$ is the principal symbol of $Q$ and
$a_{pr}$ is the principal symbol of $A$. Because $Q$ and $A$ are positive 
symmetric elliptic pseudodifferential operators, the operator $\Tilde{A}$
is elliptic as well. 
Because $\Tilde{M}$ is compact, the operator $\Tilde{A}$
acting on $L^2(\Tilde{M};\Tilde{E})$ is essentially selfadjoint. 
The restriction of $\Tilde{A}$ to the closed subspace of $G$ invariant 
$L^2$ sections $L^2(\Tilde{M};\Tilde{E})^G\overset{\sim}=L^2(M;E)$ which 
is equal to $A$ is essentially selfadjoint. The spectrum of $\Tilde{A}$ is 
real and discrete so the spectrum of $A$ is a discrete subset of $\RR$. 
\end{proof}

To show that the complex powers of the elliptic positive selfadjoint 
operator $A$ are 
pseudodifferential operators, one can apply the same local approach as in 
\cite{seeley}, and use approximations of the resolvent of $(A-\lambda)$
in order to show that $A^s$ are, up to a smoothing operator, pseudodifferential
operators of complex order $s$. Because all the proofs follow in the same 
way with minimal complications, we will not repeat them in our paper.  
For $Re(s)<\frac{-m}{d}$ the distributional 
kernel of $A^s$ is continuous and $A^s$ is of trace class.

We can relax the conditions imposed on the operator $A$, and prove an 
analogous 
statement to Theorem \ref{selfadjoint}. Suppose that $\pi$ is an Agmon 
angle for $A$. Then one can define the complex powers of $A$ as
\begin{equation}
\label{pow:neg}
A^s_{\pi}=\frac{1}{2\pi i}\int_{\wp}\lambda^s(\lambda - A)^{-1} \,
d\lambda \qquad \text{when }Re(s) < 0
\end{equation}
(where $\wp$ is a contour in the complex plane obtained by joining
two parallel half-lines to the negative real axis by a circle around
the origin) and
\begin{equation}
\label{pow:pos}
A^s_{\pi}=A^{s-k}_{\pi}A^k_{\pi} \qquad \text{for } Re(s)\ge 0
\end{equation}
for large enough $k\in \mathbb Z$ that makes $s-k < 0$.
Obviously, if $A$ is positive selfadjoint then $\pi$ is an Agmon angle for $A$.

\begin{prop}
\label{Agmon}
Let $A$ be an elliptic pseudodifferential operator of positive order $d$ 
acting on the space of sections of a vector orbibundle 
$E\overset{p}\to M$ with compact 
base space. Suppose that $\pi$ is an Agmon angle for $A$.
Then the operator $A$ acting on $L^2(M;E)$ has a discrete spectrum
and its complex powers $A^s_{\pi}$ are pseudodifferential operators of 
complex order $sd$.
\end{prop}
\begin{proof}
The proof is analogous to the proof of Theorem \ref{selfadjoint}
We construct the operator $\Tilde{A}$ acting on sections of the $G$ 
vector bundle 
$\Tilde{E}\overset {\Tilde{p}}\to \Tilde{M}$ with the help of a $G$ 
bi-invariant elliptic selfadjoint pseudodifferential operator $Q$ acting 
on $C^{\infty}(G)$. 
The principal symbol of $\Tilde{A}$ is equal to $\Tilde{a}_{pr}=q_{pr}\otimes 
Id+ a_{pr}$ with $q_{pr}$ the principal symbol of $Q$ and $a_{pr}$ 
the principal symbol of $A$.
Because $\pi$ is an Agmon angle for $A$, the spectrum of $a_{pr}(x,\xi)$ is
disjoint from the region in the complex plane 
$\mathcal{C}=\{z\,|\, arg(z)\in(\pi-\varepsilon, \pi+\varepsilon)\}
\cup\{z\,|\,|z|<\varepsilon\}$ for a small enough $\varepsilon>0$. Because
$q_{pr}(x,\xi)$ is selfadjoint the spectrum
of $\Tilde{a}(x,\xi)$ is disjoint from $\mathcal{C}$. We conclude that $\pi$ is
an Agmon angle for $\Tilde{A}$. Then $\Tilde{A}$ is elliptic and has a discrete
spectrum. $A$ being the restriction to the $G$ invariant sections will have a 
discrete spectrum as well.
The complex powers $\Tilde{A}^s_{\pi}$ are well
defined and are pseudodifferential operators of order $sd$ so the restriction 
to the $G$ invariant sections, which are equal to $A^s_{\pi}$, are 
pseudodifferential operators of order $sd$. 
\end{proof}

Throughout the rest of the chapter we suppose that $\pi$ is an Agmon
angle for the elliptic pseudodifferential operator $A$. 
The complex powers $A^s=A^s_{\pi}$ will be defined as in 
\eqref{pow:neg} and \eqref{pow:pos}.

\begin{defn} The zeta function of the operator $A$ is equal to
\begin{equation} 
\zeta_A(s)=\Tr (A^s) \qquad\text{ for }Re(s)<\frac{-m}{d}.
\end{equation}
\end{defn}

We will prove the following theorem:

\begin{thm} 
\label{zeta:Dirac:density}
The zeta function of an elliptic pseudodifferential operator $A$
acting on the sections of a vector orbibundle $E\overset{p}\to M$ 
can be extended to a meromorphic function on $\CC$ with at most
simple poles at $s=\frac{-m+k}{d}$, $k=0,1,2,\dots$. The residues are 
integrals on $M$ of Dirac-type densities which can be explicitly computed 
in terms of the asymptotic expansion of the total symbol of the operator $A$.
\end{thm}

\begin{proof}
Let $\mathcal{A}=(\Tilde{U}_{\alpha}, V,\Gamma_{\alpha},U_{\alpha},
\Pi_{\alpha},\pi_{\alpha})_{\alpha}$ be a fixed finite atlas of linear 
orbibundle charts . Using a partition of unity we can decompose the 
complex powers of the operator $A$ as
\begin{equation}
A^s=\sum_{\alpha}A^s_{\alpha} + K_s
\end{equation}
where $A^s_{\alpha}$ are pseudodifferential operators of order $sd$ 
with support inside 
$U_{\alpha}$ and $K_s$ is a smoothing operator. We can arrange for 
$A^s_{\alpha}$ to take smooth sections with support inside $U_{\alpha}$ 
into smooth sections with support inside $U_{\alpha}$. 
The trace $Tr(K_s)$ is defined for $s\in \CC$ and it is a holomorphic 
function on the complex plane. In order to prove the theorem, we will 
show that the statement in the theorem holds for the trace of the 
operators $A^s_{\alpha}$. 
We will show that $Tr(A^s_{\alpha})$ is a meromorphic function on $\CC$ with
at most simple poles at $s=\frac{-m+k}{d}$, $k=1,2,\dots$, and that 
the residues are computed as integrals on $M$ of Dirac-type densities. 
\par
$A^s_{\alpha}$ can be seen as a holomorphic family of pseudodifferential 
operators of order $sd$ acting on the smooth sections of the vector orbibundle 
$E_{|U_{\alpha}}\overset{p}\to U_{\alpha}$. For the sake of simplicity, 
we will drop the index $\alpha$ and denote $A^s_{\alpha}$ by $A_s$. There 
exists a holomorphic family of $\Gamma$ equivariant pseudodifferential 
operators $\Tilde{A}_s$ of order $sd$ acting on the sections of the 
$\Gamma$ vector bundle $\Tilde{U}\times V
\overset{pr_1}\to \Tilde{U}$ which induce $A_s$ on the $\Gamma$ 
invariant sections.
This family is unique up to smoothing operators, cf. Proposition \ref{unique}.
Then the component corresponding to the trivial representation of the trace 
of $\Tilde{A}_s$ and of the residues will be equal to respectively the
trace functional of $A_s$ and its residues. 
Theorem \ref{residues:i} shows that the component of the trivial 
representation 
of trace functional is a meromorphic function with at most simple 
poles at $s=\frac{-m+k}{d}$ for $k=0,1,2,\dots$. From Theorem 
\ref{residues:d}, the component of the residue of the trace 
functional $Tr(\Tilde{A}_s)$ corresponding to the trivial irreducible 
representation is:
\begin{equation}
\label{residues:o}
res_{|s=\frac{-m+k}{d}}Tr(\Tilde{A}_{s,\chi_0})=\frac{1}{|\Gamma|}
\sum_{\gamma\in\Gamma}\int_{\Tilde{U}^{\gamma}}\eta_k^{\gamma}
\end{equation}
where $\{\eta_k^{\gamma}\}_{\gamma\in\Gamma}$ are $\Gamma$ Dirac-type 
densities on $\Tilde{U}$ constructed from the total symbol of $\Tilde{A}_s$ 
which on the chosen linear vector orbibundle chart coincides with the 
total symbol of the pseudodifferential operator $A_s$ on $U$.

Though the statements of Theorems \ref{residues:i} and \ref{residues:d}
refer to the complex powers of a pseudodifferential operator, their proofs can 
be slightly changed and adapted to include the case of holomorphic families
of pseudodifferential operators and so the conclusions listed above are 
still true. 

The $\Gamma$ Dirac-type density $\{\eta_k^\gamma\}_{\gamma\in\Gamma}$ 
defines a Dirac-type density on the orbifold $U$ which we will denote 
by $\eta_k$. Observe that the formula \eqref{residues:o} is equivalent to
\begin{equation}
\res_{|s=\frac{-m+k}{d}}Tr(A_s)=\int_U \eta_k
\end{equation}

The Dirac-type density whose integral on $M$ is equal to the residues of 
the zeta function  $\zeta_A$ at $s=\frac{-m+k}{d}$ will be equal to 
the sum of the Dirac-type densities $\eta_k$ constructed as above on the
charts of the vector orbibundle.

\end{proof}

Using Proposition \ref{eta:strata}, we can reformulate 
Theorem \ref{zeta:Dirac:density} and express the residues of the zeta 
function of $A$ in terms of integrals of smooth densities on submanifolds of 
$\mathcal{S}(M)$ associated with the canonical stratification.

\begin{thm} 
\label{zeta:strata:density}
The zeta function of an elliptic pseudodifferential operator $A$
acting on the sections of a vector orbibundle $E\overset{p}\to M$ 
can be extended to a meromorphic function on $\CC$ with at most
simple poles at $s=\frac{-m+k}{d}$, $k=0,1,2,\dots$. 
For each $k\ge 0$ there exist densities $\eta_{\upsilon,k}$ on the strata 
$M_{\upsilon}$ of the canonical 
stratification of $M$ such that the residue of the zeta function
at $s=\frac{-m+k}{d}$ is the sum of
integrals of these densities.
On each strata $M_{\upsilon}$ the density $\eta_{\upsilon,k}$ can be 
explicitly computed in terms of the asymptotic expansion of total symbol 
of the operator $A$ in a neighborhood of $M_{\upsilon}$.
\end{thm}

\begin{proof}
For each $\upsilon\in\mathcal{O}(M)$ and $k\ge 0$, the density 
$\eta_{\upsilon,k}$
can be obtained as follows: let $x\in M_{\upsilon}$ and $(\Tilde{U}, V, 
\Gamma, U, \Pi, \pi)$ be a linear vector orbibundle chart centered at $x$.
Observe that $\Gamma\in\upsilon$ and the map $\pi: \Tilde{U}^{\Gamma} 
\to U\cap M_{\upsilon}$ is a local diffeomorphism.
Theorem \ref{zeta:Dirac:density} gives us Dirac-type densities $\eta_k$ on $M$ 
whose integrals are equal to the residues of the zeta function and can 
be computed explicitly in terms of the total symbol of $A$. These 
densities can be described above $U$ as  families of densities 
$\{\eta^{\gamma}_k\}$ such that $\eta^{\gamma}_k$ is a smooth density on 
$\Tilde{U}^{\gamma}$.  
Consider the density on $\Tilde{U}^{\Gamma}$ given by
\begin{equation}
\eta_{\Gamma,k}=\frac{1}{|\Gamma|} \sum_{\substack{\gamma\in\Gamma \\ 
dim(\Tilde{U}^{\Gamma})=dim(\Tilde{U}^{\gamma})}} \eta^{\gamma}_k
{}_{|\Tilde{U}^{\Gamma}}
\end{equation}
Define $\eta_{\upsilon,k}$ on $M_{\upsilon}$ to be the density whose 
restriction to $U$  is equal to $\pi_*(\eta_{\Gamma,k})$.

As proved in Proposition \ref{eta:strata}, the integral on $M$ of the 
Dirac-type density  $\eta_k$ and the sum of the integrals of the 
densities  $\eta_{\upsilon,k}$ on the strata 
$M_{\upsilon}$ for $\upsilon\in\mathcal{O}(M)$  are equal, so 
they are equal to the residue of $\zeta_A$ at $s=\frac{-m+k}{d}$.

\end{proof}

\vfill\eject


\appendix
\section{Proof of Theorem \ref{local}}

\begin{proof}
Let $U$ be a $\Gamma_x$-invariant open neighborhood of $x$ in $M$ such that 
$\gamma U\cap U=\emptyset$ for $\gamma \in \Gamma\backslash \Gamma_x$. Let 
$\mu':\Gamma_x
\times T_x(M)\to T_x(M)$ be the linear representation of $\Gamma_x$ on
the tangent space at $x$ generated by the action of $\Gamma_x$ on $U$.
Choose a $\Gamma_x$ invariant metric on $U$ (this can always be done by averaging 
a metric over the finite group $\Gamma_x$). Then the exponential map 
$exp:T_x(M)\to U$ is $\Gamma_x$ equivariant. Indeed, if $s(t)$ is a geodesic curve
with $s(0)=x$ and $s'(0)=X\in T_x(M)$ then, for $\gamma \in \Gamma_x$, $\gamma 
\cdot s(t)$ is again
a geodesic curve and $(\gamma \cdot s)'(0)=\mu'(\gamma )X$. So $exp(
\mu'(\gamma )X)
=\gamma \cdot s(1)=\gamma \cdot exp(X)$. One can choose a smaller $U$ such 
that the 
exponential map realizes a diffeomorphism between a $\Gamma_x$ invariant 
neighborhood of 
$0$ in $T_x(M)$ and $U$. If we fix a linear isomorphism $T_x(M)\overset
{\sim}=\mathbb{R}^m$ and transport the action of $\Gamma_x$ on $T_x(M)$  
via this isomorphism, we obtain a linear representation $\mu:\Gamma_x\times
\mathbb{R}^m\to\mathbb{R}^m$ and a $\Gamma_x$ equivariant diffeomorphism between
a neighborhood $O$ of the origin and $U$. 
\par
Using the pull-back via this diffeomorphism  we can replace the bundle 
$E_{|U}\overset{p}\to U$ by a new $\Gamma_x$ vector bundle such that the action of 
$\Gamma_x$ on the base will be the restriction of a linear representation to
an open neighborhood of the origin. The point $x$ corresponds to the origin 
$0$. For the sake of simplicity we will use the same notation for the new bundle. 
\par
Consider now a linear connection $\nabla$ in the $\Gamma$ vector bundle 
$E_{|U}\overset{p}\to U$. The group $\Gamma_x$ acts on the space of connections, and 
we can make $\nabla$ to be invariant with respect to $\Gamma_x$ by replacing it 
with the connection $\nabla+\frac 1{|\Gamma_x|}\sum_{\gamma \in \Gamma_x}
(\gamma _*\nabla-\nabla)$.  We will define
the bundle map $\Psi$ between the trivial $\Gamma$ bundle $U\times E_x \to U$ and
$E_{|U}\overset{p}\to U$ as follows:\newline
for a point $y\in U$ and a vector $X\in E_x$, let $\Psi(y,X)=Y$ be the vector 
in the fiber of $E$ above $y$ obtained by parallel transport of $X$ along 
the curve $\phi(t)=ty,\; t\in[0,1]$. The vector space $V$ as in the statement 
of the proposition will be equal to $E_x$ with the action of $\Gamma_x$ induced 
from the action of $\Gamma$ on the total space of the bundle.

$\Psi$ is a vector bundle isomorphism and we have to show that $\Psi$
is $\Gamma_x$ equivariant. If $\Tilde{\phi}(t)$ is the path in $E$ that 
realizes the parallel transport from $X$ to $Y=\Psi(y,X)$ above $\phi(t)$, 
and $\gamma \in \Gamma_x$, then, because the connection is $\Gamma_x$ invariant, 
the path $\gamma \cdot\Tilde{\phi}(t)$ realizes the parallel transport between 
$\gamma X$ and $\gamma Y$ above the curve $p(\gamma \cdot\Tilde{\phi}(t))=
\gamma\cdot\phi(t)=\gamma (ty)=t(\gamma y)$ (we used the fact that the action of 
$\Gamma_x$ on $U$ is linear). Then $\gamma Y=\Psi(\gamma y,\gamma X)$, so $\Psi$ 
is $\Gamma_x$ equivariant.

For $U\subset M$ as above 
the set $\Gamma\cdot U=\bigcup_{\gamma\in\Gamma}\gamma U\subset M$ is the 
reunion of disjoint open sets $\gamma_i U$ with $\{\gamma_1, \gamma_2, \dots\}
$ a complete system of left coset representatives for $\Gamma/\Gamma_x$. Then
$\Gamma\times_{\Gamma_x}U\overset{u}\to \Gamma\cdot U$ given by 
$u(\gamma,y)=\gamma y$ is a diffeomorphism, where $\Gamma\times_{\Gamma_x}U$ 
is the cross-product $\Gamma\times U/\sim$ , with $(\gamma \gamma',y)\sim
(\gamma,\gamma'y)$, for any $\gamma\in \Gamma$, $\gamma'\in\Gamma_x$ and $y
\in U$. Moreover, $u$ is $\Gamma$ equivariant. The action of $\Gamma$ on the  
cross-product $\Gamma\times_{\Gamma_x}U$ is by left translations. The 
restriction of the vector bundle $E_{|\Gamma U}\overset{p}\to \Gamma \cdot U$ is 
trivial and a trivialization is given by  
$$
\xymatrix{
\Gamma\times_{\Gamma_x}(O\times E_x)\ar[d]_{pr_1}\ar[r]^{\qquad \mathbf{1}\times\Psi}
& E_{|\Gamma \cdot U}\ar[d]^p\\
\Gamma\times_{\Gamma_x}O\ar[r]^{\quad \mathbf{1}\times\psi}&\Gamma \cdot U
}
$$
The bundle isomorphism is also $\Gamma$ equivariant.
\end{proof}

\section{Proof of Theorem \ref{residues:g}}

\begin{proof}
Let $\{U_{\alpha}\}$ be a finite cover of $M$ with open sets as in Proposition \ref
{local} so that the restrictions of the $\Gamma$ bundle $E\overset{p}\to M$ 
to the subsets $U_{\alpha}$ are isomorphic to $\Gamma_x$ trivial vector bundles
for some $x\in U_{\alpha}$. 
Then $\{\Gamma\cdot U_{\alpha}\}$ is a 
finite open cover with $\Gamma$ invariant open sets. Consider a partition of 
unity $\{\phi_{\beta}\}$ subordinated 
to the open cover $\{\Gamma\cdot U_{\alpha}\}$ such that the functions 
$\phi_{\beta}$ are $\Gamma$ equivariant. One can choose the open cover and the 
partition of unity such that for any $\beta$ and $\beta'$ either $supp(\phi_{\beta})
\cap supp(\phi_{\beta'})=\emptyset$ or $supp(\phi_{\beta})\cap supp(\phi_{\beta'})
\subset \Gamma\cdot U_{\alpha}$ for some $\alpha$. Then the operators 
$\phi_{\beta}\cdot A^s \cdot \phi_{\beta'}$ are either smoothing (in the first case)
or have the support and range included in the space of sections that vanish 
outside the open set $\Gamma \cdot U_{\alpha}$ (in the second case). Consider the 
trivialization by the $\Gamma$ vector bundle isomorphism $(\Gamma\times_{\Gamma_x}
(O_\alpha\times E_x), \, \Gamma\times_{\Gamma_x}O_{\alpha})\overset{\sim}\to
(E_{|\Gamma\cdot U_{\alpha}},\,\Gamma\cdot U_{\alpha})$. Let $a_s(x, \xi)$ the 
complete symbol of the operator $\phi_{\beta}\cdot A^s \cdot \phi_{\beta'}$ in this
trivialization. Then for a section $f$ with compact support in 
$\Gamma\cdot O_{\alpha}$ we have:
\begin{equation}
(\phi_{\beta}\cdot A^s \cdot \phi_{\beta'})(f)(x)=\iint e^{i<x-y,\xi>}a_s(x,\xi)f(y)
dy\, \dbar\xi
\end{equation}
where $\dbar\xi=(2\pi)^{-m}d\xi$ and the double integral is computed over 
$\RR^m\times\RR^m$.
To keep our notation simple, we will drop the indexes $\alpha$ and $\beta$ and
denote $\phi_{\beta}\cdot A^s \cdot \phi_{\beta'}$ by $A_s$.   
The action of $\mathcal{T}$ on $f$ is given by $\mathcal{T}f(x)= Tf
(t^{-1}x)$. Then
\begin{align}
(A_s\circ\mathcal{T})f(x)&=\iint e^{i<x-y,\xi>}a_s(x,\xi)\circ T  
f(t^{-1}y) dy\,\dbar\xi\\
\intertext{and after the change of variable $y\mapsto ty$}
&=\iint e^{i<x-ty,\xi>}a_s(x,\xi)\circ T  f(y)|det(t)|\,dy\,\dbar\xi
\end{align}

We can choose the trivialization in such a way that 
$|det(t)|=1$.
The distributional kernel of $A_s\circ\mathcal{T}$ will be equal to $K_s(x,y)=
\int e^{i<x-ty,\xi>}a_s(x,\xi)\circ T  \dbar\xi$. For $Re(s)< -\frac{m}{d}$ the 
trace of this operator is equal to the integral 
\begin{equation}
\label{trsT}
Tr(A_s\circ T )=\iint e^{i<x-tx,\xi>}Tr(a_s(x,\xi)\circ T )\,\dbar\xi
\,dx
\end{equation}

We will choose the open cover $O_{\alpha}$ so that either $t$ has fixed points
inside $O_{\alpha}$ or $\|x-tx\|>\varepsilon$ for all $x\in O_{\alpha}$ for some 
fixed $\varepsilon>0$. In 
the second case we will show that the trace function $Tr(A_s\circ \mathcal{T})$
can be extended to the whole complex plane. Indeed, if we denote by $\|D_{\xi}\|^2$ 
the differential operator $-(\frac{\partial^2}{\partial \xi_1^2}+\frac{\partial^2}
{\partial \xi_2^2}+\dots+\frac{\partial^2}{\partial \xi_m^2})$ we have 
$\|D_{\xi}\|^{2\nu}e^{i<x-tx,\xi>}=\|x-tx\|^{2\nu}e^{i<x-tx,\xi>}$. Then
\begin{align}
Tr(A_s\circ \mathcal{T})&=\iint e^{i<x-tx,\xi>}Tr(a_s(x,\xi)\circ T )\,\dbar\xi
\,dx=\\
&=\iint \frac{\|D_{\xi}\|^{2\nu}e^{i<x-tx,\xi>}}{\|x-tx\|^{2\nu}}Tr(a_s(x,\xi)\circ
 T )\,\dbar\xi\,dx \\
\intertext{and after integration by parts}
&\iint \frac{e^{i<x-tx,\xi>}}{\|x-tx\|^{2\nu}}\|D_{\xi}\|^{2\nu}Tr(a_s(x,\xi)\circ
 T )\,\dbar\xi\,dx
\end{align}
For a fixed half-plane $Re(s)<K$ if we choose a large enough $\nu\in\mathbb{N}$
the expression $\|D_{\xi}\|^{2\nu}(Tr(a_s(x,\xi)\circ T )$ is a symbol in 
$S^{-2m}(O)$ so it is absolutely integrable  and gives after 
integration a holomorphic function in $s$ in the half-plane $Re(s)<K$. It follows 
that on open sets where $\|x-tx\|>\varepsilon$ the trace function $Tr(A_s\circ 
\mathcal{T})$ has a holomorphic extension to the complex plane $\mathbb{C}$.

Consider now a set $O=O_\alpha$ so that $t_{|O}$ has a nonempty fixed point set.
The diffeomorphism $t$ is given by $\gamma$. Denote $N=M^{\gamma}$.
Let $x\in O$ such that $\gamma\in \Gamma_x$. As shown in Proposition 
\ref{local}, we can choose $x$ to be the origin in $O\subset\mathbb{R}^m$ and
$\Gamma_x$ act on $O$ by linear isometries. Then $t$ is the restriction of a 
linear isometry and its fixed set $N\cap O$ is the intersection between a linear 
subspace $F_t$ of $\mathbb{R}^m$ and $O$. The dimension of $F_t$ is one of the 
dimensions $n_i$ of the connected components of the fixed point set $N$. 
For the sake of simplicity we 
will denote it by $n$. Let $(x_1,x_2)$ be linear coordinates on $\mathbb{R}^m$ so 
that $x_1$ are coordinates along the fixed-point set $F_t$ and
$x_2$ are normal coordinates to $F_t$. In these new coordinates $t$ has the form 
$\bigl(\begin{smallmatrix} Id&0 \\ 0&\overline{t} \end{smallmatrix}\bigr)$ with 
$\overline{t}$ an $(m-n)$ square matrix whose 
eigenvalues are different from $1$. Let $\xi=(\xi_1,\xi_2)$ be the coordinates in 
the cotangent space corresponding to the coordinates $x=(x_1,x_2)$. Using the new 
coordinates in the equation \eqref{trsT} $Tr(A_s\circ\mathcal{T})$ will be equal to
\begin{align}
&\iint e^{i<(x_1,x_2)-(x_1,\overline{t}x_2),(\xi_1,\xi_2)>}
Tr(a_s(x_1,x_2,\xi_1,\xi_2)\circ T )\,\dbar\xi_1 \dbar\xi_2 dx_1 dx_2= \\
&=\iint e^{i<(Id-\overline{t})x_2,\xi_2>}Tr(a_s(x_1,x_2,\xi_1,\xi_2)\circ T )
\dbar\xi_1 \dbar\xi_2 dx_1 dx_2
\end{align}
Because $\overline{t}$ is a matrix whose 
eigenvalues are different from $1$, the matrix $(\overline{t}-Id)$ is non-degenerated and
if we use the change of coordinates $w=(\overline{t}-Id)x_2$ the above integral becomes
\begin{align}
&\iint e^{-i<w,\xi_2>}Tr(a_s(x_1,(\overline{t}-Id)^{-1}w,\xi_1,\xi_2)\circ T )
|det(\overline{t}-Id)|^{-1}\,\dbar\xi_1 \dbar\xi_2 dx_1 dw=\\
&=\iint e^{-i<w,\xi_2>}Tr(\overline{a}_s(x_1,w,\xi_1,\xi_2)\circ T )
|det(\overline{t}-Id)|^{-1}\,\dbar\xi_1 \dbar\xi_2 dx_1 dw=\\
&=Tr\left[\left(\iint e^{-i<w,\xi_2>}\overline{a}_s(x_1,w,\xi_1,\xi_2)
\,\dbar\xi_1 \dbar\xi_2 dx_1 dw\right)\circ  T \right]|det(\overline{t}-Id)|^{-1}
\end{align}
where $\overline{a}_s(x_1,w,\xi_1,\xi_2)=a_s(x_1,(\overline{t}-Id)^{-1}w,\xi_1,
\xi_2)$ is a classical symbol of order $sd$ in $(\xi_1,\xi_2)$ and with 
compact support in $(x_1,w)$. For $Re(s)<-\frac m d$ the integrand is absolutely 
integrable so we can apply Fubini's theorem and get:
\begin{equation}
Tr\left[\left( \iint(\iint e^{-i<w,\xi_2>}\overline{a}_s(x_1,w,\xi_1,\xi_2)
\,\dbar\xi_2 dw) \dbar\xi_1 dx_1\right)\circ  T \right]|det(\overline{t}-Id)|^{-1}
\end{equation}
We will need the following result:
\begin{lemma}
\label{reduce}
If $a_s(x,\xi)=a_s(x_1,x_2,\xi_1,\xi_2)$ is a classical matrix valued symbol of 
complex order $sd$ then, for $Re(s)<-\frac{m-n}{d}$, the expression 
\begin{equation}
\label{defn:b}
b_s(x_1,\xi_1)=\iint e^{-i<x_2,\xi_2>}a_s(x_1,x_2,\xi_1,\xi_2)\,\dbar\xi_2\,dx_2
\end{equation}
is a classical symbol of order $sd$ with the asymptotic expansion
\begin{equation}
\label{defn:b:asympt}
b_s(x_1,\xi_1)\sim\sum_{\alpha}\frac{1}{\alpha !}(D^{\alpha}_{x_2} \partial^
{\alpha}_{\xi_2}a_s)(x_1,0,\xi_1,0)
\end{equation}
If $a_s$ is a holomorphic family of symbols then there exists a holomorphic
family of symbols $\Tilde{b}_s$ defined for $s\in \CC$ such that $\Tilde{b}_s
- b_s$ is a holomorphic family of smoothing symbols.
The asymptotic expansion of $\Tilde{b}_s(x_1,\xi_1)$ is given by 
\eqref{defn:b:asympt}.
\end{lemma}
We will postpone the proof until the end of the proof of the 
main theorem.

Using the above lemma for $Re(s)<-\frac{(m-n)}{d}$, the trace of 
$A_s\circ \mathcal{T}$ becomes:
\begin{equation}
Tr(A_s\circ\mathcal{T})=Tr\left(\iint b_s(x_1,\xi_1)\,\dbar\xi_1 dx_1\circ 
 T \right)|det(\overline{t}-Id)|^{-1}
\end{equation}
where the integration is taken over the cotangent space of the fixed point 
set $N\cap O$. Because $\Tilde{b}_s-b_s$ is a smoothing symbol, the difference 
between the integral above and
\begin{equation}
Tr\left(\iint \Tilde{b}_s(x_1,\xi_1)\,\dbar\xi_1 dx_1\circ 
 T \right)|det(\overline{t}-Id)|^{-1}
\end{equation}
is a whole function. We will show that, as a function in $s\in \mathbb{C}$, the 
above expression has a meromorphic extension to the complex plane.    
Let us consider the asymptotic expansion as a classical symbol 
\begin{equation}
\Tilde{b}_s(x_1,\xi_1)\sim\sum_{i \ge 0}\Tilde{b}_{s,i}(x_1,\xi_1)
\end{equation}
with $\Tilde{b}_{s,i}(x_1,\xi_1)$ homogeneous of  degree $sd-i$ in $\xi_1$. Let 
$\psi(\xi_1)$ be a positive, real valued smooth function which vanishes on a 
neighborhood of the origin and is equal to $1$ for $\|\xi_1\|\ge 1$.
Let us fix 
a half-plane $Re(s)<K$. For a large enough $\nu$ the difference  
$\Tilde{b}_s(x_1,\xi_1)-\sum_{i=o}^{\nu} \psi(\xi_1)\Tilde{b}_{s,i} (x_1,\xi_1)$ 
is a symbol in 
$S^{-n}(N\cap O)$ for any $s$ with $Re(s)<K$. Then on this half-plane the 
difference 
\begin{equation}
\begin{split}
&Tr\left(\iint \Tilde{b}_s(x_1,\xi_1)\,\dbar\xi_1 dx_1\circ T \right)
|det(\overline{t}-Id)|^{-1}-\\
&\qquad-\sum_{i=0}^{\nu}Tr\left(\iint \psi(\xi_1)\Tilde{b}_{s,i}
(x_1,\xi_1)\,\dbar\xi_1 dx_1\circ T \right)|det(\overline{t}-Id)|^{-1}
\end{split}
\end{equation}
is holomorphic. We will have to prove the existence of a meromorphic extension
and find  the poles and residues for each of the functions 
\begin{equation}
s\mapsto Tr\left(\iint \psi(\xi_1)\Tilde{b}_{s,i}(x_1,\xi_1)\,\dbar\xi_1 dx_1\circ T 
\right)|det(\overline{t}-Id)|^{-1}
\end{equation}
Because $\psi(\xi_1)=1$ for $\|\xi_1\|\ge 1$ and the integral on the compact 
set $|\xi_1|\le 1$ yields a holomorphic function in $s$ we need to study the
meromorphic extension of 
\begin{equation}\label{trshom}
s\mapsto Tr\left(\int_{N\cap O}\int_{|\xi_1|\ge 1} \Tilde{b}_{s,i}(x_1,\xi_1)
\,\dbar\xi_1 dx_1\circ T \right)|det(\overline{t}-Id)|^{-1}
\end{equation}
Let $\xi_1=\lambda\overline{\xi}$ with $\lambda=\|\xi_1\|$ and $\overline{\xi}=
\dfrac {\xi_1}{\|\xi_1\|}
\in S^{n-1}$ be the decomposition in polar coordinates. The degree of homogeneity 
of $\Tilde{b}_{s,i}$ in $\xi_1$ is equal to $sd-i$. Then $\Tilde{b}_{s,i}(x_1,
\lambda\overline{\xi})= \lambda^{sd-i}\Tilde{b}_{s,i}(x_1,\overline{\xi})$ and 
after passing to polar coordinates, the expression \eqref{trshom} becomes:
\begin{align}
&Tr\left(\int_{N\cap O}\int_{S^{n-1}}\int_{1}^{\infty} \lambda^{sd-i}
\Tilde{b}_{s,i}(x_1,\overline{\xi})\lambda^{n-1}\,d\lambda \dbar\overline{\xi} 
dx_1\circ T \right) |det(\overline{t}-Id)|^{-1}=\\
&=Tr\left(\int_{N\cap O}\int_{S^{n-1}}\Tilde{b}_{s,i}(x_1,\overline{\xi}) 
\int_{1}^{\infty} \lambda^{sd+n-1-i}\,d\lambda \dbar\overline{\xi} 
dx_1\circ T \right)|det(\overline{t}-Id)|^{-1}=\\
&=-\frac{1}{sd+n-i}Tr\left(\int_{N\cap O}\int_{S^{n-1}}\Tilde{b}_{s,i}
(x_1,\overline{\xi})\, \dbar\overline{\xi} dx_1\circ T \right)|det(\overline{t}-Id)|^{-1}
\end{align}
The double integral defines a holomorphic function on $\mathbb{C}$, so the above 
expression has a meromorphic extension to $\mathbb{C}$ with a simple pole at 
$s=\frac{-n+i}{d}$ and residue: 
\begin{equation}\label{res}
-\frac{1}{d}Tr\left(\int_{N\cap O}\int_{S^{n-1}}\Tilde{b}_{s,i}(x_1,\overline{\xi})\, 
\dbar\overline{\xi} dx_1\circ T \right)|det(\overline{t}-Id)|^{-1}
\end{equation} 

As a consequence, the trace functions $Tr(A_s\circ \mathcal{T})$ and $Tr(A^s\circ
\mathcal{T})$ have meromorphic extensions on any half plane $Re(s)<K$ and so on 
$\mathbb{C}$. 
\par
We will proceed with the computation of the poles and residues
of $Tr(A_s\circ \mathcal{T})$. We observed that on the half space $Re(s)<K$ the 
difference 
\begin{equation}
Tr(A_s\circ \mathcal{T})-\sum_{i=0}^{\nu}Tr\left(\int_
{N\cap O} \int_{|\xi_1|\ge 1} \Tilde{b}_{s,i}(x_1,\xi_1)\,\dbar\xi_1 dx_1\circ  T 
\right)|det(\overline{t}-Id)|^{-1} 
\end{equation}
is holomorphic. Then $Tr(A_s\circ \mathcal{T})$ has 
simple poles at $s=\frac{-n+i}{d}$ for $i=0,1,\dots$ and the residue at $s=\frac
{-n+i}{d}$ is equal to the expression in \eqref{res}. Lemma \ref{reduce} gives an
asymptotic expansion of $\Tilde{b}_s(x_1,\xi_1)$ of the form:
\begin{equation}
\Tilde{b}_s(x_1,\xi_1)\sim\sum_{\alpha}\frac{1}{\alpha !}(D^{\alpha}_{w} \partial^
{\alpha}_{\xi_2}\overline{a}_s)(x_1,0,\xi_1,0)
\end{equation}
with $\overline{a}_s(x_1,w,\xi_1,\xi_2)=a_s(x_1,(\overline{t}-Id)^{-1}w,\xi_1,\xi_2)$
If $a_s\sim\sum a_{s,i}$ is the asymptotic expansion in homogeneous terms, with 
$a_{s,i}$ homogeneous of degree $sd-i$, then the homogeneous component of degree 
$sd-i$ of $\Tilde{b}_s$ will be equal to
\begin{equation}
\Tilde{b}_{s,i}(x_1,\xi_1)=\sum_{|\alpha|+k=i}\frac{1}{\alpha !}(D^{\alpha}_{w} \partial^
{\alpha}_{\xi_2}\overline{a}_{s,k})(x_1,0,\xi_1,0)
\end{equation}
with $\overline{a}_{s,k}(x_1,w,\xi_1,\xi_2)=a_{s,k}(x_1,(\overline{t}-Id)^{-1}w,\xi_1,\xi_2)$
of degree of homogeneity $sd-k$.

To conclude the proof of the theorem, consider the fixed point set $N=M^{\gamma}$ and a 
coordinate chart $O$ where $N\cap O= N_i\cap O$. Let $n_i=dim( N_i\cap O)$. We define 
the densities $\eta^{\gamma}_{i,k}$ that compute the residue of the zeta function 
at $s=\frac{-m+k}{d}$ as:
\begin{align}
\label{defn:eta}
\eta^{\gamma}_{i,k}(x_1)&=-\frac{1}{d}Tr(\int_{S^{n-1}}\Tilde{b}_{s,n_i-m+k}
(x_1,\overline{\xi})\,\dbar\overline{\xi}\circ T)\,dx_1&&\text{ if } k\ge m-n_i\\
\eta^{\gamma}_{i,k}&=0 &&\text{ if } k< m-n_i\\
\intertext{and}
\label{defn:d}
d_i^{\gamma}&=|det(\overline{t}-Id)|^{-1}&&
\end{align}
\end{proof}

We will now continue with the proof of Lemma \ref{reduce}. We will follow closely 
the ideas contained in the proof of Theorem 3.1 in \cite{shubin}.

\begin{proof}[Proof of Lemma \ref{reduce}]

We have $\|a_s(x,\xi)\| \le C(1+\|\xi\|)^{Re(s)d}$
so for $Re(s)<-\frac{m-n}{d}$ the integral defining $b_s$ is absolutely 
convergent so one can change the order of integration in \eqref{defn:b}.
Consider the Taylor expansion of $a_s(x_1,x_2,\xi_1,\xi_2)$ near $\xi_2=0$
\begin{gather}
a_s(x_1,x_2,\xi_1,\xi_2)=\sum_{|\alpha|\le Q-1}\frac{1}{\alpha!}(\partial^{\alpha}_
{\xi_2} a_s)(x_1,x_2,\xi_1,0)\xi^{\alpha}_2+R_Q \\
\intertext{with the reminder in the integral form}
r_Q=\sum_{|\alpha|=Q}\frac{Q\xi^{\alpha}_2}{\alpha !}\int_0^1(1-t)^{Q-1}(\partial
^{\alpha}_{\xi_2} a_s)(x_1,x_2,\xi_1,t\xi_2)\, dt
\end{gather}
We have
\begin{equation}
\iint e^{-i<x_2,\xi_2>} \frac{1}{\alpha!}(\partial_{\xi_2}^{\alpha}a_s)(x_1,
x_2,\xi_1,0)\xi_2^{\alpha}\,dx_2 \dbar \xi_2=
\frac{1}{\alpha!}(D^{\alpha}_{x_2}\partial^{\alpha}_{\xi_2} a_s)(x_1,0,\xi_1,0)
\end{equation}
the double integral being the composition of the Fourier transform in $x_2$
and the inverse Fourier transform in $\xi_2$ evaluated at $0$. Because
$a_s$ is a classical symbol of order $s$, the right-hand side term of 
the above equality is a classical symbol of order $s-|\alpha|$.

In order to prove that $b_s(x_1,\xi_1)$ is a symbol with the asymptotic expansion as in 
\eqref{defn:b} we need to show that the integral of the remainder 
$\iint e^{-i<x_2,\xi_2>} r_Q \,dx_2\dbar\xi_2$ is a symbol of arbitrary 
negative order if $Q$ is chosen large enough.
After changing the order of integration of $dt$ and $dx\dbar\xi_2$ we see that
it will be sufficient to provide a uniform estimate in $t\in(0,1]$ for the integrals 
\begin{equation}
R_{\alpha,t}(x_1,\xi_1)=\iint e^{-i<x_2,\xi_2>}\xi_2^{\alpha}(\partial ^{\alpha}_{\xi_2}
a_s)(x_1,x_2,\xi_1,t\xi_2)\,
dx_2\dbar\xi_2
\end{equation}
with $|\alpha|=Q$. Integrating by parts we get:
\begin{equation}
R_{\alpha,t}(x_1,\xi_1)=\iint e^{-i<x_2,\xi_2>} (D^{\alpha}_{x_2} 
\partial ^{\alpha}_{\xi_2}a_s) (x_1,x_2,\xi_1,t\xi_2)\,dx_2\dbar\xi_2
\end{equation}

Let $R_{\alpha,t}=R^1_{\alpha,t}+R^2_{\alpha,t}$ where 
$R^1_{\alpha,t}$ is the integral over the set
$D=$\linebreak
$\{(x_2,\xi_2)\,|\,\|\xi_2\|\le \|\xi_1\|\}$ and $R^2_{\alpha,t}$ 
is the integral over the complement. The volume of $D$ is 
bounded by $C\|\xi_1\|^n$ where $C$ doesn't depend on $t$, $x_1$ and $\xi_1$. 
On $D$ we have $\|(\xi_2,t\xi_2)\| \le 2\|\xi_1\|$ so the 
integrand in $R^1_{\alpha,t}$ is bounded by $C(1+\|\xi_1\|)^{Re(s)-Q}$. 
Thus
\begin{equation}
\|R^1_{\alpha,t}(x_1,\xi_1)\|\le C(1+\|\xi_1\|)^{Re(s)-Q+n}
\end{equation}
with $C$ independent of $t$, $x_1$ and $\xi_1$.
\par
By using the identity $(1+\|\xi_2\|^2)^{-\nu}(1+\|D_{x_2}\|^2)^{\nu}
e^{-i<x_2,\xi_2>}=e^{-i<x_2,\xi_2>}$ and integrating by parts, 
we can rewrite $R^2_{\alpha,t}$ as a sum of terms of the form:
\begin{equation}
\iint_{\|\xi_2\| > \|\xi_1\|} e^{-i<x_2,\xi_2>}(1+\|\xi_2\|^2)^{-\nu} 
(\partial^{\alpha}_{\xi_2} D^{\beta}_{x_2} a_s)(x_1,x_2,\xi_1,t\xi_2)
\,dx_2\dbar\xi_2
\end{equation}
with $|\beta|\le 2\nu$. Because $\|\xi_1\|<\|\xi_2\|$, the expression 
$(\partial^{\alpha}_{\xi_2} D^{\beta}_{x_2} a_s)(x_1,x_2,\xi_1,t\xi_2)$ is 
bounded from above by $C\|\xi_2\|^{Re(s)-Q}$ for $Re(s)\ge Q$ and by 
$C$ otherwise. The constant $C$ doesn't depend on $\beta$, $t$, $\xi_1$ 
and $\xi_2$. For $\nu$ large enough the previous integral is bounded by
\begin{equation}
C\int_{\|\xi_2\|>\|\xi_1\|}\!\!\!(1+\|\xi_2\|)^{-\nu'}\dbar\xi_2
\le C \|\xi_1\|^{-\nu'+n+1}\!\!\int (1+\|\xi_2\|)^{-n-1}\dbar \xi_2\le C 
\|\xi_1\|^{-\nu'+n+1}
\end{equation}
with $\nu'$ arbitrary large. Thus 
\begin{equation}
\|R_{\alpha,t}(x_1,\xi_1)\|\le C \|\xi_1\|^{Re(s)-Q+n}
\end{equation}
for an arbitrary $Q\in \NN$ with the constant $C$ independent of 
$(x_1,\xi_1)$. 
\par
Using Proposition $3.6$  and Theorem $3.1$ in \cite{shubin} we 
conclude that $b_s(x_1,\xi_1)$ is a symbol of order $s$ and has 
the asymptotic expansion given in 
\eqref{defn:b:asympt}.  If $a_s$ is a classical symbol then $b_s$
is classical as well.
\par
For the second part of the proof, let us consider the positive, real valued 
smooth functions $\psi_i(\xi_1)$, for $i\in\NN$, such that $\psi_i\equiv 0$ on
the ball of radius $i$ and $\psi_i\equiv 1$ outside the ball of radius $i+1$.
Let
\begin{equation}
\Tilde{b}_s(x_1,\xi_1)=\sum_{i=0}^{\infty}\psi_i(\xi_1)(\sum_{|\alpha|=i}
\frac{1}{\alpha !} (D^{\alpha}_{x_2} \partial^{\alpha}_{\xi_2}a_s)
(x_1,0,\xi_1,0))
\end{equation}
It is obvious that $\Tilde{b}_s-b_s$ is a smoothing operator.
If $a_s$ is a holomorphic family of symbols for $s\in \CC$ then $\Tilde{b}_s$
defined as above is a holomorphic family as well.
\end{proof}

\vfill\eject


\end{document}